\documentclass{article}
\usepackage{geometry}
\newgeometry{vmargin={20mm},hmargin={25mm,25mm}}
\usepackage[utf8]{inputenc}
\usepackage{mathtools}
\usepackage{amsmath}
\usepackage{amsthm}
\usepackage{amssymb}
\usepackage{amsfonts}
\usepackage{upgreek}
\usepackage{centernot}
\usepackage{stmaryrd}
\usepackage{comment}
\usepackage{enumerate}
\usepackage[shortlabels]{enumitem}
\usepackage[linguistics]{forest}
\usepackage{float}
\usepackage{wasysym}
\usepackage{caption}
\usepackage{subcaption}
\usepackage{tikz-qtree}
\usepackage{adjustbox}
\usepackage{graphicx}
\usepackage{epsfig}

\usepackage{accents}

\usepackage{xcolor}
\usepackage{hyperref}
\hypersetup{colorlinks,linkcolor={red},citecolor={blue},urlcolor={red},
  pdftitle={A Lagrangian Quantum Field Theory},pdfauthor={John G. Miller}}
\usepackage{xurl}

\theoremstyle{plain}
\newtheorem{theorem}{Theorem}[section]
\newtheorem{lemma}[theorem]{Lemma}
\theoremstyle{definition}
\newtheorem{defn}[theorem]{Definition}
\newtheorem{prop}[theorem]{Proposition}
\newtheorem{eg}[theorem]{Example}
\newtheorem{rmk}[theorem]{Remark}
\newtheorem{cor}[theorem]{Corollary}

\def\R{\mathbb{R}}
\def\C{\mathbb{C}}

\def\Z{\mathbb{Z}}
\def\N{\mathbb{N}}

\def\cc{\mathcal{C}}
\def\dd{\mathcal{D}}
\def\ff{\mathcal{F}}
\def\gg{\mathcal{G}}

\def\yy{\mathcal{Y}}

\def\xx{\mathcal{X}}

\def\hh{\mathcal{H}}

\def\mm{\mathcal{M}}
\def\aa{\mathcal{A}}

\def\ii{\mathcal{I}}
\usepackage{tikz-cd}
\tikzset{
  symbol/.style={
    draw=none,
    every to/.append style={
      edge node={node [sloped,allow upside down,auto=false]{$#1$}}
    }
  }
}
\usepackage{mathdots}
\usepackage{young}
\numberwithin{equation}{section}
\numberwithin{theorem}{subsection}

\usetikzlibrary{arrows.meta}
\newcommand{\kk}{\mathcal K}

\allowdisplaybreaks[2]
\author{John G. Miller}
\title{An Algebraic Framework for Lagrangian Quantum Field Theory}
\date{\today}
\begin{document}

\maketitle
\begin{abstract}
This paper provides a combinatorial, algebraic framework for understanding how specific $A_\infty$-functors of the form $\ff:Q_n\to\cc$ glue together, motivated by the gluing of Lagrangian cobordisms. The categories $Q_n$ generalise the universal exact triangle attributed to Kontsevich. The resulting structure is that of an $A_\infty$-multicategory, whose composition operations are governed by the combinatorial structure of equivalence classes of ordered tree diagrams. We detail this structure explicitly and use the resulting homotopy multicategory $\text{H}_0\text{Cob}(\cc)$ to produce results and invariants. Most notably, we realise the Grothendieck group of the derived category $\text{H}_0\text{Tw}\cc$ as an algebra over $\text{H}_0\text{Cob}(\cc)$. This can be viewed as an (algebraic) Lagrangian interpretation of a topological quantum field theory. We also introduce a family of Abelian groups $K^{Q;n}_0(\cc)$ which `filter' the Grothendieck group. These, together with the notion of constructability, allow us to detect information about $\cc$ that is not captured by $K_0(\text{H}_0\text{Tw}\cc)$.
\end{abstract}
\begingroup
\small
\tableofcontents
\endgroup

\section{Introduction}\label{sec:introduction}
Through a series of papers, Biran and Cornea (\cite{BC1}, \cite{BC2}, \cite{BC3}, \cite{BC4}) proved that Lagrangian cobordisms, specific types of Lagrangians in $M\times \C$ (where $M$ is our underlying symplectic manifold),  induce iterated cone decompositions in the derived Fukaya category $D\text{Fuk}(M)$. More recently, Bosshard proved that cobordisms in the non-compact setting induce iterated cone decompositions in the wrapped Fukaya category \cite{Bo}.
 The structure of Lagrangian cobordisms has also been studied by Nadler and Tanaka in terms of a stable infinity category \cite{NT}. The purpose of the following text is to build the correct combinatorial, algebraic framework to study Lagrangian cobordisms and their compositions. Given two Lagrangian cobordisms with a corresponding end, one can `compose' them by gluing them together along their common end. The result is a new cobordism; at the level of the derived Fukaya category, this corresponds to a concatenation of iterated mapping cones. One should have in mind the following picture

\vspace{1cm}
\tikzset{every picture/.style={line width=0.75pt}} 

\begin{tikzpicture}[x=0.75pt,y=0.75pt,yscale=-1,xscale=1]

\draw    (67,113) -- (160,113) ;
\draw    (66,72) -- (159,72) ;
\draw   (158,93) .. controls (158,81.95) and (173.67,73) .. (193,73) .. controls (212.33,73) and (228,81.95) .. (228,93) .. controls (228,104.05) and (212.33,113) .. (193,113) .. controls (173.67,113) and (158,104.05) .. (158,93) -- cycle ;
\draw    (168,78) .. controls (172,75) and (156,71) .. (159,72) ;
\draw    (171,109) .. controls (165,113) and (159,114) .. (160,113) ;
\draw    (228,91.88) -- (308,92) ;
\draw   (430,94) .. controls (430,82.95) and (445.67,74) .. (465,74) .. controls (484.33,74) and (500,82.95) .. (500,94) .. controls (500,105.05) and (484.33,114) .. (465,114) .. controls (445.67,114) and (430,105.05) .. (430,94) -- cycle ;
\draw    (500,94) -- (580,94.12) ;
\draw    (430,94) -- (342,93) ;
\draw    (342,144) .. controls (387,141) and (409,140) .. (449,113) ;
\draw [color={rgb, 255:red, 237; green, 14; blue, 14 }  ,draw opacity=1 ] [dash pattern={on 4.5pt off 4.5pt}]  (308,92) -- (342,93) ;

\draw (184,86.4) node [anchor=north west][inner sep=0.75pt]    {$V$};
\draw (453,85.4) node [anchor=north west][inner sep=0.75pt]    {$W$};
\draw (35,60.4) node [anchor=north west][inner sep=0.75pt]    {$A$};
\draw (35,102.4) node [anchor=north west][inner sep=0.75pt]    {$B$};
\draw (316,130.4) node [anchor=north west][inner sep=0.75pt]    {$D$};
\draw (318,58.4) node [anchor=north west][inner sep=0.75pt]    {$C$};
\draw (592,85.4) node [anchor=north west][inner sep=0.75pt]    {$E$};

\end{tikzpicture}

 \vspace{1cm}
 
The left cobordism (really its image under the projection $\pi:M\times \C\to \C$), denoted $V$, will induce an exact triangle in the derived Fukaya category (under suitable Floer-theoretic assumptions). Namely, the existence of $V$ implies that $
    C\cong \text{Cone}(A\xrightarrow{f}B)$
as objects in $D\text{Fuk}(M)$. Similarly, $W$, gives us another exact triangle, this time with $
    E\cong \text{Cone}(C\xrightarrow{g}D).$ 
 The morphisms $f$ and $g$, can be explicitly constructed from $V$ and $W$ respectively. This is done roughly as follows; consider $V$ and a specific curve $\gamma:(0,1)\to \C$, giving the following configuration

\vspace{1cm}
\begin{center}
\tikzset{every picture/.style={line width=0.75pt}} 

\begin{tikzpicture}[x=0.75pt,y=0.75pt,yscale=-1,xscale=1]

\draw    (172,109) -- (265,109) ;
\draw    (171,68) -- (264,68) ;
\draw   (263,89) .. controls (263,77.95) and (278.67,69) .. (298,69) .. controls (317.33,69) and (333,77.95) .. (333,89) .. controls (333,100.05) and (317.33,109) .. (298,109) .. controls (278.67,109) and (263,100.05) .. (263,89) -- cycle ;
\draw    (273,74) .. controls (277,71) and (261,67) .. (264,68) ;
\draw    (276,105) .. controls (270,109) and (264,110) .. (265,109) ;
\draw    (333,87.88) -- (413,88) ;
\draw    (201,13) -- (200,163) ;
\draw  [dash pattern={on 0.84pt off 2.51pt}]  (201,13) -- (201,-1) ;
\draw  [dash pattern={on 0.84pt off 2.51pt}]  (200,177) -- (200,163) ;
\draw [color={rgb, 255:red, 230; green, 13; blue, 13 }  ,draw opacity=1 ]   (202,68.5) .. controls (235.02,84.31) and (210.95,101.42) .. (202.59,106.95) ;
\draw [shift={(201,108)}, rotate = 326.31] [color={rgb, 255:red, 230; green, 13; blue, 13 }  ,draw opacity=1 ][line width=0.75]    (10.93,-3.29) .. controls (6.95,-1.4) and (3.31,-0.3) .. (0,0) .. controls (3.31,0.3) and (6.95,1.4) .. (10.93,3.29)   ;

\draw (289,82.4) node [anchor=north west][inner sep=0.75pt]    {$V$};
\draw (140,56.4) node [anchor=north west][inner sep=0.75pt]    {$A$};
\draw (140,98.4) node [anchor=north west][inner sep=0.75pt]    {$B$};
\draw (420,74.4) node [anchor=north west][inner sep=0.75pt]    {$C$};
\draw (210,7.4) node [anchor=north west][inner sep=0.75pt]    {$\gamma $};

\end{tikzpicture}
\end{center}
 \vspace{1cm}

The red arrow, represents (again the image under $\pi$ of ) a $J$-holomorphic strip, running between intersection points of \(V\cap(N\times\gamma)\) near the ends \(A\) and \(B\). Here, $N$ is an arbitrary `test Lagrangian'.  In particular, for each such $N$, we obtain a chain map $f^N:\text{CF}(N,A)\to \text{CF}(N,B)$. These maps $f^N$ can be  interpreted as a morphism of Yoneda modules, $\bar{f}:\yy(A)\to \yy(B)$, where $\yy(X):=\text{CF}(-,X)$. The Yoneda lemma allows us to then view $\bar{f}$ as $\bar{f}\simeq \yy(f)$ for some $f\in \text{CF}(A,B)$, this is our morphism defining the mapping cone.  When we `glue' the two cobordism $V$ and $W$ along their common end, in this example that being $C$ where the gluing is represented by the red dashed line connecting the two copies of $C$, we build a new cobordism $W\# V$, now with ends $A,B,D$ and $E$. This now represents the `iterated mapping cone' $E\cong \text{Cone}(A\to \text{Cone}(B\to D))$ (omitting grading/shifts) where every arrow can be viewed similarly to the above $f$. What is less well understood is what this operation is algebraically at the chain level. 
 
This operation, we will show, has a concrete description in terms of `gluings' of certain $A_\infty$ functors. Indeed, there are strictly unital $A_\infty$-categories which we denote $Q_n$, for any $n\in \N$, such that $A_\infty$-functors $\ff:Q_n \to \cc$ represent an iterated cone decomposition with `length' $n$ inside of the derived category $D\cc$. These categories were first introduced by Kontsevich (as attributed by Seidel in \cite{Se}), at least for the case of $n=2$, in which Seidel gives an explicit proof for the fact that functors $\ff:Q_2\to \cc$ represent exact triangles (again see \cite{Se}). Furthermore, Bosshard used the fact that these categories represent what they refer to as `universal cone decompositions' along with an identification of $Q_n$ with the full generating subcategory on `linking disks' in the partially wrapped Fukaya category of $\C$ with $n+1$ stops, to prove the generalisation of Biran and Cornea's result to the non-compact setting {\cite{Bo}}.

Slightly more concretely, each $Q_n$, has objects $\text{Obj}(Q_n)=\{0,\dots,n\}$, with hom-complexes generated by a unique closed morphism ${a_{i+1}:i\to i+1}$ between successive objects modulo $n+1$, as well as strict units. Every $\mu_k$ will vanish other than $\mu_{n+1}(a_{i+1},\dots,a_i)=1_i$ (and the obvious binary composition with a unit). One can realise that Lagrangian cobordisms have an almost identical structure to such a functor (though the Fukaya category need not be strictly unital), indeed an $A_\infty$-functor $\ff:Q_{n}\to \cc$ consists of a family of objects $\ff_0(i)$ and for every sequence of successive morphisms $(a_{i+1},\dots,a_{i+k})$ we obtain a morphism $\ff_k(a_{i+1},\dots,a_{i+k}):\ff_0(i)\to \ff_0(i+k)$ (again mod $n+1$). By generalising the construction of the morphism $f$ described before, we can construct corresponding morphisms induced by our cobordism. One should have in mind the following  

\vspace{1cm}
\begin{center}

\tikzset{every picture/.style={line width=0.75pt}} 

\begin{tikzpicture}[x=0.75pt,y=0.75pt,yscale=-1,xscale=1]

\draw    (243,118) -- (336,118) ;
\draw    (242,77) -- (335,77) ;
\draw   (334,98) .. controls (334,86.95) and (349.67,78) .. (369,78) .. controls (388.33,78) and (404,86.95) .. (404,98) .. controls (404,109.05) and (388.33,118) .. (369,118) .. controls (349.67,118) and (334,109.05) .. (334,98) -- cycle ;
\draw    (344,83) .. controls (348,80) and (332,76) .. (335,77) ;
\draw    (347,114) .. controls (341,118) and (335,119) .. (336,118) ;
\draw    (404,96.88) -- (484,97) ;
\draw   (281,94.25) .. controls (281,47.44) and (318.94,9.5) .. (365.75,9.5) .. controls (412.56,9.5) and (450.5,47.44) .. (450.5,94.25) .. controls (450.5,141.06) and (412.56,179) .. (365.75,179) .. controls (318.94,179) and (281,141.06) .. (281,94.25) -- cycle ;
\draw [color={rgb, 255:red, 241; green, 17; blue, 17 }  ,draw opacity=1 ]   (283.5,77.5) .. controls (298.23,89.38) and (293.09,102.6) .. (288.26,114.19) ;
\draw [shift={(287.5,116)}, rotate = 292.62] [color={rgb, 255:red, 241; green, 17; blue, 17 }  ,draw opacity=1 ][line width=0.75]    (10.93,-3.29) .. controls (6.95,-1.4) and (3.31,-0.3) .. (0,0) .. controls (3.31,0.3) and (6.95,1.4) .. (10.93,3.29)   ;
\draw [color={rgb, 255:red, 74; green, 144; blue, 226 }  ,draw opacity=1 ]   (283.5,77.5) .. controls (335.5,92) and (320.5,152.5) .. (360.5,157.5) .. controls (397.44,157.01) and (424.19,123.04) .. (448.87,98.13) ;
\draw [shift={(450,97)}, rotate = 135] [color={rgb, 255:red, 74; green, 144; blue, 226 }  ,draw opacity=1 ][line width=0.75]    (10.93,-3.29) .. controls (6.95,-1.4) and (3.31,-0.3) .. (0,0) .. controls (3.31,0.3) and (6.95,1.4) .. (10.93,3.29)   ;

\draw (360,91.4) node [anchor=north west][inner sep=0.75pt]    {$V$};
\draw (211,65.4) node [anchor=north west][inner sep=0.75pt]    {$A$};
\draw (211,107.4) node [anchor=north west][inner sep=0.75pt]    {$B$};
\draw (491,88.4) node [anchor=north west][inner sep=0.75pt]    {$C$};
\draw (436,14.4) node [anchor=north west][inner sep=0.75pt]    {$\gamma $};

\end{tikzpicture}
\end{center}
\vspace{1cm}
The red arrow gives us a closed degree zero morphism $f_1:A\to B$, the blue, a degree one morphism $f_2:A\to C$. The point is, if we construct all such arrows, that must run counterclockwise, and in general allow loops, we can associate the following diagram (on the left)

\begin{equation*}
\tikzset{
    every loop/.append style={min distance=15mm,looseness=10}
}
\begin{tikzcd}[row sep=2em,column sep=2em]
    A
    \ar[loop above]
    \ar[dd]
    \ar[drr,bend left]
    & & \\
    & &
    C
    \ar[loop right]
    \ar[ull]
    \ar[dll,bend left]
    \\
    B
    \ar[loop below]
    \ar[uu,bend left]
    \ar[urr]
    & &
\end{tikzcd}
\qquad
\begin{tikzcd}[row sep=2em,column sep=2em]
    \ff_0(0)
    \ar[loop above,"\ff_3(a_1{,}a_2{,}a_0)"]
    \ar[dd,"\ff_1(a_1)"]
    \ar[drr,bend left,"\ff_2(a_1{,}a_2)"]
    & & \\
    & &
    \ff_0(2)
    \ar[loop right,"\ff_3(a_0{,}a_1{,}a_2)"]
    \ar[ull,"\ff_1(a_0)"]
    \ar[dll,bend left,"\ff_2(a_0{,}a_1)"]
    \\
    \ff_0(1)
    \ar[loop below,"\ff_3(a_2{,}a_0{,}a_1)"]
    \ar[uu,bend left,"\ff_2(a_2{,}a_0)"]
    \ar[urr,"\ff_1(a_2)"]
    & &
\end{tikzcd}
\end{equation*}

The diagram on the right depicts the image of some $\ff:Q_2 \to \cc$ (with $\cc$ representing the Fukaya category of $M$, both diagrams only show up to order three morphisms). Hopefully the similarity is clear. The construction in \cite{BC4}
gives related morphisms and homotopies at the level of
$A_\infty$-modules. The resulting exact triangle can be encoded as
a functor $\ff:Q_2\to\cc$, using the universal exact triangle
described in \cite{Se}.

Thus we ask the following: given a strictly unital $A_\infty$-category $\cc$ and two strictly unital $A_\infty$-functors $\ff:Q_n \to \cc$ and $\gg:Q_m \to \cc$ such that there exist some $i\in Q_n$ and $j\in Q_m$ with $\ff_0(i)=\gg_0(j)$, can we construct a new functor $\ff\#_{ij}\gg:Q_{n+m-1}\to \cc$ which corresponds to a concatenation of iterated cone decompositions in the derived category of $\cc$?

One should realise that there is a subtle difference between an iterated cone decomposition in $D\cc$ and a functor $\ff:Q_n \to \cc$. Indeed, such functors can still exist even if $\cc$ is not pre-triangulated; moreover, as we will discuss later, a functor $\ff:Q_n \to \cc$ need not decompose as a sequence of gluings of functors $\gg:Q_2\to \cc$ (i.e., exact triangles). In fact, in the symplectic setting, with $\cc=\text{Fuk}(M)$ for some $M$, this appears to be much closer to the embedded geometry of Lagrangian cobordisms. From this we are able to build a ($A_\infty$) multicategory (coloured operad) of such functors, which we denote $\text{Cob}(\cc)$.

A related operadic structure appears in a remark in \cite{Ta}, which describes planar coloured operads associated with cobordisms and cone decompositions, and relates them to the Waldhausen $S$-construction studied by Dyckerhoff and Kapranov \cite{DK}. In the stable setting, these operations are given by filtrations, with composition given by refining them. Here, we give an explicit combinatorial construction at the $A_\infty$-level, using functors $Q_n\to\cc$ as multimorphisms and labelled trees to describe their gluings and coherent higher homotopies. We do not assume that $\cc$ is pre-triangulated, and our multimorphisms need not decompose into triangles within $\cc$.

This also gives a useful way to view the Grothendieck group: as an algebra over the homotopy multicategory $\text{H}_0\text{Cob}(\cc)$, with each multimorphism encoding the relation that its output class is the sum of its input classes. With this, we can also define a family of groups $K_0^{Q;n}(\cc)$, for $n\geq0$, by imposing only the relations coming from functors $Q_m\to\cc$ with $m\leq n$, together with the shift relations. These form a sequence of surjective homomorphisms whose colimit is $K_0(\text{H}_0\text{Tw}\cc)$. The intermediate groups record which relations are generated by operations with at most $n$ inputs, and can retain information about $\cc$ that is lost in passing to its pre-triangulated completion.

We begin the paper by discussing some background material, explicitly defining the categories $Q_n$ and exploring some of their basic properties. Note that for all following statements we work over $k=\Z/2\Z$. In Section~\ref{sec:functors}, we describe the process for building a glued functor $\ff\#_{ij}\gg:Q_{n+m-1}\to \cc$. This involves two steps: first, the construction of an $A_\infty$-functor $\ff\vee_{ij}\gg:Q_n\bigvee_{ij}Q_m \to \cc$, with $Q_n\bigvee_{ij}Q_m$ denoting the pushout over the inclusion of the point category $*$ under the inclusions as object $i$ and object $j$ of $Q_n$ and $Q_m$, respectively. Second, we construct a functor $\Phi_{ij}:Q_{n+m-1}\to Q_n\bigvee_{ij}Q_m$; the resulting gluing will be defined by
\begin{equation}
     \ff\#_{ij}\gg:=\Phi_{ij}^*(\ff\vee_{ij} \gg):Q_{n+m-1}\to \cc.
    \end{equation}
We show that this result is the `correct' functor; for example, if $n=m=1$, the construction of $\ff\#_{ij}\gg:Q_1\to \cc$ can be interpreted as the homotopy equivalences being closed under composition. 
When we then consider sequences of such gluings, we quickly realise the structure is somewhat complicated. To begin, we consider a sequence of `parallel gluings', i.e., many $\ff^i$-s glued to a single $\gg$. We show the following:
\begin{theorem}
    Let $\gg:Q_m\to \cc$ and consider a family of functors $\ff^I:=\{\ff^{i}:Q_{n_i}\to \cc\}_{i\in I\subset \{0,\dots,m-1\}}$ with $\ff^i_0(n_i)=\gg_0(i)$. Then there exists a family of functors 
    \begin{equation}\Big\{\#^p(\ff^I;\gg):Q_{\sum_{i\in I}n_i +m-|I|}\to \cc\Big\}_{p\in \text{Ord}(I)}\end{equation}
   which, with corresponding higher homotopies, satisfy the cubical homotopy relations. 
\end{theorem}
Here, $\text{Ord}(I)$ denotes the set of all ordered sequences $(i_1,\dots,i_{|I|})$ with $i_j\in I$ appearing exactly once. 
One can unpack this statement as the existence of a family of functors (and higher homotopies) $\Gamma^I$ and, for every such $I$, a family of functors (and higher homotopies) $\Gamma^{A|B}$ with $A,B\subset I$, with $A\cup B=I$ and $A\cap B=\emptyset$, such that
\begin{equation}
    d^{\text{Fun}}\Gamma^{I}=\sum_{\substack{A\cup B=I,\\ A\cap B=\emptyset,\\ A,B\neq\emptyset}} \Gamma^{A|B}{+C_\Gamma^I}
\end{equation}
where $\Gamma^{A|B}$ corresponds to a composition of $\Gamma^B$ with $\Gamma^A$, {and $C_\Gamma^I$ contains the additional composition terms, as in~\eqref{eq:cubicalarity}}; in particular, $\Gamma^{\{i\}|\{j\}}$ corresponds to $(\ff^i\#\gg)\#_{}\ff^j$ and $\Gamma^{\{j\}|\{i\}}$ to $\ff^i\#(\gg\#_{}\ff^j)$ and $\Gamma^{\{i,j\}}$ is a $1$-homotopy between them.

Following this, we then consider what we refer to as `series gluings', which correspond to a sequence of functors $(\ff^i)_{i=0,\dots,k}$ being glued sequentially, i.e., with each functor attached to the next. In particular, we have $\ff^i_0(n_i)=\ff^{i+1}_0(j)$ for some $j\neq n_{i+1}$. In this setting, we prove the following:
\begin{theorem}\label{intro:Main theroem for series}
   For every $k\geq 2$, the family of functors \[
\#^v(\ff^0,\ldots,\ff^k):
Q_{\sum_{i=0}^k n_i-k}\longrightarrow\cc,
\qquad v\in V(k),
\]
admits a coherent family of $l$-homotopies, for $1\leq l \leq k-1$, that correspond to the fillings of the cells of the Stasheff associahedra $K(k+1)$ (see \cite{St1,St2}).
\end{theorem}
Here, $V(k)$ denotes the set of vertices of $K(k+1)$. In explicit terms, this amounts to the existence of a family of functors (and higher homotopies) $\Theta^{\sigma}$ in which $\sigma$ denotes a generic `bracketing' of the list $0,1,\dots,k$. In particular, for every such $k$ we show that there is $\Theta^{01...k}$, corresponding to the empty bracketing, for which
\begin{equation*}
    {\mu_1^{\text{Fun}}\Theta^{0\dots k}=\sum_{\ell=2}^{k}\sum_{r=0}^{k+1-\ell} \Theta^{0\dots (r\,\dots\,r+\ell-1)\dots k}+C^{0\dots k}}
\end{equation*}
{Here, $C^{0\dots k}$ contains the composition terms described in Section~\ref{subsec:series-gluings}.} In this setup, one would find $\Theta^{(01)2}$ corresponds to $(\ff^0\#\ff^1)\#\ff^2$, $\Theta^{0(12)}$ to $\ff^0\#(\ff^1\#\ff^2)$ and $\Theta^{012}$ a homotopy between them.

The combinatorics of both parallel and series gluings can be reinterpreted in terms of `ordered trees'. Here, trees are fixed planar trees in which every edge is given a corresponding natural number, which will represent the order in which we attach the branch. Multiple branches can carry the same number, which corresponds to a homotopy with degree given by the multiplicities. This construction is very similar to that of Boardman--Vogt \cite{BV}, but a discretised version, whereas they allow the edges to be labeled continuously.
The common idea is to use trees to organise compositions and their coherent higher homotopies, with Berger and Moerdijk extending the Boardman--Vogt construction to coloured operads \cite{BM}. We show that if one considers a general sequence of gluings of functors, i.e., a mix of both series and parallel, that the corresponding combinatorial structure of the homotopies is that of these trees. Explicitly:
\begin{theorem}\label{mainmaintheorem}
   Given a collection of $\ff^{\mathbf{\alpha}}:Q_{n_\mathbf{\alpha}}\to \cc$ with $\mathbf{\alpha}\in \mathcal{I}$ for some fixed indexing set $\mathcal{I}$, there is a correspondence between all orders of gluings of functors, corresponding higher homotopies, and the set of equivalence classes of ordered trees, $\overline{\text{Tree}(\ii)}$.  More precisely, to each tree class $\tau$
of degree zero we assign the functor
\[
    \Theta^\tau:=\#^\tau(\ff^{\mathcal I}):
    Q_{d(\mathcal I)}\longrightarrow\cc,
\]
and to each tree class $\tau$ of degree $q>0$ we assign a
$q$-dimensional homotopy $\Theta^\tau$. These assignments respect
the equivalence and zero relations and, for $q>0$, satisfy
\[
    \mu_1^{\mathrm{Fun}}\Theta^\tau
    =\Theta^{\mu\tau}{+C^\tau}.
\]

\end{theorem}

This structure allows us to consider an `$A_\infty$-multicategory' associated with any strictly unital $A_\infty$-category $\cc$, that we denote by $\text{Cob}(\cc)$, with objects (colours) given by the objects of $\cc$ and multimorphisms (with $n$ inputs) given by functors $\ff:Q_n\to \cc$. The structure of the compositions of such functors is entirely controlled by that of Theorem \ref{mainmaintheorem}.

In Section~\ref{sec:multicategory} we then study the resulting structure of the homotopy multicategory $\text{H}_0\text{Cob}(\cc)$ given by taking homotopy classes of multimorphisms in $\text{Cob}(\cc)$. The key construction is that of an Abelian group 
\begin{equation}
        K^Q_0(\cc):=\frac{\Z\cdot \text{Obj}(\cc)}{\sim}
    \end{equation}
    where $\sim$ denotes the collection of relations 
    \begin{equation}
        \chi(\ff):=[\ff_0(0)]+\sum_{i=2}^{n}[\ff_0(i)]-[\ff_0(1)] \qquad [A[s]]=(-1)^s[A]
    \end{equation}
    for any $\ff:Q_n\to \cc$, for any $n\geq 0$. The signs here are chosen to correspond with the chosen grading of $Q_n$ which we fix in Section~\ref{subsec:cyclic-quivers}.
We take $\chi(\ff)=[\ff_0(0)]-[\ff_0(1)]$ for $\ff:Q_1\to\cc$, and for $n=0$ we take $\chi(\ff)=[\ff_0(0)]$. We assume that $\cc$ has shifts, and impose the shift relation $[A[s]]=(-1)^s[A]$. With this, we prove that this group is in fact isomorphic to the usual Grothendieck group of the derived category of $\cc$, namely, there is an isomorphism of groups \begin{equation}
        \lambda:K^Q_0(\cc) \xrightarrow{\cong} K_0(\text{H}_0\text{Tw}\cc)
    \end{equation}
    with $\lambda([A])=[A]$.
Here, we take the derived category of $\cc$ to be $\text{H}_0\text{Tw}\cc$. Thus one recovers its standard Grothendieck group from $\text{H}_0\text{Cob}(\cc)$.
 This gives rise to an interesting idea: 
\begin{theorem}\label{intro:maintheoremalg}
    There is an algebra $\aa:\text{H}_0\text{Cob}(\cc)\to \dd$ which realises the Grothendieck group of $\text{H}_0\text{Tw}\cc$.
\end{theorem}
Here, $\dd$ is a category that we will construct specifically for this purpose in Section~\ref{subsec:cob-category}.
The algebra realises the group in the sense that each relation $\chi(\ff)$ corresponds to $\aa(\ff)$.
We can view the standard Grothendieck group as an algebra over the multicategory $\text{H}_0\text{Cob}(\cc)$. From this there are many interesting questions one can study; in particular, what other such algebras exist? This is briefly discussed in Section~\ref{sec:questions}. Along this line of study we show that the Grothendieck group, which we now view as $K^Q_0(\cc)$, admits a filtration by Abelian groups. Here, by a `filtration' we mean a sequence of surjective homomorphisms between Abelian groups whose colimit is $K^Q_0(\cc)$. In particular, we define
 
 \begin{equation}
        K_0^{Q;n}(\cc):=\frac{\Z\cdot \text{Obj}(\cc)}{\sim_n}
    \end{equation}
    where $\sim_n$ denotes the collection of relations 
    \begin{equation}
        \chi(\ff):=[\ff_0(0)]+\sum_{i=2}^{m}[\ff_0(i)]-[\ff_0(1)] \qquad [A[s]]=(-1)^s[A]
    \end{equation}
    for any $\ff:Q_m\to \cc$, for any $ 0\leq m\leq n$. With this, one finds that there exist surjective homomorphisms
    \begin{equation}
        \lambda_{n,n+1}:K_0^{Q;n}(\cc)\to K_0^{Q;n+1}(\cc)
    \end{equation}
and one can recover 
    \begin{equation}
        K^Q_0(\cc) \cong \text{colim}_{n\geq 0}\left(\to K^{Q;n}_0(\cc) \xrightarrow{\lambda_{n,n+1}}K^{Q;n+1}_0(\cc)\to \dots \right).
    \end{equation}
The remainder of Section~\ref{sec:multicategory} is dedicated to exploring the relation between these groups and a notion we refer to as (weak) $n$-constructability; $\cc$ being $n$-constructable if every multimorphism in $\text{Cob}(\cc)$ (or $\text{H}_0\text{Cob}(\cc)$) decomposes as a sequence of compositions of multimorphisms $\gg^i:Q_{k_i}\to \cc$ with $k_i\leq n$. We show that if $\cc$ is weakly $n$-constructable, then \begin{equation}
    \kappa(\cc)\leq n,
\end{equation}
where $\kappa(\cc)$ denotes the least $m$ such that $\ker(\lambda_m:K^{Q;m}_0(\cc)\to K^Q_0(\cc))=0$.
If no such finite $m$ exists, we set $\kappa(\cc)=\infty$. This invariant, together with the groups \(K^{Q;n}_0(\mathcal C)\), can distinguish \(\mathcal C\) from its pre-triangulated completion \(\operatorname{Tw}\mathcal C\). Indeed, we show that for \(\mathcal C=Q_n\), with \(n>2\), we have
\[
\kappa(Q_n)=n>2,
\qquad
\kappa(\operatorname{Tw}Q_n)=2.
\]
Thus these invariants can distinguish categories that have equivalent derived categories.

Finally, in Section~\ref{sec:questions} we comment on how some of the results relate to other works and discuss some interesting possible future directions of study.

\par\smallskip
{\centering\textbf{Acknowledgements}\par}
\smallskip

Parts of the contents of this paper appeared (in a preliminary version) as a chapter of my doctoral thesis whilst at Universit\'e de Montr\'eal; this can be found in \cite{MiTh}.
 During this time I was supervised by Octav Cornea, whom I would like to thank for their many discussions on this topic and for getting me involved in this area of study. I would also like to acknowledge and thank my co-supervisor Egor Shelukhin for their guidance throughout my PhD, as well as Paul Biran for their useful discussions during my time at ETH. 

I would also like to remark that this text was proofread with GPT. From the resulting feedback, I made some minor corrections, most notably to indexing in equations for $A_\infty$-relations and filled some small gaps in proofs, along with some general formatting. The most notable gap appeared in my original proof of Lemma~\ref{lemphihomology}. Though the proof given in this paper follows the same idea, it has been corrected. The overall theory and content are my own.

\subsection{Cyclic Quiver Categories}\label{subsec:cyclic-quivers}

\textbf{A note on the setting:} For the contents of this paper, we choose to work with $\Z$-graded, strictly unital, $A_\infty$-categories over the field $k=\Z/2\Z$. 

There exists a family of strictly unital $A_\infty$-categories $Q_n$ for $n\in \N$ (where we believe $0\in \N$). The objects are given by the collection of natural numbers less than or equal to $n$:
\begin{equation}\text{Obj}(Q_n):=\{0,\dots,n\} .\end{equation}
For $n\geq1$, the morphism complexes are given by 
\begin{equation}
    Q_n(i,j):=\begin{cases}k \cdot \langle a_{i+1}\rangle & j=i+1 \mod n+1\\
    k\cdot \langle 1_i\rangle & i=j\\
    0 & \text{else}
    \end{cases},
\end{equation}
where $1_i$ is the (strict) unit for the object $i\in Q_n$, and (following a slight adjustment to the convention of Bosshard \cite{Bo})
\begin{equation}
    |a_{i+1}|=\begin{cases}
        0 & i= 0 \text{ or }1,\\-1 & \text{otherwise}.
    \end{cases}
\end{equation}
i.e., the maps out of $0\in Q_n$ and $1\in Q_n$ are of trivial degree, and all other maps have degree $-1$. The $A_\infty$-structure relations on $Q_n$ are such that
\begin{equation}
    \mu_k(a_{i+1},\dots,a_{i+k})=\begin{cases}
        0 & k\neq n+1\\
        1_i & k=n+1
    \end{cases} ,
\end{equation}
where indices are taken modulo $n+1$.
For $n=0$, the category $Q_0$ has a single object $0$, with $Q_0(0,0)=k\langle1_0,a_0\rangle$, where $1_0$ is the degree-zero unit and $a_0:0\to0$ is a non-identity morphism of degree $1$ satisfying $\mu_1(a_0)=1_0$.
 Note that in our conventions we have
\begin{equation*}
    \mu_k:\cc_{n_1}(A_0,A_1)\otimes \dots \otimes \cc_{n_k}(A_{k-1},A_k)\to \cc_{\sum_{i=1}^{k}n_i+k-2}(A_0,A_k)
\end{equation*}
in particular, for $n\geq1$,
\begin{equation*}
|\mu_{n+1}(a_{i+1},\dots,a_i)|=\sum_{j=0}^{n}|a_j|+(n+1)-2=-(n-1)+(n-1)=0.
\end{equation*}

\begin{eg}
    For example, $Q_2$ can be viewed as the diagram 
\begin{equation*}
    \begin{tikzcd}
        & 1\ar[dr,"a_2"] & \\
        0\ar[ur,"a_1"] & &2\ar[ll,"a_0"]
    \end{tikzcd}
\end{equation*}
where $|a_1|=|a_2|=0$ and $|a_0|=-1$. If we view $Q_n$ inside its shift completion, we can regard the diagram above as the sequence
\begin{equation}
    \begin{tikzcd}
        0\ar[r,"a_1"] & 1\ar[r,"a_2"] & 2\ar[r,"{a_0}"] & {[1]0}
    \end{tikzcd}
\end{equation}
where $\mu_1(a_{i+1})=0$, $\mu_2(a_{i+1},a_{i+2})=0$ and $\mu_3(a_{i+1},a_{i+2},a_i)=1_i$, which is a realisation of the `universal exact triangle'.
 
For $n=3$, we may think of $Q_3$ as 
\begin{equation*}
    \begin{tikzcd}
        0\ar[r,"a_1"] & 1\ar[d,"a_2"]\\ 3\ar[u,"a_0"]& 2\ar[l,"a_3"]
    \end{tikzcd}
\end{equation*}
which unravels to a sequence:
\begin{equation}\label{seqq3}
    \begin{tikzcd}
        0\ar[r,"a_1"] & 1\ar[r,"a_2"] & 2\ar[r,"a_3"] & {[1]3}\ar[r,"{[1]a_0}"] & {[2]0}.
    \end{tikzcd}
\end{equation}
In general, one should have in mind the $(n+1)$-gon:
\begin{equation}
\begin{tikzcd}
  & & 0 \ar[drr,"a_1"]& &\\
  n\ar[urr,"a_0"]& & && 1 \ar[d,"a_2"]\\
  n-1\ar[u,"a_n"]& & & &2\ar[dl,"a_3"]\\
  &\dots \ar[ul,"a_{n-1}"]& & \ar[ll,dashed]\dots & \\
\end{tikzcd} 
\end{equation}
and its unravelling to 
\begin{equation*}
    \begin{tikzcd}
        0\ar[r,"a_1"] & 1\ar[r,"a_2"] & 2\ar[r,"a_3"] & {[1]3}\ar[r,"{[1]a_4}"] &\dots \ar[r,"{[n-3]a_n}"] & {[n-2]n} \ar[r,"{[n-2]a_0}"]&{[n-1] 0}
    \end{tikzcd}
\end{equation*}
\end{eg}

\begin{rmk}
    Since the categories are cyclic, it is worth noting that one also finds rotations of the above sequences; for example, rotating (to the left) the sequence for $Q_3$ in equation \ref{seqq3} gives 
    \begin{equation*}
         \begin{tikzcd}
       1\ar[r,"a_2"] & 2\ar[r,"a_3"] & {[1]3}\ar[r,"{[1]a_0}"] & {[2]0} \ar[r,"{[2]a_1}"]& {[2]1}.
    \end{tikzcd}
    \end{equation*}
\end{rmk}

The purpose of these categories is to describe (iterated) mapping cones in $A_\infty$-categories. As mentioned in the introduction, this notion can be attributed to Kontsevich, as done by Seidel in \cite{Se}.
 To understand this, consider a (strictly unital) functor $\ff:Q_n \to \cc$. It is defined via a family of objects $\ff_0(i)\in \cc$ indexed over $\{0,\dots,n\}$, and morphisms $\ff_1(a_{i+1}):\ff_0(i)\to \ff_0(i+1)$ of degree $|\ff_1(a_{i+1})|=|a_{i+1}|$ as well as higher homotopies $\ff_k(a_{i+1},\dots,a_{i+k}):\ff_0(i)\to \ff_0(i+k)$ such that, for $k\geq1$,
\begin{align*}
    &\sum_{\substack{1\leq l\leq k\\ r_1+\cdots+r_l=k\\ r_1,\dots,r_l\geq1}} \mu_l(\ff_{r_1}(a_{i+1},\dots,a_{i+r_1}),\dots,\ff_{r_l}(a_{i+k+1-r_l},\dots,a_{i+k}))\\
    =&\sum_{s=1}^{k}\sum_{j=0}^{k-s} \ff_{k-s+1}(a_{i+1},\dots,a_{i+j},\mu_s(a_{i+j+1},\dots,a_{i+j+s}),a_{i+j+s+1},\dots,a_{i+k})\\
    =&\begin{cases}
        \ff_1(1_i) & k=n+1\\
        0 & k\neq n+1
    \end{cases} .
\end{align*}
Here, object and arrow indices are taken modulo $n+1$.

\begin{rmk}
    The degrees of these morphisms are given by
    \begin{equation}
        |\ff_k(a_{i+1},\dots, a_{i+k})|=(k-1) +\sum_{t=1}^{k}|a_{i+t}|
    \end{equation}

\end{rmk}

    If we consider the case of $n=1$, we can view $Q_1$ as the diagram 
    
    \begin{equation*}
        \begin{tikzcd}
            0\ar[r,bend left, "a_1"] & 1\ar[l,bend left,"a_0"]
        \end{tikzcd}
    \end{equation*}
    where $|a_0|=|a_1|=0$, $\mu_2(a_1,a_0)=1_0$ and $\mu_2(a_0,a_1)=1_1$.  Thus, a functor $\ff:Q_1 \to \cc$ defines a diagram 
\begin{equation}
    \begin{tikzcd}
        \ff_0(0)\ar[r,bend left,"\ff_1(a_1)"]\ar[loop left,"{\ff_2(a_1,a_0)}"] & \ff_0(1)\ar[loop right,"{\ff_2(a_0,a_1)}"]\ar[l,bend left,"\ff_1(a_0)"]
    \end{tikzcd}
\end{equation}
(along with some higher order terms), where $|\ff_1(a_1)|=0=|\ff_1(a_0)|$ and $|\ff_2(a_1,a_0)|=|\ff_2(a_0,a_1)|=1$ such that 
\begin{equation}
    \mu_1\ff_2(a_1,a_0)+\mu_2(\ff_1(a_1),\ff_1(a_0))= \ff_1(1_0)=1_{\ff_0(0)} .
\end{equation}
That is, $\ff_2(a_1,a_0)$ defines a homotopy $\ff_2(a_1,a_0):\mu_2(\ff_1(a_1),\ff_1(a_0)) \implies 1_{\ff_0(0)}$. In particular, we see that $\ff_1(a_1)$ defines a homotopy equivalence $\ff_0(0) \simeq \ff_0(1)$ with homotopy inverse $\ff_1(a_0)$. Thus, functors $\ff:Q_1 \to \cc$ characterise homotopy equivalences in $\cc$.

\begin{eg}
    If we consider $n=2$, we find that such a functor $\ff:Q_2 \to \cc$ defines a diagram 
    \begin{equation}
    \begin{tikzcd}
        & \ff_0(0)\ar[loop,swap,"\ff_3(a_1{,}a_2{,}a_0)"]\ar[dr,"\ff_1(a_1)"]\ar[dl,bend right=70,swap,"\ff_2(a_1{,}a_2)"]&\\
        \ff_0(2)\ar[loop left,"\ff_3(a_0{,}a_1{,}a_2)"]\ar[ur,"\ff_1(a_0)"] \ar[rr,bend right =60,swap,"\ff_2(a_0{,}a_1)"]& & \ff_0(1)\ar[loop right,"\ff_3(a_2{,}a_0{,}a_1)"]\ar[ul,bend right =70,swap,"\ff_2(a_2{,}a_0)"]\ar[ll,"\ff_1(a_2)"]
    \end{tikzcd}
\end{equation}

\begin{align*}
    \mu_1(\ff_3(a_0,a_1,a_2))+\mu_2(\ff_1(a_0),\ff_2(a_1,a_2))\quad&\\
    \quad+\mu_2(\ff_2(a_0,a_1),\ff_1(a_2))+\mu_3(\ff_1(a_0),\ff_1(a_1),\ff_1(a_2))=&1_{\ff_0(2)}\\
    \mu_1(\ff_2(a_0,a_1))+\mu_2(\ff_1(a_0),\ff_1(a_1))=&0\\
    \mu_1(\ff_2(a_1,a_2))+\mu_2(\ff_1(a_1),\ff_1(a_2))=&0.
    \end{align*}
These state exactly that $\ff_0(2)$ is the mapping cone of the morphism $\ff_1(a_1)$. Seidel \cite{Se} shows that these relations are sufficient to prove that the triangle
\begin{equation}
    \ff_0(0) \xrightarrow{\ff_1(a_1)}\ff_0(1) \xrightarrow{\ff_1(a_2)}\ff_0(2) \xrightarrow{\ff_1(a_0)}[1]\ff_0(0)
\end{equation}
is exact in $\cc$. Note, furthermore, that the cyclic nature of $Q_2$ proves that all rotations of this are also exact.
 
\end{eg}

\begin{rmk}
    Here, `exact' means that the image of this triangle under $\iota:\cc\to\text{Tw}\cc$ is homotopy equivalent to a standard mapping-cone triangle. In particular, omitting $\iota$ from the notation, we have
    \[
        \ff_0(2)\simeq\ff_0(1)\oplus\ff_0(0)[1],
    \]
    where the right-hand side carries the mapping cone twisting determined by $\ff_1(a_1)$.
\end{rmk}

\begin{eg}
    We also consider the case of $Q_0$. This is the category with one object, $0$, and $Q_0(0,0)$ generated by the unit $1_0$ and a single non-identity morphism $a_0$ of degree $|a_0|=1$, such that $\mu_1(a_0)=1_0$; in particular, $0$ is contractible with contraction given by $a_0$. A functor from this category, $\ff:Q_0 \to \cc$, defines an object $\ff_0(0)$ and a contraction $\ff_1(a_0):\ff_0(0)\to [-1]\ff_0(0)$.

\end{eg}
One viewpoint to have in mind is to think of a functor $\ff:Q_2 \to \cc$ as a hole with a filling inside of $\cc$. The generalisations, i.e., functors $\ff:Q_n \to \cc$, provide generalisations of exact triangles and can be thought of as higher dimensional holes inside of $\cc$. Such functors can in fact be seen as iterated cones in the derived category: 
\begin{align*}
    \ff:Q_0 \to \cc \hspace{0.5cm} & \hspace{0.5cm} \text{Contractible objects }\ff_0(0), \\
    \ff:Q_1 \to \cc \hspace{0.5cm} & \hspace{0.5cm} \text{Homotopy equivalences } \ff_0(0) \simeq \ff_0(1),\\
    \ff:Q_2 \to \cc \hspace{0.5cm} & \hspace{0.5cm} \text{Exact triangles } \ff_0(0) \to \ff_0(1)\to \ff_0(2)\to [1]\ff_0(0),\\
    \vdots\\
    \ff:Q_n \to \cc \hspace{0.5cm} & \hspace{0.5cm} \text{Generalised cone decompositions } \ff_0(0) \to \hdots \to [n-2]\ff_0(n) \to [n-1]\ff_0(0).\\
    \vdots
\end{align*}

\begin{rmk}
    A key point is that, even when a morphism $f\in\cc(A,B)$ is closed and of degree zero, it need not extend to a functor $\ff:Q_n\to\cc$ with $\ff_1(a_1)=f$. We do not assume that $\cc$ is pre-triangulated or satisfies any analogous completion property.

\end{rmk}

\subsection{Multicategories}
For references on multicategories, see \cite{Le}, \cite{LV}, \cite{BM} and \cite{Tr}. We use the terminology of multicategories, rather than the perhaps more common terminology of coloured operads. Throughout, our multicategories are non-symmetric, with ordered input lists.

A \textbf{multicategory}, $\mm$, consists of a collection of objects $\text{Obj}(\mm)$ (referred to as colours in the language of coloured operads) and a family of multimorphisms. That is, for every collection $A_0,\dots,A_n,B\in\text{Obj}(\mm)$, with $n\geq0$, we have a collection of operations called multimorphisms. We denote this collection by $\mm(A_0,\dots,A_n:B)$, and call an element $f\in\mm(A_0,\dots,A_n:B)$ an $(n+1)$-ary multimorphism with inputs $A_0,\dots,A_n$ and output $B$.
\begin{rmk}
    We also allow empty input lists. For every $B\in\mm$, the collection $\mm(;B)$ consists of nullary multimorphisms, with no inputs and output $B$. In a composition, inserting a nullary multimorphism at an input $B_i$ replaces that input by the empty list. A nullary multimorphism itself has no input slots at which to compose.
\end{rmk}
One can view a multimorphism $f\in \mathcal{M}(A_0,\dots,A_n:B)$ as either of the following diagrams

\begin{center}
    
\tikzset{every picture/.style={line width=0.75pt}} 

\begin{tikzpicture}[x=0.75pt,y=0.75pt,yscale=-1,xscale=1]

\draw   (320,73.88) -- (280.2,109.11) -- (279.8,39.12) -- cycle ;
\draw    (188,93) -- (281,93) ;
\draw    (186,73) -- (279,73) ;
\draw    (187,52) -- (280,52) ;
\draw    (320,73.88) -- (400,74) ;
\draw  [dash pattern={on 0.84pt off 2.51pt}]  (205,52) -- (205,73) ;
\draw  [dash pattern={on 0.84pt off 2.51pt}]  (205,73) -- (205,94) ;
\draw    (178,225) -- (271,225) ;
\draw    (176,205) -- (269,205) ;
\draw    (177,184) -- (270,184) ;
\draw  [dash pattern={on 0.84pt off 2.51pt}]  (195,184) -- (195,205) ;
\draw  [dash pattern={on 0.84pt off 2.51pt}]  (195,205) -- (195,226) ;
\draw   (269,205) .. controls (269,193.95) and (284.67,185) .. (304,185) .. controls (323.33,185) and (339,193.95) .. (339,205) .. controls (339,216.05) and (323.33,225) .. (304,225) .. controls (284.67,225) and (269,216.05) .. (269,205) -- cycle ;
\draw    (279,190) .. controls (283,187) and (267,183) .. (270,184) ;
\draw    (282,221) .. controls (276,225) and (270,226) .. (271,225) ;
\draw    (339,203.88) -- (419,204) ;

\draw (138.15,34.4) node [anchor=north west][inner sep=0.75pt]  [xslant=-0.04]  {$ \begin{array}{l}
A_{0}\\
\end{array}$};
\draw (149,59.4) node [anchor=north west][inner sep=0.75pt]    {$A_{i}$};
\draw (147,83.4) node [anchor=north west][inner sep=0.75pt]    {$A_{n}$};
\draw (404,65.4) node [anchor=north west][inner sep=0.75pt]    {$B$};
\draw (288,65.4) node [anchor=north west][inner sep=0.75pt]    {$f$};
\draw (136,166.4) node [anchor=north west][inner sep=0.75pt]  [xslant=-0.04]  {$ \begin{array}{l}
A_{0}\\
\end{array}$};
\draw (139,191.4) node [anchor=north west][inner sep=0.75pt]    {$A_{i}$};
\draw (137,215.4) node [anchor=north west][inner sep=0.75pt]    {$A_{n}$};
\draw (423,195.4) node [anchor=north west][inner sep=0.75pt]    {$B$};
\draw (295,197.4) node [anchor=north west][inner sep=0.75pt]    {$f$};

\end{tikzpicture}
\end{center}
where the first gives the more classical depiction of a multimorphism (as in \cite{Le}), the second gives a depiction which resembles the `spider' diagrams used by Cornea and Biran to represent Lagrangian cobordisms. 

Our multimorphisms come equipped with a family of composition operations. For nonempty input lists, with $I\subseteq\{0,\dots,n\}$ and $n_i\geq0$ for each $i\in I$, we write
\begin{equation}
\begin{aligned}
    \circ_I:\quad &\mm(B_0,\dots,B_n:C)\times
    \prod_{i\in I}\mm(A^i_0,\dots,A^i_{n_i}:B_i)\\
    &\longrightarrow\mm(\mathbf{L}_0,\dots,\mathbf{L}_n:C),
\end{aligned}
\end{equation}
where the factors indexed by $I$ are ordered by increasing $i$, and
\[
    \mathbf{L}_i=
    \begin{cases}
        (A^i_0,\dots,A^i_{n_i}) & i\in I,\\
        (B_i) & i\notin I.
    \end{cases}
\]
The inputs in the target are obtained by concatenating the lists $\mathbf{L}_0,\dots,\mathbf{L}_n$ in that order. In other words, we replace each chosen input $B_i$ by the inputs of the multimorphism attached there. For a nullary multimorphism at input $B_i$, the corresponding factor is $\mm(;B_i)$ and $\mathbf{L}_i$ is the empty list.
We can depict this as:

\tikzset{every picture/.style={line width=0.75pt}} 

\begin{tikzpicture}[x=0.75pt,y=0.75pt,yscale=-1,xscale=1]

\draw    (87,113) -- (180,113) ;
\draw    (85,93) -- (178,93) ;
\draw    (86,72) -- (179,72) ;
\draw  [dash pattern={on 0.84pt off 2.51pt}]  (104,72) -- (104,93) ;
\draw  [dash pattern={on 0.84pt off 2.51pt}]  (104,93) -- (104,114) ;
\draw   (178,93) .. controls (178,81.95) and (193.67,73) .. (213,73) .. controls (232.33,73) and (248,81.95) .. (248,93) .. controls (248,104.05) and (232.33,113) .. (213,113) .. controls (193.67,113) and (178,104.05) .. (178,93) -- cycle ;
\draw    (188,78) .. controls (192,75) and (176,71) .. (179,72) ;
\draw    (191,109) .. controls (185,113) and (179,114) .. (180,113) ;
\draw    (248,91.88) -- (328,92) ;
\draw   (432,160) .. controls (432,148.95) and (447.67,140) .. (467,140) .. controls (486.33,140) and (502,148.95) .. (502,160) .. controls (502,171.05) and (486.33,180) .. (467,180) .. controls (447.67,180) and (432,171.05) .. (432,160) -- cycle ;
\draw    (502,160) -- (582,160.12) ;
\draw    (85,181) -- (178,181) ;
\draw    (83,161) -- (176,161) ;
\draw    (84,140) -- (177,140) ;
\draw  [dash pattern={on 0.84pt off 2.51pt}]  (102,140) -- (102,161) ;
\draw  [dash pattern={on 0.84pt off 2.51pt}]  (102,161) -- (102,182) ;
\draw   (176,161) .. controls (176,149.95) and (191.67,141) .. (211,141) .. controls (230.33,141) and (246,149.95) .. (246,161) .. controls (246,172.05) and (230.33,181) .. (211,181) .. controls (191.67,181) and (176,172.05) .. (176,161) -- cycle ;
\draw    (186,146) .. controls (190,143) and (174,139) .. (177,140) ;
\draw    (189,177) .. controls (183,181) and (177,182) .. (178,181) ;
\draw    (246,159.88) -- (326,160) ;
\draw    (84,250) -- (177,250) ;
\draw    (82,230) -- (175,230) ;
\draw    (83,209) -- (176,209) ;
\draw  [dash pattern={on 0.84pt off 2.51pt}]  (101,209) -- (101,230) ;
\draw  [dash pattern={on 0.84pt off 2.51pt}]  (101,230) -- (101,251) ;
\draw   (175,230) .. controls (175,218.95) and (190.67,210) .. (210,210) .. controls (229.33,210) and (245,218.95) .. (245,230) .. controls (245,241.05) and (229.33,250) .. (210,250) .. controls (190.67,250) and (175,241.05) .. (175,230) -- cycle ;
\draw    (185,215) .. controls (189,212) and (173,208) .. (176,209) ;
\draw    (188,246) .. controls (182,250) and (176,251) .. (177,250) ;
\draw    (245,228.88) -- (325,229) ;
\draw    (328,92) .. controls (377,91) and (453,135) .. (450,141) ;
\draw    (326,160) -- (432,160) ;
\draw    (325,229) .. controls (378,227) and (412,207) .. (452,177) ;
\draw  [dash pattern={on 0.84pt off 2.51pt}]  (288,91.94) -- (286,159.94) ;
\draw  [dash pattern={on 0.84pt off 2.51pt}]  (287,160.94) -- (285,228.94) ;
\draw    (329,28) .. controls (393,28) and (477,122) .. (467,140) ;
\draw    (467,180) .. controls (464,216) and (376,282) .. (336,282) ;
\draw  [dash pattern={on 0.84pt off 2.51pt}]  (322,28.94) -- (320,91) ;
\draw  [dash pattern={on 0.84pt off 2.51pt}]  (325,229) -- (323,279) ;

\draw (328,64.4) node [anchor=north west][inner sep=0.75pt]    {$B_{i_{0}}$};
\draw (205,84.4) node [anchor=north west][inner sep=0.75pt]    {$f^{i_{0}}$};
\draw (202,153.4) node [anchor=north west][inner sep=0.75pt]    {$f^{i_{l}}$};
\draw (201,218.4) node [anchor=north west][inner sep=0.75pt]    {$f^{i_{k}}$};
\draw (43,78.4) node [anchor=north west][inner sep=0.75pt]    {$A_{j_{0}}^{i_{0}}$};
\draw (328,134.4) node [anchor=north west][inner sep=0.75pt]    {$B_{i_{l}}$};
\draw (326,204.4) node [anchor=north west][inner sep=0.75pt]    {$B_{i_{k}}$};
\draw (43,148.4) node [anchor=north west][inner sep=0.75pt]    {$A_{j_{l}}^{i_{l}}$};
\draw (43,213.4) node [anchor=north west][inner sep=0.75pt]    {$A_{j_{k}}^{i_{k}}$};
\draw (461,154.4) node [anchor=north west][inner sep=0.75pt]    {$g$};
\draw (587,149.4) node [anchor=north west][inner sep=0.75pt]    {$C$};
\draw (325,6.4) node [anchor=north west][inner sep=0.75pt]    {$B_{0}$};
\draw (331,260.4) node [anchor=north west][inner sep=0.75pt]    {$B_{n}$};

\end{tikzpicture}

For every object $A\in\mm$ (from now on, we write $A\in\mm$ instead of $A\in\text{Obj}(\mm)$), we require an identity $1$-ary morphism $1_A\in\mm(A:A)$. Inserting an identity at an input leaves that input unchanged. More explicitly, for $g\in\mm(B_0,\dots,B_n:C)$ and $J\subseteq I$, with $f_i=1_{B_i}$ for all $i\in J$, we require
\[
    \circ_I\bigl(g;(f_i)_{i\in I}\bigr)
    =
    \circ_{I\setminus J}\bigl(g;(f_i)_{i\in I\setminus J}\bigr),
\]
where $\circ_{\emptyset}(g)=g$.
We also require the unit condition at the output:
\[
    \circ_{\{0\}}(1_B;f)=f
    \qquad\text{for }f\in\mm(A_0,\dots,A_n:B).
\]
See \cite{Le}.
Finally, we require the composition operations to be associative: the order in which we carry out the compositions along a fixed tree does not matter. See the references above for this property in explicit terms.

Much like operads, one can study representations of (algebras over) multicategories. A \textbf{$\cc$-algebra over $\mm$}, where $\cc$ is a symmetric monoidal category, consists of an assignment $\aa:\text{Obj}(\mm)\to\text{Obj}(\cc)$. That is, to every object $A\in\mm$ we associate an object $\aa(A)\in\cc$. We also associate with every multimorphism $f\in\mm(A_0,\dots,A_n:B)$ a morphism in $\cc$,
\begin{equation}
    \aa(f):\aa(A_0)\otimes_\cc\dots\otimes_\cc\aa(A_n)\longrightarrow\aa(B),
\end{equation}
where $\otimes_\cc$ is the product defining the monoidal structure of $\cc$. We require that $\aa$ commutes with the composition operations on $\mm$ and preserves identities, so that $\aa(1_A)=1_{\aa(A)}$. For a nullary multimorphism, the source of $\aa(f)$ is the monoidal unit of $\cc$.
\begin{rmk}
    A non-symmetric operad is simply a non-symmetric multicategory with one object, and an algebra over such an operad is an algebra over a multicategory with one object.
\end{rmk}

The rest of the paper will be in the setting of $A_\infty$-categories, and the structure we construct will be an \textbf{$A_\infty$-multicategory}. Its composition operations are associative up to homotopy, with coherent higher homotopies, rather than associative on the nose.

More specifically, in this paper we view an $A_\infty$-multicategory as the structure described by the trees in Theorem~\ref{mainmaintheorem}. Trees of degree zero describe compositions, while trees of positive degree describe the corresponding higher homotopies. These satisfy the boundary relations given by the tree differential, together with the equivalence and zero relations of that construction. The formal definitions of these trees and their equivalence and zero relations will be given in Sections~\ref{subsec:ordered-trees}--\ref{subsec:all-gluings}.

In this setting, we need to consider two types of composition. We can compose in \textbf{parallel}, attaching many multimorphisms to one multimorphism (think $(f_0,\dots,f_k)\to g$), or in \textbf{series}, attaching multimorphisms in a chain (think $f_0\to f_1\to\dots\to f_k\to g$). We can also combine the two.

Our construction has many similarities with the $W$-construction of Boardman and Vogt \cite{BV} and its extension to coloured operads by Berger and Moerdijk \cite{BM}. The coloured viewpoint is particularly natural here, since the objects of our multicategory serve as colours. In our interpretation, the edges carry discrete labels recording the stages at which branches are attached. Equal labels represent simultaneous attachments and the corresponding higher homotopies. The tree differential separates a simultaneous stage into two, and we impose the equivalence and zero relations described later. In the classical topological $W$-construction, the internal edges instead carry continuous length parameters in $[0,1]$.
Suppose we have a multimorphism at the base which we wish to compose with a family of other multimorphisms. For example, take multimorphisms labelled $i,j,k$, where $i$ and $j$ can be attached to the base, and $k$ can be attached to $i$. We can view this in terms of trees:

\begin{equation*}
    \begin{tikzcd}
k& & & &  \\
        & i\ar[ul,no head, "0"]& & j&\\
        & & \ar[ul,no head,"0"]\cdot\ar[ur,no head,swap,"0"] & &
    \end{tikzcd}
\end{equation*}

This tree should be viewed as attaching all branches simultaneously, at time $t=0$. We can then consider the six trees that correspond to a single attaching followed by a pair attached simultaneously, or vice versa:
\begin{equation*}
\begin{gathered}
\begin{tikzcd}[ampersand replacement=\&,row sep=small,column sep=small]
k \& \& \& \& \\
\& i\ar[ul,no head,"0"] \& \& j \& \\
\& \& \cdot\ar[ul,no head,"0"]\ar[ur,no head,swap,"1"] \& \&
\end{tikzcd}
\qquad
\begin{tikzcd}[ampersand replacement=\&,row sep=small,column sep=small]
k \& \& \& \& \\
\& i\ar[ul,no head,"0"] \& \& j \& \\
\& \& \cdot\ar[ul,no head,"1"]\ar[ur,no head,swap,"0"] \& \&
\end{tikzcd}
\qquad
\begin{tikzcd}[ampersand replacement=\&,row sep=small,column sep=small]
k \& \& \& \& \\
\& i\ar[ul,no head,"1"] \& \& j \& \\
\& \& \cdot\ar[ul,no head,"0"]\ar[ur,no head,swap,"0"] \& \&
\end{tikzcd}
\\[1.5em]
\begin{tikzcd}[ampersand replacement=\&,row sep=small,column sep=small]
k \& \& \& \& \\
\& i\ar[ul,no head,"1"] \& \& j \& \\
\& \& \cdot\ar[ul,no head,"1"]\ar[ur,no head,swap,"0"] \& \&
\end{tikzcd}
\qquad
\begin{tikzcd}[ampersand replacement=\&,row sep=small,column sep=small]
k \& \& \& \& \\
\& i\ar[ul,no head,"1"] \& \& j \& \\
\& \& \cdot\ar[ul,no head,"0"]\ar[ur,no head,swap,"1"] \& \&
\end{tikzcd}
\qquad
\begin{tikzcd}[ampersand replacement=\&,row sep=small,column sep=small]
k \& \& \& \& \\
\& i\ar[ul,no head,"0"] \& \& j \& \\
\& \& \cdot\ar[ul,no head,"1"]\ar[ur,no head,swap,"1"] \& \&
\end{tikzcd}
\end{gathered}
\end{equation*}
 The original tree (with all attachments simultaneous at $t=0$) can be viewed as a $2$-homotopy, whose boundary is the sum of the resulting trees with attachments at two distinct times. Each of these trees can in turn be viewed as a $1$-homotopy, whose boundary gives the two orders of attaching the branches that share a time. The relations of an $A_\infty$-multicategory, in the sense used here, are governed by these tree boundaries.
\begin{rmk}
    In the above depiction of the six trees, the tree  
    \begin{equation*}
       \begin{tikzcd}[ampersand replacement=\&,row sep=small,column sep=small]
k \& \& \& \& \\
\& i\ar[ul,no head,"0"] \& \& j \& \\
\& \& \cdot\ar[ul,no head,"1"]\ar[ur,no head,swap,"0"] \& \&
\end{tikzcd}
    \end{equation*}
    will actually correspond to a trivial homotopy. Indeed, attaching $\cdot -j$ and $i-k$ occur at the same time, and happen away from one another, this will be made explicit when we impose an equivalence relation of the set of ordered trees.
\end{rmk}
\begin{rmk}
    With these discrete labels, we view attaching a single node, along a branch with time $t\in\N$, as a single composition. When we attach two nodes at the same time, we view this as `attaching a $1$-homotopy', attaching 3 nodes as `attaching a $2$-homotopy' and so on.
\end{rmk}

\section{The Gluing}\label{sec:functors}

The goal is to create purely combinatorial machinery that tells us how to glue two functors $\ff:Q_n \to \cc$ and $\gg:Q_m \to \cc$, with $n,m>0$, together along a common end (assume for now $\ff_0(i)=\gg_0(j)$) to form a new functor $\ff\#_{ij}\gg:Q_{n+m-1}\to \cc$. The cases $n=0$ or $m=0$ will be handled separately.

\begin{rmk}
    One can think of this as saying that we can concatenate two holes of dimensions $n$ and $m$ in $\cc$ which meet along a vertex (object) to form a new hole of dimension $n+m-1$.
\end{rmk}

\subsection{The pushout functor}
This will work as follows: first consider the pushout over the inclusions of the point category $*$ into $Q_n$ (as object $i$) and into $Q_m$ (as object $j$):
\begin{equation}
    \begin{tikzcd}
        & Q_n\ar[dr,dashed]& \\
        *\ar[ur,"i"]\ar[dr,swap,"j"]& & Q_n \bigvee_{ij} Q_m\\
        & Q_m\ar[ur,dashed]&
    \end{tikzcd}
\end{equation}
Here, $*$ denotes the strictly unital $k$-linear $A_\infty$-category with one object, also denoted $*$, and morphism space $\operatorname{Hom}_*(*,*)=k\cdot1_*$ concentrated in degree zero.
Its operations are $\mu_2(\lambda1_*,\eta1_*)=(\lambda\eta)1_*$ for $\lambda,\eta\in k$, and $\mu_d=0$ for $d\ne2$.

Explicitly, we have
\begin{align*}
    i_0(*)=&i\in Q_n\\
    i_1(1_*)=&1_i\in Q_n\\
    i_k(1_*,\dots,1_*)=&0 \in Q_n \quad \forall k>1.
\end{align*}
and similarly for $j:* \to Q_m$. One should note that $i,j$ are both `strict' strictly unital functors. It is known that such a pushout is well defined in the category of strict strictly unital functors, i.e., we obtain a new pushout category whose universal property applies to strict functors; see {\cite{Or} Theorem~4.5}. Ideally, we would like a diagram
\begin{equation}
    \begin{tikzcd}
        & Q_n\ar[dr]\ar[rrrd,bend left,"\ff"]& \\
        *\ar[ur,"i"]\ar[dr,swap,"j"]& & Q_n \bigvee_{ij} Q_m \ar[rr,dashed,"\ff\vee_{ij}\gg"] & & \cc\\
        & Q_m\ar[ur]\ar[rrru,bend right,swap,"\gg"]&
    \end{tikzcd}
\end{equation}
where $Q_n\bigvee_{ij}Q_m$ denotes the pushout in the category of strictly unital $A_\infty$-categories and strict strictly unital functors. This construction gives us the category. However, the given functors $\ff$ and $\gg$ need not be strict, so the strict-pushout universal property does not itself provide an extension $\ff\vee_{ij}\gg:Q_n\bigvee_{ij}Q_m\to\cc$. The remainder of this section is devoted to constructing such a strictly unital $A_\infty$-functor explicitly: we specify its components so that it restricts to $\ff$ and $\gg$ on the two factors, and verify the $A_\infty$-functor identities.

We should first take a look at the structure of $Q_n \bigvee_{ij} Q_m$. Writing $P:=Q_n\bigvee_{ij}Q_m$ for this strict pushout, we have
\begin{align*}
\operatorname{Obj}(P)&=\bigl(\operatorname{Obj}(Q_n)\sqcup\operatorname{Obj}(Q_m)\bigr)/(i\sim j),\\
|\operatorname{Obj}(P)|&=n+m+1.
\end{align*}

For example, one should picture $Q_2 \vee_{20} Q_2$ as
\begin{equation}
    \begin{tikzcd}
        &1 \ar[dr]& & 1'\ar[dr]\\
        0 \ar[ur]& &2=0'\ar[ll]\ar[ur] && 2'\ar[ll]
    \end{tikzcd}
\end{equation}

For $n,m\geq1$, the non-identity morphisms in the two factors entering the common object $i=j$ are generated by $a_i:i-1\to i$ in $Q_n$ and $b_j:j-1\to j$ in $Q_m$. Indices are taken modulo $n+1$ in the first factor and modulo $m+1$ in the second.

Assume $x\in P$ lives inside the copy of $Q_n$ and $y\in P$ lives inside the copy of $Q_m$. Then $P(x,y)$ is the $k$-linear span of formal iterated operations on composable morphisms from the two factors, with source $x$ and target $y$, modulo the $A_\infty$-relations, the relations in each factor and strict unitality. Thus, we keep track not only of the underlying path, but also of the nested $\mu$-operations applied to it.

For example, in $Q_2\vee_{20}Q_2$, there is a path given in red:
\begin{equation}
    \begin{tikzcd}
        &1 \ar[dr,red]& & 1'\ar[dr,red]\\
        0 \ar[ur,red]& &2=0'\ar[ll]\ar[ur,red] && 2'\ar[ll]
    \end{tikzcd}
\end{equation}
which gives us an element
\begin{equation*}
  \mu_4(a_1,a_2,b_1,b_2)\in  Q_2\bigvee_{20} Q_2(0,2')
\end{equation*}
with degree given by $|\mu_4(a_1,a_2,b_1,b_2)|=4-2+|a_1|+|a_2|+|b_1|+|b_2|=2$.
The $A_\infty$-relations are `freely taken' but restrict to those of $Q_n$ and $Q_m$, so we have
\begin{align*}
\mu_1\mu_4(a_1,a_2,b_1,b_2)
&=\mu_3(\mu_2(a_1,a_2),b_1,b_2)+\mu_3(a_1,\mu_2(a_2,b_1),b_2)\\
&\quad+\mu_3(a_1,a_2,\mu_2(b_1,b_2))
+\mu_2(\mu_3(a_1,a_2,b_1),b_2)\\
&\quad+\mu_2(a_1,\mu_3(a_2,b_1,b_2))\\
&=\mu_3(a_1,\mu_2(a_2,b_1),b_2)\\
&\quad+\mu_2(\mu_3(a_1,a_2,b_1),b_2)
+\mu_2(a_1,\mu_3(a_2,b_1,b_2)).
\end{align*}

All the above terms define elements of $Q_2\bigvee_{20} Q_2(0,2')$.
Furthermore, we can look at the morphism $\mu_3(a_1,a_2,b_1)$ and the relations give
\begin{align*}
    \mu_1\mu_3(a_1,a_2,b_1)=&\mu_2(\mu_2(a_1,a_2),b_1)+\mu_2(a_1,\mu_2(a_2,b_1))\\
    =&\mu_2(a_1,\mu_2(a_2,b_1)).
\end{align*}
So how do we construct our desired map $\ff\vee_{ij}\gg:Q_n\bigvee_{ij}Q_m\to\cc$? Of course, we require that for any $k\geq1$ and any $\mathbf{a}=(a_{l+1},\dots,a_{l+k})$, a sequence of composable morphisms belonging entirely to $Q_n$, we have
\begin{equation}
    (\ff\vee_{ij}\gg)_k(\mathbf{a})=\ff_k(\mathbf{a})
\end{equation}
and similarly, for any $\mathbf{b}=(b_{l+1},\dots,b_{l+k})$ of composable morphisms entirely in $Q_m$, we have
\begin{equation}
    (\ff\vee_{ij}\gg)_k(\mathbf{b})=\gg_k(\mathbf{b})
\end{equation}
Furthermore, in order for this $\ff\vee_{ij}\gg$ to really represent the desired gluing, we want
\begin{equation}
(\ff\vee_{ij}\gg)_1(\mu_2(a_{i},b_{j+1})):=\mu_2(\ff_1(a_i),\gg_1(b_{j+1})).
\end{equation}
Note that the $A_\infty$-functor relation gives
\begin{align*}
&(\ff\vee_{ij}\gg)_2(\mu_1(a_i),b_{j+1})+(\ff\vee_{ij}\gg)_2(a_i,\mu_1(b_{j+1}))+(\ff\vee_{ij}\gg)_1(\mu_2(a_{i},b_{j+1}))\\
&=\mu_2(\ff_1(a_i),\gg_1(b_{j+1}))+\mu_1(\ff\vee_{ij}\gg)_2(a_i,b_{j+1}).\\[0.5ex]
\implies &(\ff\vee_{ij}\gg)_1(\mu_2(a_{i},b_{j+1}))
=\mu_2(\ff_1(a_i),\gg_1(b_{j+1}))+\mu_1(\ff\vee_{ij}\gg)_2(a_i,b_{j+1})
\end{align*}
Thus, we can take
\begin{equation}
    (\ff\vee_{ij}\gg)_2(a_i,b_{j+1}):=0
\end{equation}

\begin{rmk}Taking this component to be zero is a choice, rather than a requirement of the functor relation; as we will see, this choice is well suited to the construction.\end{rmk}

So let's look at how we can extend from here. Consider the sequence $(a_{i-1},a_i,b_{j+1})$ (here we assume $n\geq2$). We require that
\begin{align*}
&\mu_1(\ff\vee_{ij}\gg)_3(a_{i-1},a_i,b_{j+1})+\mu_2((\ff\vee_{ij}\gg)_2(a_{i-1},a_i),(\ff\vee_{ij}\gg)_1(b_{j+1}))\\&+\mu_2((\ff\vee_{ij}\gg)_1(a_{i-1}),(\ff\vee_{ij}\gg)_2(a_i,b_{j+1}))\\
&+\mu_3((\ff\vee_{ij}\gg)_1(a_{i-1}),(\ff\vee_{ij}\gg)_1(a_i),(\ff\vee_{ij}\gg)_1(b_{j+1}))\\
=&(\ff\vee_{ij}\gg)_1(\mu_3(a_{i-1},a_i,b_{j+1}))+(\ff\vee_{ij}\gg)_2(\mu_2(a_{i-1},a_i),b_{j+1})+(\ff\vee_{ij}\gg)_2(a_{i-1},\mu_2(a_i,b_{j+1}))
\end{align*}
Here, every $a_k$ and $b_k$ is a cycle. Plugging in the definitions
 of $(\ff\vee_{ij}\gg)_1$ and $(\ff\vee_{ij}\gg)_2$  then gives
\begin{align*}
    &\mu_1(\ff\vee_{ij}\gg)_3(a_{i-1},a_i,b_{j+1})+\mu_2(\ff_2(a_{i-1},a_i),\gg_1(b_{j+1}))
    +\mu_3(\ff_1(a_{i-1}),\ff_1(a_i),\gg_1(b_{j+1}))\\
    =&(\ff\vee_{ij}\gg)_1(\mu_3(a_{i-1},a_i,b_{j+1}))+(\ff\vee_{ij}\gg)_2(a_{i-1},\mu_2(a_i,b_{j+1}))
\end{align*}
We can now set
\begin{align}(\ff\vee_{ij}\gg)_3(a_{i-1},a_i,b_{j+1})&:=0\\
(\ff\vee_{ij}\gg)_1(\mu_3(a_{i-1},a_i,b_{j+1}))&:=\mu_3(\ff_1(a_{i-1}),\ff_1(a_i),\gg_1(b_{j+1}))\label{eq:pushout-first-three}\\
(\ff\vee_{ij}\gg)_2(a_{i-1},\mu_2(a_i,b_{j+1}))&:=\mu_2(\ff_2(a_{i-1},a_i),\gg_1(b_{j+1}))\label{eq:pushout-mixed-two}
\end{align}
Once these choices have been made, the functor relation forces the value in \eqref{eq:pushout-mixed-two} on the pair $(a_{i-1},\mu_2(a_i,b_{j+1}))$. We check its compatibility with the functor relation and determine the remaining value on the nested composite.
Notice that
\begin{align*}
\mu_1(\ff\vee_{ij}\gg)_2(a_{i-1},\mu_2(a_i,b_{j+1}))
&=\mu_1\mu_2(\ff_2(a_{i-1},a_i),\gg_1(b_{j+1}))\\
&=\mu_2(\mu_1\ff_2(a_{i-1},a_i),\gg_1(b_{j+1}))\\
&=\mu_2(\mu_2(\ff_1(a_{i-1}),\ff_1(a_i)),\gg_1(b_{j+1}))\\
&=\mu_1\mu_3(\ff_1(a_{i-1}),\ff_1(a_i),\gg_1(b_{j+1}))\\
&\quad+\mu_2(\ff_1(a_{i-1}),\mu_2(\ff_1(a_i),\gg_1(b_{j+1})))\\
&=\mu_1(\ff\vee_{ij}\gg)_1(\mu_3(a_{i-1},a_i,b_{j+1}))\\
&\quad+\mu_2((\ff\vee_{ij}\gg)_1(a_{i-1}),(\ff\vee_{ij}\gg)_1(\mu_2(a_i,b_{j+1}))).
\end{align*}
The value of $(\ff\vee_{ij}\gg)_1(\mu_3(a_{i-1},a_i,b_{j+1}))$ was fixed in \eqref{eq:pushout-first-three}; this forces the value on the nested composite. We know that
\begin{align*}
    \mu_1(\ff\vee_{ij}\gg)_1(\mu_3(a_{i-1},a_i,b_{j+1}))=&(\ff\vee_{ij}\gg)_1(\mu_1\mu_3(a_{i-1},a_i,b_{j+1}))\\
    =&(\ff\vee_{ij}\gg)_1(\mu_2(a_{i-1},\mu_2(a_i,b_{j+1})))
\end{align*}
i.e.,
\begin{align*}
     &\mu_1 (\ff\vee_{ij}\gg)_2(a_{i-1},\mu_2(a_i,b_{j+1}))+\mu_2((\ff\vee_{ij}\gg)_1(a_{i-1}),(\ff\vee_{ij}\gg)_1(\mu_2(a_i,b_{j+1})))\\=&(\ff\vee_{ij}\gg)_1(\mu_2(a_{i-1},\mu_2(a_i,b_{j+1})))
\end{align*}
that is, $(\ff\vee_{ij}\gg)_2(a_{i-1},\mu_2(a_i,b_{j+1})):=\mu_2(\ff_2(a_{i-1},a_i),\gg_1(b_{j+1}))$ is well defined. So, to sum up, so far we have
\begin{align*}
     (\ff\vee_{ij}\gg)_k(\mathbf{a})&:=\ff_k(\mathbf{a})\\
     (\ff\vee_{ij}\gg)_k(\mathbf{b})&:=\gg_k(\mathbf{b})\\
     (\ff\vee_{ij}\gg)_1(\mu_2(a_{i},b_{j+1}))&:=\mu_2(\ff_1(a_i),\gg_1(b_{j+1})),\\
     (\ff\vee_{ij}\gg)_1(\mu_3(a_{i-1},a_i,b_{j+1}))&:=\mu_3(\ff_1(a_{i-1}),\ff_1(a_i),\gg_1(b_{j+1}))\\
(\ff\vee_{ij}\gg)_2(a_i,b_{j+1})&:=0\\
     (\ff\vee_{ij}\gg)_2(a_{i-1},\mu_2(a_i,b_{j+1}))&:=\mu_2(\ff_2(a_{i-1},a_i),\gg_1(b_{j+1}))\\
     (\ff\vee_{ij}\gg)_3(a_{i-1},a_i,b_{j+1})&:=0
\end{align*}

If we were to guess how to extend this, we would start by taking
\begin{equation*}
    (\ff\vee_{ij}\gg)_1(\mu_k(a,\dots ,a,b,\dots, b))=\mu_k(\ff_1(a),\dots,\ff_1(a),\gg_1(b),\dots,\gg_1(b))
\end{equation*}
and
\begin{equation*}
    (\ff\vee_{ij} \gg)_{k}(a,\dots,a,b,\dots,b)=0
\end{equation*}
for sequences of $k\geq2$ arrows of the displayed form, starting in $Q_n$ and ending in $Q_m$ (note that we have dropped the indices to avoid a notational overload). Let's check a couple of easy examples.
First, consider the sequence $(a_i,b_{j+1},b_{j+2})$. Assume $n\geq1$ and $m\geq2$ here, so that $\mu_2(b_{j+1},b_{j+2})=0$. Our chosen definitions give
\begin{align*}
    (\ff\vee_{ij}\gg)_1(\mu_3(a_i,b_{j+1},b_{j+2}))&=\mu_3(\ff_1(a_i),\gg_1(b_{j+1}),\gg_1(b_{j+2}))\\
    (\ff\vee_{ij}\gg)_{3}(a_i,b_{j+1},b_{j+2})&=0\\
 (\ff\vee_{ij}\gg)_1(\mu_2(a_i,b_{j+1}))&=\mu_2(\ff_1(a_i),\gg_1(b_{j+1}))\\
 (\ff\vee_{ij}\gg)_2(a_{i},b_{j+1})&=0
\end{align*}
Thus we need to be able to recover $(\ff\vee_{ij}\gg)_2(\mu_2(a_i,b_{j+1}),b_{j+2})$ from these. We can: the following relation must hold.
\begin{align*}
    \mu_1 (\ff\vee_{ij}\gg)_{3}(a_i,b_{j+1},b_{j+2})&+\mu_2((\ff\vee_{ij}\gg)_1(a_i),(\ff\vee_{ij}\gg)_2(b_{j+1},b_{j+2}))\\&+\mu_2((\ff\vee_{ij}\gg)_2(a_i,b_{j+1}),(\ff\vee_{ij}\gg)_1(b_{j+2}))\\
    &+\mu_3((\ff\vee_{ij}\gg)_1(a_i),(\ff\vee_{ij}\gg)_1(b_{j+1}),(\ff\vee_{ij}\gg)_1(b_{j+2}))\\
    =&(\ff\vee_{ij}\gg)_1(\mu_3(a_i,b_{j+1},b_{j+2}))+(\ff\vee_{ij}\gg)_2(a_i,\mu_2(b_{j+1},b_{j+2}))\\
    &+(\ff\vee_{ij}\gg)_2(\mu_2(a_i,b_{j+1}),b_{j+2})
\end{align*}

Using only what we already know or have assumed, this simplifies to
\begin{align*}
  (\ff\vee_{ij}\gg)_2(\mu_2(a_i,b_{j+1}),b_{j+2})= &\mu_2(\ff_1(a_i),\gg_2(b_{j+1},b_{j+2}))
\end{align*}
and so we have no choice in what we define to be $(\ff\vee_{ij}\gg)_2(\mu_2(a_i,b_{j+1}),b_{j+2})$. Another term remains to be found, namely $(\ff\vee_{ij}\gg)_1(\mu_2(\mu_2(a_{i},b_{j+1}),b_{j+2}))$. Fortunately, the functor equation gives
\begin{align*}
\mu_1(\ff\vee_{ij}\gg)_2(\mu_2(a_i,b_{j+1}),b_{j+2})
&+\mu_2((\ff\vee_{ij}\gg)_1(\mu_2(a_i,b_{j+1})),(\ff\vee_{ij}\gg)_1(b_{j+2}))\\
&=(\ff\vee_{ij}\gg)_1(\mu_2(\mu_2(a_i,b_{j+1}),b_{j+2})).
\\
(\ff\vee_{ij}\gg)_1(\mu_2(\mu_2(a_i,b_{j+1}),b_{j+2}))
&=\mu_1\mu_2(\ff_1(a_i),\gg_2(b_{j+1},b_{j+2}))\\
&\quad+\mu_2(\mu_2(\ff_1(a_i),\gg_1(b_{j+1})),\gg_1(b_{j+2})).
\end{align*}

\begin{defn}\label{initdef}
  For a one-crossing sequence $(a_{i-l},\dots,a_i,b_{j+1},\dots,b_{j+s})$, assume $0\leq l\leq n-1$, $1\leq s\leq m$, $0\leq r\leq l$, and $1\leq q\leq s$. Set
   \begin{align*}
&(\mathcal F\vee_{ij}\mathcal G)_{r+s-q+1}
\Bigl(
a_{i-l},\ldots,a_{i-l+r-1},
\mu_{l-r+1+q}(a_{i-l+r},\ldots,a_i,b_{j+1},\ldots,b_{j+q}),
b_{j+q+1},\ldots,b_{j+s}
\Bigr)
\\
&\quad:=
\sum_{\substack{p_1+\cdots+p_{l-r+1}=l+1\\p_1,\ldots,p_{l-r+1}\geq1}}\,
\sum_{\substack{t_1+\cdots+t_q=s\\t_1,\ldots,t_q\geq1}}
\mu_{l-r+1+q}\Bigl(
\\
&\qquad\mathcal F_{p_1}(a_{i-l},\ldots,a_{i-l+p_1-1}),
\ldots,
\mathcal F_{p_{l-r+1}}(a_{i-p_{l-r+1}+1},\ldots,a_i),
\\
&\qquad
\mathcal G_{t_1}(b_{j+1},\ldots,b_{j+t_1}),
\ldots,
\mathcal G_{t_q}(b_{j+s-t_q+1},\ldots,b_{j+s})
\Bigr).
\end{align*}
   Here, $l+1$ and $s$ are the total numbers of $a$- and $b$-arrows, respectively; $r$ is the number of $a$-arrows before the internal $\mu$, $l-r+1$ and $q$ are the numbers of $a$- and $b$-arrows inside it, and $s-q$ is the number of $b$-arrows after it.
\end{defn}

\begin{rmk}
    If $r=0$, $q=s$, and $p=l+s+1$, all terms are inside the internal $\mu$, and we are considering $(\ff\vee_{ij}\gg)_1(\mu_p(a_{i-l},\dots,b_{j+s}))$. All the $p_j$ and $t_j$ are forced to be exactly one, so we recover
    \begin{equation*}
        (\ff\vee_{ij}\gg)_{1}(\mu_{p}(a_{i-l},\dots,b_{j+s}))=\mu_{p}(\ff_1(a_{i-l}),
        \dots,\gg_1(b_{j+s})).
    \end{equation*}
\end{rmk}

We then need to consider sequences which contain more than one nested $\mu$, still working with the one-crossing sequences and bounds of Definition \ref{initdef}. Notice that the only time our desired map will be applied to a non-zero morphism is for sequences of the form
\begin{equation*}
    (a,\dots,a,\mu(a,\dots,a,\mu(a,\dots\dots \mu(a,\dots,a,b,\dots,b)\dots\dots,b),b,\dots,b),b,\dots,b)
\end{equation*}
as all other sequences will contain some term of the form $\mu(a,\dots,a)=0$ or $\mu(b,\dots,b)=0$.

 In order to give explicit maps, we need to first introduce some notation.

We can view all of our sequences of elements in terms of trees. Given a fixed initial and terminal morphism, i.e., $a_{i-l}$ and $b_{j+s}$, any non-trivial morphism in $P(i-l-1,(j+s)')$ (restricting to those that only make one jump $Q_n \to Q_m$) is represented by an expression in the arrows $a_r$ and $b_t$ for $i-l\leq r\leq i$ and $j+1\leq t\leq j+s$, grouped via nested $\mu$s. Here, the prime marks an object in the second factor.

 We can think of such an object as its underlying tree operation applied to the sequence $(a_{i-l},\dots,b_{j+s})$. For example, consider the sequence $(a_{i-1},\mu_2(a_i,b_{j+1}))$. We can view this as
\begin{equation}
    (a_{i-1},T(a_i,b_{j+1}))
\end{equation}
where $T(x,y)=\mu_2(x,y)$. For a more complicated example, consider
\[
    (a_{i-2},\mu_3(a_{i-1},\mu_2(a_i,b_{j+1}),b_{j+2}),b_{j+3}).
\]
We can write this as
\begin{equation}
    (a_{i-2},T(a_{i-1},a_i,b_{j+1},b_{j+2}),b_{j+3})
\end{equation}
where
\begin{equation*}
    T(x,y,z,w)=\mu_3(x,\mu_2(y,z),w).
\end{equation*}
Thus, any non-trivial sequence that we are considering (i.e., one with no $0$s) can be written simply as
\begin{equation*}
(a_{i-l},\ldots,a_{i-l+r-1},T(a_{i-l+r},\ldots,a_i,b_{j+1},\ldots,b_{j+q}),
b_{j+q+1},\ldots,b_{j+s}),
\end{equation*}

If one then studies Definition \ref{initdef}, the equation can be written simply as
\begin{align}\label{non1terms}
&(\mathcal F\vee_{ij}\mathcal G)_{r+s-q+1}
\Bigl(
a_{i-l},\ldots,a_{i-l+r-1},
T(a_{i-l+r},\ldots,a_i,b_{j+1},\ldots,b_{j+q}),
b_{j+q+1},\ldots,b_{j+s}
\Bigr)
\\
&\quad:=
\sum_{\substack{p_1+\cdots+p_{l-r+1}=l+1\\p_1,\ldots,p_{l-r+1}\geq1}}\,
\sum_{\substack{t_1+\cdots+t_q=s\\t_1,\ldots,t_q\geq1}}
T\Bigl(
\notag\\
&\qquad\mathcal F_{p_1}(a_{i-l},\ldots,a_{i-l+p_1-1}),
\ldots,
\mathcal F_{p_{l-r+1}}(a_{i-p_{l-r+1}+1},\ldots,a_i),
\notag\\
&\qquad
\mathcal G_{t_1}(b_{j+1},\ldots,b_{j+t_1}),
\ldots,
\mathcal G_{t_q}(b_{j+s-t_q+1},\ldots,b_{j+s})
\Bigr).\notag
\end{align}
where $T(x_1,\dots,x_{l-r+1+q}):=\mu_{l-r+1+q}(x_1,\dots,x_{l-r+1+q})$.

 In fact, we will see that this formula is `correct' for any tree $T$ and with $r+s-q\geq1$. By this we mean the following: if we define $(\ff\vee_{ij}\gg)_k$ as above for any $k>1$, then there is a natural way to define the remaining values of $(\ff\vee_{ij}\gg)_1$ such that the functor equations are satisfied.

Let us consider a `simple' example: the sequence $(a_{i-1},a_i,b_{j+1},b_{j+2})$ (here we assume $n,m\geq2$). There are a total of 12 ways to break this sequence using $\mu$s such that no term is zero (i.e., we cannot use $(\mu_2(a_{i-1},a_i),b_{j+1},b_{j+2})$).
For sequences with at most one $\mu$:
\begin{align*}
(\ff\vee_{ij}\gg)_4(a_{i-1},a_i,b_{j+1},b_{j+2})&=0,\\
(\ff\vee_{ij}\gg)_3(a_{i-1},\mu_2(a_i,b_{j+1}),b_{j+2})
&=\mu_2(\ff_2(a_{i-1},a_i),\gg_2(b_{j+1},b_{j+2})),\\
(\ff\vee_{ij}\gg)_2(a_{i-1},\mu_3(a_i,b_{j+1},b_{j+2}))
&=\mu_3(\ff_2(a_{i-1},a_i),\gg_1(b_{j+1}),\gg_1(b_{j+2})),\\
(\ff\vee_{ij}\gg)_2(\mu_3(a_{i-1},a_i,b_{j+1}),b_{j+2})
&=\mu_3(\ff_1(a_{i-1}),\ff_1(a_i),\gg_2(b_{j+1},b_{j+2})),\\
(\ff\vee_{ij}\gg)_1(\mu_4(a_{i-1},a_i,b_{j+1},b_{j+2}))
&=\mu_4(\ff_1(a_{i-1}),\ff_1(a_i),\gg_1(b_{j+1}),\gg_1(b_{j+2})).
\end{align*}

Using the stated formula, we will then also have
\begin{align*}
    (\ff\vee_{ij}\gg)_2(a_{i-1},\mu_2(\mu_2(a_i,b_{j+1}),b_{j+2}))=&\mu_2(\mu_2(\ff_2(a_{i-1},a_i),\gg_1(b_{j+1})),\gg_1(b_{j+2}))\\
     (\ff\vee_{ij}\gg)_2(\mu_2(a_{i-1},\mu_2(a_i,b_{j+1})),b_{j+2})=&\mu_2(\ff_1(a_{i-1}),\mu_2(\ff_1(a_i),\gg_2(b_{j+1},b_{j+2})))
\end{align*}
Thus, the remaining terms are\nopagebreak
\begin{align*}
    (\ff\vee_{ij}\gg)_1(\mu_3(a_{i-1},\mu_2(a_i,b_{j+1}),b_{j+2}))=&?\\
    (\ff\vee_{ij}\gg)_1(\mu_2(a_{i-1},\mu_3(a_i,b_{j+1},b_{j+2})))=&?\\
    (\ff\vee_{ij}\gg)_1(\mu_2(\mu_3(a_{i-1},a_i,b_{j+1}),b_{j+2}))=&?\\
    (\ff\vee_{ij}\gg)_1(\mu_2(\mu_2(a_{i-1},\mu_2(a_i,b_{j+1})),b_{j+2}))=&?\\
    (\ff\vee_{ij}\gg)_1(\mu_2(a_{i-1},\mu_2(\mu_2(a_i,b_{j+1}),b_{j+2})))=&?
\end{align*}

The claim is that all of these terms are determined by our previous choices. Indeed, first consider $(\ff\vee_{ij}\gg)_1(\mu_2(a_{i-1},\mu_3(a_i,b_{j+1},b_{j+2})))$. The functor relations give
\begin{align*}
\mu_1(\ff\vee_{ij}\gg)_2(a_{i-1},\mu_3(a_i,b_{j+1},b_{j+2}))
&+\mu_2((\ff\vee_{ij}\gg)_1(a_{i-1}),(\ff\vee_{ij}\gg)_1(\mu_3(a_i,b_{j+1},b_{j+2})))\\
&=(\ff\vee_{ij}\gg)_1(\mu_2(a_{i-1},\mu_3(a_i,b_{j+1},b_{j+2})))\\
&\quad+(\ff\vee_{ij}\gg)_2(a_{i-1},\mu_1\mu_3(a_i,b_{j+1},b_{j+2})).
\end{align*}

Notice that all terms other than $(\ff\vee_{ij}\gg)_1(\mu_2(a_{i-1},\mu_3(a_i,b_{j+1},b_{j+2})))$ are already known (with $(\ff\vee_{ij}\gg)_2(a_{i-1},\mu_1\mu_3(a_i,b_{j+1},b_{j+2}))=(\ff\vee_{ij}\gg)_2(a_{i-1},\mu_2(\mu_2(a_i,b_{j+1}),b_{j+2}))$). One can check this is true for all of the above equations. In fact, what you will find is that
\begin{equation}\label{num1term}
\begin{aligned}
&(\mathcal F\vee_{ij}\mathcal G)_1
\Bigl(\mu_{r+s-q+1}\Bigl(
a_{i-l},\ldots,a_{i-l+r-1},
T(a_{i-l+r},\ldots,a_i,b_{j+1},\ldots,b_{j+q}),
b_{j+q+1},\ldots,b_{j+s}
\Bigr)\Bigr)
\\
&\quad:=\mu_1(\mathcal F\vee_{ij}\mathcal G)_{r+s-q+1}
\Bigl(a_{i-l},\ldots,a_{i-l+r-1},
T(a_{i-l+r},\ldots,a_i,b_{j+1},\ldots,b_{j+q}),
b_{j+q+1},\ldots,b_{j+s}\Bigr)
\\
&\quad+(\mathcal F\vee_{ij}\mathcal G)_{r+s-q+1}
\Bigl(a_{i-l},\ldots,a_{i-l+r-1},
\mu_1T(a_{i-l+r},\ldots,a_i,b_{j+1},\ldots,b_{j+q}),
b_{j+q+1},\ldots,b_{j+s}\Bigr)
\\
&\qquad+\sum\mu_{\alpha+\beta+1}\Biggl(
\mathcal F_{p_1}(a_{i-l},\ldots,a_{i-l+p_1-1}),\ldots,
\mathcal F_{p_\alpha}(a_{i-l+r-u-p_\alpha},\ldots,a_{i-l+r-u-1}),
\\
&\hspace{1cm}(\mathcal F\vee_{ij}\mathcal G)_{u+v+1}
\Bigl(a_{i-l+r-u},\ldots,a_{i-l+r-1},
\\
&\hspace{2cm}T(a_{i-l+r},\ldots,a_i,b_{j+1},\ldots,b_{j+q}),
b_{j+q+1},\ldots,b_{j+q+v}\Bigr),
\\
&\hspace{2cm}\mathcal G_{t_1}(b_{j+q+v+1},\ldots,b_{j+q+v+t_1}),\ldots,
\mathcal G_{t_\beta}(b_{j+s-t_\beta+1},\ldots,b_{j+s})
\Biggr)
\\
&\qquad+\sum(\mathcal F\vee_{ij}\mathcal G)_{r+s-q+1-u-v}
\Bigl(a_{i-l},\ldots,a_{i-l+r-u-1},
\\
&\hspace{1cm}\mu_{u+v+1}\Bigl(
a_{i-l+r-u},\ldots,a_{i-l+r-1},
\\
&\hspace{2cm}T(a_{i-l+r},\ldots,a_i,b_{j+1},\ldots,b_{j+q}),
b_{j+q+1},\ldots,b_{j+q+v}\Bigr),
\\
&\hspace{2cm}b_{j+q+v+1},\ldots,b_{j+s}\Bigr).
\end{aligned}
\end{equation}

It may look horrible, but one should realise it is simply the result of the $A_\infty$-functor relations.

\subsection{The gluing functor}
The aim of this section is to construct a `gluing functor'
\begin{equation}
    \Phi_{ij}:Q_{n+m-1}\to Q_n \bigvee_{ij}Q_m
\end{equation}
which will give the instructions for gluing together $\ff:Q_n \to \cc$ and $\gg:Q_m \to \cc$ along $\ff_0(i)=\gg_0(j)$ to produce
\begin{equation}
    \Phi_{ij}^*(\ff\vee_{ij}\gg):Q_{n+m-1}\to \cc.
\end{equation}
With this in mind, we will need to adjust our target $Q_n\bigvee_{ij}Q_m$ to account for the correct gradings. Write $c_r:r-1\to r$ for the arrows of $Q_{n+m-1}$, with indices modulo $n+m$. Indices in the two target factors are taken modulo $n+1$ and $m+1$, respectively. The object map is
\begin{equation}\label{phiobjects}
    (\Phi_{ij})_0(r):=\begin{cases}
        r_{Q_n},&0\leq r<i,\\
        (j+r-i+1)_{Q_m},&i\leq r<i+m,\\
        (r-m+1)_{Q_n},&i+m\leq r<n+m,
    \end{cases}
\end{equation}
Thus, the object map skips the common object $i_{Q_n}=j_{Q_m}$ and retains exactly $n+m$ objects. Let $f_r$ denote an arrow from $r-1$ to $r$ in either factor, i.e.,
\[
    f_r:=\begin{cases}
        a_r,&r-1,r\in Q_n\subset Q_n\bigvee_{ij}Q_m,\\
        b_r,&r-1,r\in Q_m\subset Q_n\bigvee_{ij}Q_m,
    \end{cases}
\]
where $r$ is read in the original indexing of the relevant factor. At the joined point, we write
\begin{equation}\label{phisplicearrows}
    f_i:=\mu_2(a_i,b_{j+1}),\qquad
    f_{i+m}:=\mu_2(b_j,a_{i+1}).
\end{equation}
Notice that these morphisms are cycles.

To define the chosen grading, shift the retained object indexed by $r$ by an integer $s_r$, and give the common object shift zero in both factors. In our homological convention, a morphism of original degree $d$ has degree $d+t-s$ from $x[s]$ to $y[t]$. Choose the $s_r$ so that
\begin{equation}\label{phishifts}
    |f_r|+s_r-s_{r-1}=|c_r|,
\end{equation}
where $|f_r|$ denotes the degree before shifting. Such a collection of shifts exists by fixing $s_0=0$ and solving recursively around the cycle. Note also that
\[
    \sum_{r=0}^{n+m-1}|f_r|=(1-n)+(1-m)=2-n-m=\sum_{r=0}^{n+m-1}|c_r|.
\]
Equivalently, shift the objects in each factor first, using the same shift at $i\sim j$, and take their strict pushout as defined earlier. We use the same labels for these objects, morphisms, and operations; in this grading $|f_r|=|c_r|$. When composing with the earlier pushout functor, one can apply the functor with the same shifts; for example, the object $x[s]$ is sent to $\bigl((\ff\vee_{ij}\gg)_0(x)\bigr)[s]$.

To construct a valid $\Phi_{ij}$, we need to make use of the following result:

\begin{lemma}\label{lemphihomology}
The homology of the morphism complexes in the pushout $Q_n\bigvee_{ij}Q_m$ is given by
\[
\begin{aligned}
    H_*\bigl(Q_n\bigvee_{ij}Q_m\bigr)\bigl(r,t\bigr)=\begin{cases}
        k\cdot[1_r],&t=r,\\
        k\cdot[f_{r+1}],&t=r+1,\quad r,t\in Q_n\text{ or }r,t\in Q_m,\\
        k\cdot[f_i],&r=(i-1)_{Q_n},\ t=(j+1)_{Q_m},\\
        k\cdot[f_{i+m}],&r=(j-1)_{Q_m},\ t=(i+1)_{Q_n},\\
        0,&\text{otherwise}.
    \end{cases}
\end{aligned}
\]
\end{lemma}

\begin{proof}
Let $z\in\bigl(Q_n\bigvee_{ij}Q_m\bigr)(x,y)_u$ be a cycle. We can write $z$ as a finite sum of trees,
\[
    z=\sum_{\ell=1}^{s}T_\ell(f_{\ell,1},\dots,f_{\ell,r_\ell}),
\]
where the inputs of each tree are composable original $a$- and $b$-arrows. When $f_i$ or $f_{i+m}$ occurs, we expand it as the corresponding $\mu_2$ in \eqref{phisplicearrows}, so that both original arrows are counted. We will say the length of a tree is its number of original input arrows, i.e., the total number of $a$- and $b$-arrows appearing as inputs.

The aim is to write
\begin{equation}\label{cyclewordreduction}
    z=h+\mu_1K,
\end{equation}
where $h$ is either zero or one of the listed identities, arrows $a_l$ or $b_s$, or composites $f_i$ or $f_{i+m}$, with the corresponding endpoints and degree. The idea is to recursively write $z=h_k+\mu_1K_k$ for $k\geq0$, with the largest tree length in $h_k$ strictly decreasing at each step. We will show that this process terminates with the desired expression $z=h+\mu_1K$. We first shorten the trees wherever possible, then separate off these representatives and deal with the remaining cycle.

For example, for $n\geq2$, we have
\[
    \mu_2\bigl(a_{i-1},\mu_2[a_i,b_{j+1}]\bigr)
    =0 +\mu_1\mu_3(a_{i-1},a_i,b_{j+1}).
\]
Or, for example, with $m=1$,
\[
    \mu_2\bigl(\mu_2[a_i,b_{j+1}],b_j\bigr)
    =a_i+\mu_1\mu_3(a_i,b_{j+1},b_j).
\]
In these examples, the correction uses the same original input sequence, while the remaining term has fewer inputs. We use the same idea for a general cycle. Whenever a full-loop operation $\mu_{n+1}(a_{l+1},\dots,a_l)$ or $\mu_{m+1}(b_{l+1},\dots,b_l)$ occurs, we replace it by the corresponding unit and apply strict unitality. This shortens the affected tree or makes it vanish.

For a general such $z=\sum T$, after contracting full loops in $Q_n$ and $Q_m$, we find that different trees can have different lengths, i.e., different numbers of original input arrows. Thus, we can write
\begin{equation}
    z=\sum_{l=0}^{r}z^{(l)},
\end{equation}
where each $z^{(l)}=\sum T^{(l)}$ is a sum of terms of $z$, each of length $l$, and $r$ is the largest length. We count identities as having length zero. We can rewrite this as
\[
    z=z^{(r)}+z^{(<r)},
\]
where $z^{(<r)}=\sum_{l=0}^{r-1}z^{(l)}$ contains the terms of length strictly less than $r$. Explicitly,
\[
    z^{(r)}=\sum T(f_1,\dots,f_r).
\]

Applying $\mu_1$ cannot increase the total number of original inputs. Since $z$ is a cycle, we find
\[
    \mu_1z=\sum_{l=0}^{r}(\mu_1z)^{(l)}=0.
\]
In particular, we have
\[
    (\mu_1z^{(r)})^{(r)}=(\mu_1z)^{(r)}=0.
\]
We now show that, for $r\geq3$, we can find a sum of trees $K^{(r)}$ with $r$ original inputs such that
\[
    (\mu_1K^{(r)})^{(r)}=z^{(r)}.
\]
Once this is established, we can reduce $z$ to 
\[
    z+\mu_1K^{(r)}=z^{(<r)}+(\mu_1K^{(r)})^{(<r)}.
\]
Repeating this process removes all terms of length at least three. It then remains to check that the terms of lengths zero, one and two reduce to zero or the listed representative with the required endpoints and degree.

Fix a sequence $(f_1,\dots,f_r)$ with $c$ crossings. We induct on the number of crossings and, for a fixed number of crossings, on the number $V$ of vertices, proving the statement for all input sequences at each step. We may assume that every tree in $z^{(r)}$ has $V$ vertices, since $\mu_1$ at this length increases the number of vertices by one. With no crossings, every tree of length $r\geq3$ either reduces to a shorter expression or vanishes. If $V=1$ and $r\geq3$, the only tree is $\mu_r(f_1,\dots,f_r)$. At each crossing, its differential contains the nonzero term
\[
    \mu_{r-1}\Bigl(f_1,\dots,f_{j-1},\mu_2[f_j,f_{j+1}],f_{j+2},\dots,f_r\Bigr).
\]
Thus, there is no nonzero cycle in this case.

Assume the claim holds for fewer than $c$ crossings and, with $c$ crossings, for fewer than $V$ vertices. Every inner subtree has at most the same number of crossings and fewer vertices. We can write
\[
    z^{(r)}=\sum_{\alpha}\mu_{r_\alpha}(x^\alpha_1,\dots,x^\alpha_{r_\alpha})
    =\sum_{q=2}^{r}z^{(r)}_{(q)},
\]
where $z^{(r)}_{(q)}$ collects the terms whose outer operation is $\mu_q$. Let $p$ be the largest $q$ for which $z^{(r)}_{(q)}\neq0$. Splitting an outer operation decreases its number of inputs. Thus, since $(\mu_1z^{(r)})^{(r)}=0$, we must have that
\[
    \bigl(\mu_1z^{(r)}_{(p)}\bigr)^{(r)}_{(p)}
    =\left(\sum_{\alpha:\,r_\alpha=p}\sum_{j=1}^{p}
    \mu_p\Bigl(x^\alpha_1,\dots,x^\alpha_{j-1},
    [\mu_1x^\alpha_j],x^\alpha_{j+1},\dots,x^\alpha_p\Bigr)\right)^{(r)}=0.
\]
For each $i$, fix the input arrows $f^i_1,\dots,f^i_{l_i}$ of $x_i^\alpha$, while still summing over all terms with these input sequences. Since the sum above vanishes and we work over a field, we may assume, up to adding a boundary at these fixed input lengths, that $(\mu_1x_j^\alpha)^{(l_j)}=0$ for every $\alpha,j$ ( see for example {\cite{LV}}).

Each $x_i^\alpha$ is a sum of trees with $v_i^\alpha$ vertices, and $\sum_{i=1}^{p}v_i^\alpha=V-1$. Since each $v_i^\alpha<V$, the inductive assumption allows us to write
\[
    x_i^\alpha=(\mu_1K_i^\alpha)^{(l_i)}
\]
whenever $l_i\geq3$. The corresponding summand of $z^{(r)}_{(p)}$ satisfies
\[
    \Bigl(\mu_1\mu_p(x_1^\alpha,\dots,K_i^\alpha,\dots,x_p^\alpha)\Bigr)^{(r)}_{(p)}
    =\mu_p(x_1^\alpha,\dots,x_i^\alpha,\dots,x_p^\alpha).
\]
Apply this once to each summand having an input of length at least three, and let $K$ be the sum of the resulting primitives. Then we can write
\begin{align*}
    (z^{(r)}+\mu_1K)^{(r)}
    &=\sum_{\substack{\alpha:\,r_\alpha=p\\l_1,\dots,l_p\leq2}}
    \mu_p(x_1^\alpha,\dots,x_p^\alpha)
    +\sum_{q=2}^{p-1}(z^{(r)}+\mu_1K)^{(r)}_{(q)}.
\end{align*}
The inputs in the first sum are original arrows, $\mu_2(a_i,b_{j+1})$, or $\mu_2(b_j,a_{i+1})$. Next, consider $\mu_p(x_1^\alpha,\dots,x_p^\alpha)$ with each input of length $l_i\leq2$. We can write
\begin{align*}
    \mu_p(x_1^\alpha,\dots,x_p^\alpha)
    &=\mu_{r-s}\Bigl(f_1,\dots,f_{j_1-1},\mu_2[f_{j_1},f_{j_1+1}],f_{j_1+2},\dots,\\
    &\qquad f_{j_s-1},\mu_2[f_{j_s},f_{j_s+1}],f_{j_s+2},\dots,f_r\Bigr).
\end{align*}
We denote this tree by $T_J$, where $J=\{j_1,\dots,j_s\}$ records the crossings paired by $\mu_2$. Here $s=r-p\leq c$. After cancelling repeated terms modulo two, we can write
\[
    \sum_{\substack{\alpha:\,r_\alpha=p\\l_1,\dots,l_p\leq2}}\mu_p(x^\alpha_1,\dots,x^\alpha_p)
    =\sum_{|J|=r-p}T_J,
\]
where the sum runs only over the trees that remain. Choose some crossing $j_1'$ and set
\[
    H_0=\sum_{\substack{|J|=r-p\\j'_1\in J}}T_{J\setminus\{j'_1\}}.
\]
Every term of $H_0$ has $p+1=r-s+1$ outer inputs. Write
\[
    H_0=\sum_\alpha\mu_{p+1}(x_1^\alpha,\dots,x_{p+1}^\alpha),
\]
where these new inputs have lengths $l_i\in\{1,2\}$. Explicitly, we have replaced
\[
    \mu_p\Bigl(\ldots,\mu_2[f_{j'_1},f_{j'_1+1}],\ldots\Bigr)
    \quad\text{by}\quad
    \mu_{p+1}\Bigl(\ldots,f_{j'_1},f_{j'_1+1},\ldots\Bigr),
\]
with the other inputs unchanged. So we have
\begin{align*}
    z^{(r)}
    &=\sum_{|J|=r-p}T_J
    +\sum_{q=2}^{p-1}z^{(r)}_{(q)}.
\end{align*}
Adding $(\mu_1H_0)^{(r)}$, we obtain
\begin{align*}
    (z^{(r)}+\mu_1H_0)^{(r)}
    &=\sum_{|J|=r-p}T_J
    +\sum_{\substack{|J|=r-p\\j'_1\in J}}(\mu_1T_{J\setminus\{j'_1\}})^{(r)}\\
    &\quad+\sum_{q=2}^{p-1}z^{(r)}_{(q)}.
\end{align*}
Expanding the middle sum gives
\begin{align*}
    &\sum_{\substack{|J|=r-p\\j'_1\in J}}(\mu_1T_{J\setminus\{j'_1\}})^{(r)}\\
    &=\sum_{\substack{|J|=r-p\\j'_1\in J}}\Bigl[
    T_J+\sum_{\substack{j\text{ a crossing}\\j\notin J}}
    T_{(J\setminus\{j'_1\})\cup\{j\}}\Bigr]\\
    &\quad+\sum_\alpha\sum_{\substack{1\leq t\leq p\\l_t+l_{t+1}\geq3}}
    \mu_p\Bigl(x_1^\alpha,\dots,x_{t-1}^\alpha,
    \mu_2[x_t^\alpha,x_{t+1}^\alpha],x_{t+2}^\alpha,\dots,x_{p+1}^\alpha\Bigr)\\
    &\quad+\sum_{q=2}^{p-1}(\mu_1H_0)^{(r)}_{(q)}.
\end{align*}
Comparing the coefficients of $T_{S\cup\{j'_1\}}$ in $(\mu_1z^{(r)})^{(r)}=0$, with $|S|=r-p$ and $j'_1\notin S$, shows that the first sum on the RHS equals $z^{(r)}_{(p)}$. Indeed, we have
\begin{align*}
    &\sum_{\substack{|J|=r-p\\j'_1\in J}}\Bigl[
    T_J+\sum_{\substack{j\text{ a crossing}\\j\notin J}}
    T_{(J\setminus\{j'_1\})\cup\{j\}}\Bigr]\\
    &=\sum_{\substack{|J|=r-p\\j'_1\in J}}T_J
    +\sum_{\substack{|S|=r-p\\j'_1\notin S}}
    \Biggl[\sum_{\substack{j\in S\\T_{(S\setminus\{j\})\cup\{j'_1\}}\text{ occurs}}}1\Biggr]T_S\\
    &=\sum_{\substack{|J|=r-p\\j'_1\in J}}T_J
    +\sum_{\substack{|S|=r-p\\j'_1\notin S\\T_S\text{ occurs}}}T_S\\
    &=\sum_{|J|=r-p}T_J=z^{(r)}_{(p)}.
\end{align*}
Here, `occurs' means appears in $z^{(r)}_{(p)}$. Thus,
\begin{align*}
    (z^{(r)}+\mu_1H_0)^{(r)}
    &=\sum_\alpha\sum_{\substack{1\leq t\leq p\\l_t+l_{t+1}\geq3}}
    \mu_p\Bigl(x_1^\alpha,\dots,x_{t-1}^\alpha,
    \mu_2[x_t^\alpha,x_{t+1}^\alpha],x_{t+2}^\alpha,\dots,x_{p+1}^\alpha\Bigr)\\
    &\quad+\sum_{q=2}^{p-1}(z^{(r)}+\mu_1H_0)^{(r)}_{(q)}.
\end{align*}
The first sum may be nonzero, but each $x_i^\alpha$ is closed, so
\[
    \mu_1\mu_2(x_t^\alpha,x_{t+1}^\alpha)=0.
\]
Each of these $\mu_2$ terms has $l_t+l_{t+1}\geq3$ original inputs and fewer than $V$ vertices. Thus, by induction, we can write
\[
    \mu_2(x_t^\alpha,x_{t+1}^\alpha)
    =(\mu_1K_t^\alpha)^{(l_t+l_{t+1})}.
\]
Set
\[
    K'=\sum_\alpha\sum_{\substack{1\leq t\leq p\\l_t+l_{t+1}\geq3}}
    \mu_p\Bigl(x_1^\alpha,\dots,x_{t-1}^\alpha,K_t^\alpha,
    x_{t+2}^\alpha,\dots,x_{p+1}^\alpha\Bigr).
\]
We can then write
\[
    \left(z^{(r)}+\mu_1(H_0+K')\right)^{(r)}
    =\sum_{q=2}^{p-1}\sum_\alpha\mu_q(x^\alpha_1,\dots,x^\alpha_q).
\]
In particular, every outer operation now has fewer than $p$ inputs. We repeat this argument, decreasing the largest $q$ for which $z^{(r)}_{(q)}\neq0$, until no terms remain. Adding all the boundary corrections gives
\[
    (\mu_1K^{(r)})^{(r)}=z^{(r)}.
\]
We can then replace $z$ by
\[
    z_1:=z+\mu_1K^{(r)}
    =z^{(<r)}+(\mu_1K^{(r)})^{(<r)}.
\]
Then $\mu_1z_1=0$, $[z_1]=[z]$, and every tree in $z_1$ has fewer than $r$ original input arrows. We repeat this procedure until all remaining trees have at most two original inputs. These reduced trees must be identities, individual arrows, or the composites $\mu_2(a_i,b_{j+1})$ and $\mu_2(b_j,a_{i+1})$, i.e., the listed representatives. Keeping track of the corrections, we can thus write
\[
    z=h+\mu_1K,
\]
where $h$ is either zero or the listed representative with the required endpoints and degree. It remains to check that the listed representatives are cycles and give nonzero homology classes. Indeed, $\mu_1(1_i)=0$, $\mu_1(a_r)=0$ and $\mu_1(b_r)=0$, and
\[
    \mu_1\mu_2(a_i,b_{j+1})=\mu_1\mu_2(b_j,a_{i+1})=0.
\]
If $h=\mu_1L$, reducing $L$ as above preserves $\mu_1L=h$ and leaves only closed terms with at most two original inputs. Thus $h=0$, proving that the listed representatives give nonzero homology classes.
\end{proof}

\begin{rmk}
    This shows that the cycles that are not boundaries are sparse, and it will make constructing $\Phi_{ij}$ much simpler. Indeed, we find
    \begin{align}
    &H_*\bigl(Q_n\bigvee_{ij}Q_m\bigr)\bigl((\Phi_{ij})_0(r),(\Phi_{ij})_0(t)\bigr)\\
    &\quad=\begin{cases}
        k\cdot[1_{(\Phi_{ij})_0(r)}],&t=r,\\
        k\cdot[f_{r+1}],&t=r+1\pmod{n+m},\\
        0,&\text{otherwise}.
    \end{cases}\notag
\end{align}
\end{rmk}

We can now show the following.
\begin{lemma}
    There exists a well-defined strictly unital $A_\infty$-functor $\Phi_{ij}:Q_{n+m-1}\to Q_n\bigvee_{ij}Q_m$ with $(\Phi_{ij})_0$ given by \eqref{phiobjects}, and with
    \begin{equation}
        (\Phi_{ij})_1(c_{r+1}):=\begin{cases}
            f_{r+1},&r+1\neq i,i+m,\\
            f_i,&r+1=i,\\
            f_{i+m},&r+1=i+m.
        \end{cases}
    \end{equation}
    Here, $0\leq r<n+m$, and the cases are read modulo $n+m$.
\end{lemma}

\begin{proof}
    First of all, we set $\Phi_k:=(\Phi_{ij})_k$ and take
    \begin{equation}
        \Phi_k(c_{i-r},\dots,c_i,\dots, c_{i+s}):=\mu_{k+1}(a_{i-r},\dots, a_i,b_{j+1},\dots,b_{j+s+1})
    \end{equation}
    where $k=r+s+1$, $0\leq r\leq n-1$, and $0\leq s\leq m-1$, and similarly for sequences containing a single $c_{i+m}$ term. We also set
    \[
        \Phi_k(c_{v+1},\dots,c_{v+k})=0
    \]
    for $k\geq2$ whenever the images of the input arrows under $\Phi_1$ lie entirely in $Q_n$ or entirely in $Q_m$, i.e., the sequence contains neither $c_i$ nor $c_{i+m}$. Furthermore, we take $\Phi_k(\dots,1,\dots)=0$ for $k>1$ and $\Phi_1(1)=1$. It is a simple check that these definitions satisfy the $A_\infty$-functor relations on the sequences just described. For sequences whose images make no crossings, this is clear. For a single crossing through $c_i$, we have
\begin{align*}
    &\sum_{q=1}^{k}\sum_{\substack{p_1+\cdots+p_q=k\\p_1,\dots,p_q\geq1}}
    \mu_q\Bigl(\Phi_{p_1}[c_{i-r},\dots,c_{i-r+p_1-1}],\dots,
    \Phi_{p_q}[c_{i+s-p_q+1},\dots,c_{i+s}]\Bigr)\\
    &=\sum_{u=0}^{r}\sum_{v=0}^{s}
    \mu_{k-u-v}\Bigl(\Phi_1(c_{i-r}),\dots,\Phi_1(c_{i-u-1}),\\
    &\qquad \Phi_{u+v+1}[c_{i-u},\dots,c_{i+v}],
    \Phi_1(c_{i+v+1}),\dots,\Phi_1(c_{i+s})\Bigr)\\
    &=\sum_{u=0}^{r}\sum_{v=0}^{s}
    \mu_{k-u-v}\Bigl(a_{i-r},\dots,a_{i-u-1},\\
    &\qquad \mu_{u+v+2}[a_{i-u},\dots,a_i,b_{j+1},\dots,b_{j+v+1}],
    b_{j+v+2},\dots,b_{j+s+1}\Bigr).
\end{align*}
The case containing only the joining arrow $c_{i+m}$ is similar, and the relations for sequences containing neither joining arrow hold directly.

We now construct $\Phi$ by induction on the number of crossings in the image and, for a fixed number of crossings, on the input length. To treat a sequence with $c+1$ crossings and $k$ inputs, assume that $\Phi$ has been constructed and satisfies the functor equations on all sequences with at most $c$ crossings, and on all sequences with $c+1$ crossings and fewer than $k$ inputs. Thus, every component on a proper subsequence is already defined.

We need $\Phi_{ij}$ to satisfy
\begin{align*}
    &\sum_{q=1}^{k}\sum_{\substack{p_1+\cdots+p_q=k\\p_1,\dots,p_q\geq1}}
    \mu_q\Bigl(\Phi_{p_1}[c_{v+1},\dots,c_{v+p_1}],\dots,
    \Phi_{p_q}[c_{v+k-p_q+1},\dots,c_{v+k}]\Bigr)\\
    &\quad+\sum_{s=1}^{k}\sum_{r=0}^{k-s}
    \Phi_{k-s+1}\Bigl(c_{v+1},\dots,c_{v+r},
    \mu_s[c_{v+r+1},\dots,c_{v+r+s}],
    c_{v+r+s+1},\dots,c_{v+k}\Bigr)=0.
\end{align*}
For $k<n+m$, every internal $\mu$ operation vanishes. For $k>n+m$, the only nonzero ones are $\mu_{n+m}$, which produce units, and the resulting component $\Phi_{k-(n+m)+1}$ vanishes by strict unitality. When $k=n+m$, the second sum is $\Phi_1(1_v)=1_{\Phi_0(v)}$.

Thus, we need to choose $\Phi_k$ such that
\[
    \mu_1\Phi_k(c_{v+1},\dots,c_{v+k})=W(c_{v+1},\dots,c_{v+k}),
\]
where
\begin{align*}
    W(c_{v+1},\dots,c_{v+k})
    &:=\sum_{q=2}^{k}\sum_{\substack{p_1+\cdots+p_q=k\\p_1,\dots,p_q\geq1}}
    \mu_q\Bigl(\Phi_{p_1}[c_{v+1},\dots,c_{v+p_1}],\dots,\\
    &\qquad \Phi_{p_q}[c_{v+k-p_q+1},\dots,c_{v+k}]\Bigr)
    +\begin{cases}
        1_{\Phi_0(v)},&k=n+m,\\
        0,&k\neq n+m.
    \end{cases}
\end{align*}
Every term in $W$ is already known: each block has fewer than $k$ inputs and at most $c+1$ crossings. It remains to show that $W$ is a boundary. First, writing each $\Phi_{p_j}$ for its value on the corresponding block, the $A_\infty$-relations give
\begin{align*}
    \mu_1W
    &=\sum_{q=2}^{k}\sum_{\substack{p_1+\cdots+p_q=k\\p_1,\dots,p_q\geq1}}\sum_{j=1}^{q}
    \mu_q\Bigl(\Phi_{p_1},\dots,\Phi_{p_{j-1}},\\
    &\qquad\Bigl[\mu_1\Phi_{p_j}
    +\sum_{\ell=2}^{p_j}\sum_{\substack{u_1+\cdots+u_\ell=p_j\\u_1,\dots,u_\ell\geq1}}
    \mu_\ell[\Phi_{u_1},\dots,\Phi_{u_\ell}]
    +\delta_{p_j,n+m}1\Bigr],\\
    &\qquad\Phi_{p_{j+1}},\dots,\Phi_{p_q}\Bigr)=0.
\end{align*}
Each bracket is the functor equation on $p_j<k$ inputs, so it vanishes by the inductive assumption. The added unit terms cancel: indeed, for $k>n+m$, strict unitality reduces their sum to
\begin{align*}
    &\Phi_{k-n-m}(c_{v+n+m+1},\dots,c_{v+k})\\
    &\qquad+\Phi_{k-n-m}(c_{v+1},\dots,c_{v+k-n-m})=0,
\end{align*}
since the indices are read modulo $n+m$ and we work in characteristic two.

For $k\neq n+m$, a summand with $q$ blocks has degree
\begin{align*}
    |W|&=q-2+\sum_{s=1}^{q}
    \bigl|\Phi_{p_s}(c_{v+p_1+\cdots+p_{s-1}+1},\dots,c_{v+p_1+\cdots+p_s})\bigr|\\
    &=k-2+\sum_{r=1}^{k}|c_{v+r}|.
\end{align*}
By Lemma \ref{lemphihomology}, the relevant morphism space has zero homology when $2\leq k<n+m$. For $k>n+m$, $W$ could represent a non-trivial homology class only when $k=\alpha(n+m)$ or $k=\alpha(n+m)+1$. In the first case, $\alpha>1$, and
\[
    |W|=\alpha(n+m)-2-\alpha((n+m)-2)=2\alpha-2\neq0,
\]
since all but two $c_p$ in each full cycle have degree $-1$. In the second case, $\alpha\geq1$, and
\[
    |W|=2\alpha-1+|c_{v+1}|\neq|c_{v+1}|.
\]
Thus, $W$ is a boundary whenever $k\neq n+m$.
It remains then to check the case of $k=n+m$, here we make the following choice on all sequences with two crossings. Consider $(c_{i-r},\dots,c_{i+m+s})$, where $0\leq r,s\leq n-1$ and $k=r+m+s+1$. We set
\begin{align*}
    &\Phi_k(c_{i-r},\dots,c_{i+m+s})\\
    &\quad:=\sum_{u=0}^{r}\sum_{v=0}^{m}
    \mu_{k-u-v+1}\Bigl(a_{i-r},\dots,a_{i-u-1},\\
    &\qquad\mu_{u+v+2}[a_{i-u},\dots,a_i,b_{j+1},\dots,b_{j+v+1}],
    b_{j+v+2},\dots,b_{j+m+1},a_{i+1},\dots,a_{i+s+1}\Bigr).
\end{align*}
We make the `left' choice, with the first crossing inside the inner $\mu$. Applying the $A_\infty$-relations, the terms with two nested operations both crossing only from $Q_n$ to $Q_m$ cancel in pairs. Contracting the full $b$-cycle gives the last term below:
\begin{align*}
    &\mu_1\Phi_k(c_{i-r},\dots,c_{i+m+s})\\
    &=\sum_{q=2}^{k}\sum_{\substack{p_1+\cdots+p_q=k\\p_1,\dots,p_q\geq1}}
    \mu_q\Bigl(\Phi_{p_1}[c_{i-r},\dots,c_{i-r+p_1-1}],\dots,\\
    &\qquad\Phi_{p_q}[c_{i+m+s-p_q+1},\dots,c_{i+m+s}]\Bigr)
    +\mu_{r+s+2}(a_{i-r},\dots,a_{i+s+1}).
\end{align*}
The relations in $Q_n$ give
\[
    \mu_{r+s+2}(a_{i-r},\dots,a_{i+s+1})
    =\begin{cases}
        1_{\Phi_0(i-r-1)},&r+s=n-1,\\
        0,&r+s\neq n-1.
    \end{cases}
\]
Since $r+s=n-1$ implies that $k=n+m$, this gives $\mu_1\Phi_k=W$. For sequences with more than two crossings, we have $k>n+m$ thus this completes the inductive construction.

\end{proof}

\begin{rmk}
    We have chosen to work from the left. At each stage of the construction, we choose how to break up the compositions, each choice corresponds to vertices of the Stasheff polyhedra. For example, with $n=m=1$, $i=1$ and $j=0$, we know
    \begin{align*}
        &\mu_1\Phi_2(c_1,c_0)+\mu_2(\Phi_1(c_1),\Phi_1(c_0))=\Phi_1(1_0)\\
        \implies\quad&\mu_1\Phi_2(c_1,c_0)+\mu_2(\mu_2(a_1,b_1),\mu_2(b_0,a_0))=1_0
    \end{align*}
    If we choose to break up the middle term from the left, then we find
    \begin{equation*}
        \mu_1\Phi_2(c_1,c_0)=\mu_1\mu_3\bigl(\mu_2(a_1,b_1),b_0,a_0\bigr)+\mu_1 \mu_2\bigl(\mu_3(a_1,b_1,b_0),a_0\bigr)
    \end{equation*}
    Thus, we can set
    \begin{equation*}
        \Phi_2(c_1,c_0):=\mu_3\bigl(\mu_2(a_1,b_1),b_0,a_0\bigr)+\mu_2\bigl(\mu_3(a_1,b_1,b_0),a_0\bigr)
    \end{equation*}
For the remainder of the text, we write $\Phi_{ij}$ for a fixed choice of such a functor.
\end{rmk}

\subsection{The pullback}
We can now define the operation of gluing two functors together along a common object.
\begin{defn}
    Given two $A_\infty$-functors $\ff:Q_n\to\cc$ and $\gg:Q_m\to\cc$, with $\ff_0(i)=\gg_0(j)$, we define their gluing (along $\ff_0(i)=\gg_0(j)$) as
    \begin{equation}
        \ff\#_{ij}\gg:=\Phi_{ij}^*(\ff\vee_{ij}\gg):Q_{n+m-1}\to\cc,
    \end{equation}
    where $\Phi_{ij}$ and $\ff\vee_{ij}\gg$ are the explicit functors constructed above.
\end{defn}

We should explain with a couple of examples what this resulting functor actually is.
\begin{eg}
Assume to start that we have $\ff:Q_1\to\cc$ and $\gg:Q_1\to\cc$, with $\ff_0(1)=\gg_0(0)$. These represent homotopy equivalences in $\cc$, so the resulting functor $\ff\#_{10}\gg$ should represent the fact that the composition of two homotopy equivalences is again a homotopy equivalence. Indeed it does: we have
\begin{align*}
    (\ff\#_{10}\gg)_1(c_1)
    &=(\ff\vee_{10}\gg)_1\bigl((\Phi_{10})_1(c_1)\bigr)\\
    &=(\ff\vee_{10}\gg)_1\bigl(\mu_2(a_1,b_1)\bigr)\\
    &=\mu_2\bigl(\ff_1(a_1),\gg_1(b_1)\bigr).
\end{align*}
Similarly, one can check that $(\ff\#_{10}\gg)_1(c_0)=\mu_2(\gg_1(b_0),\ff_1(a_0))$. Next, one finds
\begin{align*}
    &\mu_1(\ff\#_{10}\gg)_2(c_1,c_0)\\
    &=\mu_1(\ff\vee_{10}\gg)_1\Bigl[
        \mu_3\bigl(\mu_2(a_1,b_1),b_0,a_0\bigr)
        +\mu_2\bigl(\mu_3(a_1,b_1,b_0),a_0\bigr)\Bigr]\\
    &=\mu_1\Bigl[
        \ff_2(a_1,a_0)
        +\mu_2\bigl(\mu_2(\ff_1(a_1),\gg_2(b_1,b_0)),\ff_1(a_0)\bigr)\\
    &\qquad+\mu_3\bigl(\mu_2(\ff_1(a_1),\gg_1(b_1)),\gg_1(b_0),\ff_1(a_0)\bigr)\\
    &\qquad+\mu_2\bigl(\mu_3(\ff_1(a_1),\gg_1(b_1),\gg_1(b_0)),\ff_1(a_0)\bigr)\Bigr],
\end{align*}
where we have used the fact that, by our definition,
\[
    (\ff\vee_{10}\gg)_2\bigl(\mu_2(a_1,b_1),\mu_2(b_0,a_0)\bigr)=0.
\]
Now use the functor equations
\begin{align*}
    \mu_1\ff_2(a_1,a_0)
    &=1_{\ff_0(0)}+\mu_2\bigl(\ff_1(a_1),\ff_1(a_0)\bigr),\\
    \mu_1\gg_2(b_1,b_0)
    &=1_{\gg_0(0)}+\mu_2\bigl(\gg_1(b_1),\gg_1(b_0)\bigr),
\end{align*}
to simplify down to
\begin{align*}
    &\mu_1(\ff\#_{10}\gg)_2(c_1,c_0)\\
    &\quad=1_{\ff_0(0)}
    +\mu_2\Bigl(\mu_2(\ff_1(a_1),\gg_1(b_1)),
    \mu_2(\gg_1(b_0),\ff_1(a_0))\Bigr).
\end{align*}
A similar calculation holds for $\mu_1(\ff\#_{10}\gg)_2(c_0,c_1)$. Thus we recover that the composition of two homotopy equivalences is a homotopy equivalence.
\end{eg}

\begin{eg}
Consider $\ff:Q_2\to\cc$ and $\gg:Q_1\to\cc$, with $\ff_0(2)=\gg_0(0)$. Then the gluing $\ff\#_{20}\gg:Q_2\to\cc$ will realise the property that, given an exact triangle $\Delta$, we can obtain a new exact triangle by replacing an object with a homotopy equivalent object. More generally, triangles isomorphic to exact triangles are also exact. One should picture
\begin{equation*}
    \begin{tikzcd}[column sep=large,row sep=large]
        & \ff_0(1)\ar[dr,"\ff_1(a_2)"] & &\\
        \ff_0(0)\ar[ur,"\ff_1(a_1)"] & &
        \ff_0(2)=\gg_0(0)\ar[ll,"\ff_1(a_0)"]
        \ar[r,bend left,"\gg_1(b_1)"] &
        \gg_0(1)\ar[l,bend left,"\gg_1(b_0)"]
    \end{tikzcd}
\end{equation*}
being glued to form an exact triangle:
\begin{equation*}
    \begin{tikzcd}[column sep=large,row sep=large]
        & \ff_0(1)\ar[dr,"{(\ff\#_{20}\gg)_1(c_2)}"] &\\
        \ff_0(0)\ar[ur,"{(\ff\#_{20}\gg)_1(c_1)}"] & &
        \gg_0(1)\ar[ll,"{(\ff\#_{20}\gg)_1(c_0)}"]
    \end{tikzcd}
\end{equation*}
\end{eg}

\begin{eg}
 If we glue two functors $\ff,\gg:Q_2\to\cc$, we recover a form of the octahedral axiom.
For $\ff_0(2)=\gg_0(0)$, one should picture
\begin{equation*}
    \begin{tikzcd}[column sep=2.5em,row sep=6em]
        \ff_0(0)\ar[rr,"\ff_1(a_1)"] & &
        \ff_0(1)\ar[dl,"\ff_1(a_2)"]
        \ar[dd,"{(\ff\#_{20}\gg)_1(c_2)}"]\\
        & \ff_0(2)=\gg_0(0)\ar[ul,"\ff_1(a_0)"]
        \ar[dr,"\gg_1(b_1)"] &\\
        \gg_0(2)\ar[ur,"\gg_1(b_0)"]
        \ar[uu,"{(\ff\#_{20}\gg)_1(c_0)}"] & &
        \gg_0(1)\ar[ll,"\gg_1(b_2)"]
    \end{tikzcd}
\end{equation*}
The outer square gives $\ff\#_{20}\gg:Q_3\to\cc$, with the two original triangles inside.
To see the other pair of triangles, consider
\[
    \mu_2\bigl(\ff_1(a_2),\gg_1(b_1)\bigr)
    =(\ff\#_{20}\gg)_1(c_2):\ff_0(1)\to\gg_0(1).
\]
Extend this morphism to an exact triangle (assuming that $\cc$ is pre-triangulated), represented by a functor $\hh:Q_2\to\cc$. Writing $X=\hh_0(2)$, the octahedral axiom gives a second exact triangle, which we represent by $\kk:Q_2\to\cc$. With the shifts written explicitly, these triangles are
\begin{align*}
    \hh:\quad &\ff_0(1)\xrightarrow{(\ff\#_{20}\gg)_1(c_2)}\gg_0(1)
    \longrightarrow X\longrightarrow\ff_0(1)[1],\\
    \kk:\quad &\ff_0(0)[1]\longrightarrow X\longrightarrow\gg_0(2)
    \longrightarrow\ff_0(0)[2].
\end{align*}
After the correct rotations and shifts, these rearrange into the same outer square in $\text{H}_0\text{Tw}\cc$, now with $X$ as the common object. The point is that no choice of $X$ is needed for the gluing, this information is already contained in the square and the higher components of $\ff\#_{20}\gg$.
\end{eg}

\section{The homotopy structure}\label{sec:operadic-structure}
\subsection{Two gluings}
If one were to consider three functors, for simplicity take $\ff,\gg,\hh:Q_1\to\cc$ and assume they can be glued together, i.e., $\ff_0(1)=\gg_0(0)$ and $\gg_0(1)=\hh_0(0)$. Then we can compose the three homotopy equivalences $\ff_1(a_1)$, $\gg_1(b_1)$ and $\hh_1(c_1)$. However, there are two ways we can make this composition: either $\mu_2(\ff_1(a_1),\mu_2(\gg_1(b_1),\hh_1(c_1)))$ or $\mu_2(\mu_2(\ff_1(a_1),\gg_1(b_1)),\hh_1(c_1))$. Fortunately, they are homotopic, via
\begin{equation}
\begin{aligned}
    &\mu_1\mu_3(\ff_1(a_1),\gg_1(b_1),\hh_1(c_1))\\
    &\quad=\mu_2(\ff_1(a_1),\mu_2(\gg_1(b_1),\hh_1(c_1)))
    +\mu_2(\mu_2(\ff_1(a_1),\gg_1(b_1)),\hh_1(c_1)).
\end{aligned}
\end{equation}
We can understand this in terms of our earlier constructions. It starts with realising that there are two orders in which to take a double pushout. Fix strictly positive integers $n,m,\ell$, and assume we want to glue two functors $\ff:Q_n\to\cc$ and $\hh:Q_\ell\to\cc$ to a single $\gg:Q_m\to\cc$, with $\ff_0(i)=\gg_0(j)$ and $\gg_0(i')=\hh_0(j')$. Consider the double pushout
\begin{equation}
\begin{tikzcd}[column sep=large,row sep=large]
    & *\ar[d,"i'"]\ar[r,"j'"] & Q_\ell\ar[d]\\
    Q_n\ar[r] & Q_n\bigvee_{ij}Q_m\ar[r] &
    (Q_n\bigvee_{ij}Q_m)\bigvee_{i'j'}Q_\ell\\
    *\ar[u,"i"]\ar[r,"j"] & Q_m\ar[u] &
\end{tikzcd}
\end{equation}
We first build $Q_n\bigvee_{ij}Q_m$, then use this to build $(Q_n\bigvee_{ij}Q_m)\bigvee_{i'j'}Q_\ell$. However, we could also swap the order and consider
\begin{equation}
\begin{tikzcd}[column sep=large,row sep=large]
    & *\ar[d,"j"]\ar[r,"i"] & Q_n\ar[d]\\
    Q_\ell\ar[r] & Q_m\bigvee_{i'j'}Q_\ell\ar[r] &
    Q_n\bigvee_{ij}(Q_m\bigvee_{i'j'}Q_\ell)\\
    *\ar[u,"j'"]\ar[r,"i'"] & Q_m\ar[u] &
\end{tikzcd}
\end{equation}
\begin{rmk}
The two double pushouts are strictly isomorphic. In fact, we may as well take them to be identical. This follows from the definition of our pushout category. Both categories have the same finite set of objects, which we can order in the same way. Furthermore, since morphisms and the $A_\infty$-structures are taken freely modulo the relations imposed by $Q_n,Q_m$ and $Q_\ell$, it is clear that the two constructions agree.
\end{rmk}

So where exactly does the discrepancy appear? It first appears in the construction of
\[
    L=(\ff\vee_{ij}\gg)\vee_{i'j'}\hh,
    \qquad R=\ff\vee_{ij}(\gg\vee_{i'j'}\hh).
\]
The construction extends in the obvious way from a single pushout, at the outer gluing, we treat the inner pushout as one category and its associated functor as one functor.
Indeed, go back to the example $n=m=\ell=1$, denoting the arrows in the three factors by $a$'s, $b$'s and $c$'s, respectively. For $R=\ff\vee_{10}(\gg\vee_{10}\hh)$, we build $\gg\vee_{10}\hh$ first and think of it as one functor on one category. On the sequence $(\mu_2(a_1,b_1),b_0,\mu_2(b_1,c_1))$, one finds
\begin{align*}
    &R_3(\mu_2(a_1,b_1),b_0,\mu_2(b_1,c_1))\\
    &\quad=\mu_2\bigl(\ff_1(a_1),(\gg\vee_{10}\hh)_3(b_1,b_0,\mu_2(b_1,c_1))\bigr)\\
    &\quad=\mu_2\bigl(\ff_1(a_1),\mu_2(\gg_3(b_1,b_0,b_1),\hh_1(c_1))\bigr).
\end{align*}
Whereas, for $L=(\ff\vee_{10}\gg)\vee_{10}\hh$, we first think of $\ff\vee_{10}\gg$ as one functor and obtain
\begin{align*}
    &L_3(\mu_2(a_1,b_1),b_0,\mu_2(b_1,c_1))\\
    &\quad=\mu_2\bigl((\ff\vee_{10}\gg)_3(\mu_2(a_1,b_1),b_0,b_1),\hh_1(c_1)\bigr)\\
    &\quad=\mu_2\bigl(\mu_2(\ff_1(a_1),\gg_3(b_1,b_0,b_1)),\hh_1(c_1)\bigr).
\end{align*}
Notice then that
\begin{align*}
    &(R-L)_3(\mu_2(a_1,b_1),b_0,\mu_2(b_1,c_1))\\
    &\quad=\mu_1\mu_3\bigl(\ff_1(a_1),\gg_3(b_1,b_0,b_1),\hh_1(c_1)\bigr)\\
    &\qquad+\mu_3\bigl(\ff_1(a_1),\mu_1\gg_3(b_1,b_0,b_1),\hh_1(c_1)\bigr).
\end{align*}
The point is that the difference $R-L$ need not be zero, but we will show below that it is a boundary with respect to $\mu_1^{\mathrm{Fun}}$.
We also note that the two orders can agree on particular inputs. For example,
\begin{align*}
    R_1\bigl(\mu_3(a_1,b_1,c_1)\bigr)
    &=\mu_3\bigl(\ff_1(a_1),(\gg\vee_{10}\hh)_1(b_1),(\gg\vee_{10}\hh)_1(c_1)\bigr)\\
    &=\mu_3\bigl(\ff_1(a_1),\gg_1(b_1),\hh_1(c_1)\bigr),
\end{align*}
and similarly,
\begin{align*}
    L_1\bigl(\mu_3(a_1,b_1,c_1)\bigr)
    &=\mu_3\bigl((\ff\vee_{10}\gg)_1(a_1),(\ff\vee_{10}\gg)_1(b_1),\hh_1(c_1)\bigr)\\
    &=\mu_3\bigl(\ff_1(a_1),\gg_1(b_1),\hh_1(c_1)\bigr).
\end{align*}

\begin{lemma}
Let
\[
    L=(\ff\vee_{ij}\gg)\vee_{i'j'}\hh,
    \qquad
    R=\ff\vee_{ij}(\gg\vee_{i'j'}\hh)
\]
be the chosen strictly unital $A_\infty$-functors from
$Q_n\bigvee_{ij}Q_m\bigvee_{i'j'}Q_\ell$ to $\cc$.
Then $L$ and $R$ are homotopic. More precisely, there exists a
pre-natural transformation
\[
    \Psi\in
    \operatorname{Fun}_{A_\infty}\!
    \bigl(Q_n\bigvee_{ij}Q_m\bigvee_{i'j'}Q_\ell,\cc\bigr)(L,R)_0
\]
with $\Psi_0=0$, vanishing on both two-factor pushouts, such that
\[
    \mu_1^{\mathrm{Fun}}\Psi=R-L.
\]

\end{lemma}

\begin{rmk}
We identify the two double pushouts using the same reduced trees and
the same object shifts in both orders. Thus we may use the same
attachment indices in the two expressions. The parentheses keep track of which order is taken.
\end{rmk}

\begin{proof}
Write $D_d=R_d-L_d=R_d+L_d$. We will take $\Psi_k=0$ for $k\neq1$ and construct $\Psi_1$ recursively. With this choice, the equation to be satisfied, for $d\geq1$, is
\begin{equation}\label{whatpsishouldsat}
\begin{aligned}
    D_d(x_1,\ldots,x_d)
    &=\Psi_1\bigl(\mu_d(x_1,\ldots,x_d)\bigr)\\
    &\quad+\sum\mu_{\alpha+1+\beta}\Bigl(
      L_{p_1}(x_1,\ldots,x_{p_1}),\ldots,
      L_{p_\alpha}(x_{v-p_\alpha},\ldots,x_{v-1}),\Psi_1(x_v),\\
    &\qquad R_{t_1}(x_{v+1},\ldots,x_{v+t_1}),\ldots,
      R_{t_\beta}(x_{d-t_\beta+1},\ldots,x_d)\Bigr).
\end{aligned}
\end{equation}
The above is simply the definition of a homotopy of functors with the additional assumption $\Psi_k=0$ for $k\geq 2$ (see \cite{Or} for the definition of homotopy here). First notice that $L$ and $R$ agree on every sequence contained entirely in
$Q_n\bigvee_{ij}Q_m$. Indeed, on such a sequence every subsequence evaluated
by $\gg\vee_{i'j'}\hh$ lies in $Q_m$, where this functor restricts
to $\gg$. Thus both orders give $\ff\vee_{ij}\gg$. Similarly,
both orders restrict to $\gg\vee_{i'j'}\hh$ on
$Q_m\bigvee_{i'j'}Q_\ell$. We set $\Psi_1=0$ on these two
subcategories.

Now suppose that $x_1,\ldots,x_k$ are individual $a$-, $b$-, or
$c$-arrows. For $k\geq2$ we have
\[
    D_k(x_1,\ldots,x_k)=0,
    \qquad
    D_1\bigl(\mu_k(x_1,\ldots,x_k)\bigr)=0.
\]
For a sequence using all three factors, these equalities follow directly from the definition of $(-)\vee_{ij}(-)$.  We set $\Psi_1(x_v)=0$ on individual arrows, so equation~\eqref{whatpsishouldsat} gives
\[
    \Psi_1\bigl(\mu_k(x_1,\ldots,x_k)\bigr)=0.
\]

We next insert a single operation $T$ into the sequence.
The first possible nonzero value of $\Psi_1$ is obtained by
applying one further operation to the resulting tuple.
To understand this, consider the sequence
\[
    (a_{i-l},\ldots,a_i,\,
     b_{j+1},\ldots,b_{j+B},\,
     c_{j'+1},\ldots,c_{j'+s}),
\]
its length is $k=l+1+B+s$. First suppose that $T$ crosses only the first gluing
point.
Take $T=\mu_{\alpha+q}$ and set
\[
    u=T(a_{r+1},\ldots,a_i,b_{j+1},\ldots,b_{j+q}),
\]
Writing $d=k-\alpha-q+1$, we can set
\[
    \bar D
    :=D_d(a_{i-l},\ldots,a_r,u,
          b_{j+q+1},\ldots,b_{j+B},
          c_{j'+1},\ldots,c_{j'+s}).
\]
For $L$, the input $u$ already belongs to the first pushout
$Q_n\bigvee_{ij}Q_m$. The remaining $c$-arrows are separate inputs
of the other factor, so $L_d=0$. Thus $\bar D$
is given by $R_d$, which we can calculate to give
\begin{align*}
    \bar D
    &=
    \sum_{e=1}^{q-1}
    \sum_{\substack{
      p_1+\cdots+p_\alpha=l+1\\
      t_1+\cdots+t_e=B\\
      t_{e+1}+\cdots+t_q=s\\
      p_1,\ldots,p_\alpha,t_1,\ldots,t_q\geq1}}
    T\Bigl(
      \ff_{p_1}(a_{i-l},\ldots,a_{i-l+p_1-1}),
      \ldots,
      \ff_{p_\alpha}(a_{i-p_\alpha+1},\ldots,a_i),\\
    &\qquad
      \gg_{t_1}(b_{j+1},\ldots,b_{j+t_1}),
      \ldots,
      \gg_{t_e}(b_{j+B-t_e+1},\ldots,b_{j+B}),\\
    &\qquad
      \hh_{t_{e+1}}(c_{j'+1},\ldots,c_{j'+t_{e+1}}),
      \ldots,
      \hh_{t_q}(c_{j'+s-t_q+1},\ldots,c_{j'+s})
      \Bigr).
\end{align*}
Every input of the tuple has zero
$\Psi_1$, including $u$, so equation~\eqref{whatpsishouldsat}
gives
\[
    \Psi_1\Bigl(\mu_d\bigl(
      a_{i-l},\ldots,a_r,u,
      b_{j+q+1},\ldots,b_{j+B},
      c_{j'+1},\ldots,c_{j'+s}\bigr)\Bigr):=\bar D.
\]
If $T$ crosses only the second gluing point, put
\[
    u=T(b_{j+h+1},\ldots,b_{j+B},c_{j'+1},\ldots,c_{j'+v}).
\]
Here $d=k-B+h-v+1$. Now $u$ belongs to $Q_m\bigvee_{i'j'}Q_\ell$, so
\[
    R_d(a_{i-l},\ldots,a_i,b_{j+1},\ldots,b_{j+h},u,c_{j'+v+1},\ldots,c_{j'+s})=0.
\]
Also, $\Psi_1(u)=0$, so we can set
\begin{align}
    &\Psi_1\bigl(\mu_d(a_{i-l},\ldots,a_i,b_{j+1},\ldots,b_{j+h},
      u,c_{j'+v+1},\ldots,c_{j'+s})\bigr)\notag\\
    &\quad:=D_d(a_{i-l},\ldots,a_i,b_{j+1},\ldots,b_{j+h},
      u,c_{j'+v+1},\ldots,c_{j'+s}).
\end{align}

Finally, suppose that $T$ crosses both gluing points and put
\[
u=T(a_{r+1},\ldots,a_i,b_{j+1},\ldots,b_{j+B},c_{j'+1},\ldots,c_{j'+v}).
\]
Here $d=k-\alpha-B-v+1$. For $d\geq2$, both orders may contribute to $D_d$. However, $T$ is still a single operation on individual arrows and we have already set $\Psi_1(u)=0$, so we similarly obtain
\begin{align}
    &\Psi_1\bigl(\mu_d(a_{i-l},\ldots,a_r,u,c_{j'+v+1},\ldots,c_{j'+s})\bigr)\notag\\
    &\quad:=D_d(a_{i-l},\ldots,a_r,u,c_{j'+v+1},\ldots,c_{j'+s}).
\end{align}

The same assignment works for any composable sequence with one
inserted operation, regardless of the number of jumps between the three factors. Each input
is then either an arrow or a single operation on arrows, so every
insertion term involving $\Psi_1(x_v)$ in
equation~\eqref{whatpsishouldsat} vanishes.

We then induct on the total number of internal vertices, counting each operation $\mu_d$, $d\geq2$, as one vertex.
Suppose that $\Psi_1$ has been defined on all trees with at most $k$ internal vertices. A tree with $k+1$ internal vertices is obtained by applying $\mu_d$, with $d\geq2$, to a tuple
\[
    (y_1,\ldots,y_d)
    =\bigl(T_1(x_1,\ldots,x_{s_1}),\ldots,
    T_d(x_{s_{d-1}+1},\ldots,x_{s_d})\bigr),
\]
where the $x_i$ are individual arrows and the $T_i$ have $k$ internal vertices in total. We allow a $T_i$ with no internal vertices, so that $y_i$ is simply an individual arrow. Each $\Psi_1(y_i)$ is therefore known by induction. Using equation~\eqref{whatpsishouldsat}, we set
\begin{align*}
    \Psi_1\bigl(\mu_d(y_1,\ldots,y_d)\bigr)
    &:=D_d(y_1,\ldots,y_d)\\
    &\quad+\sum_{v=1}^{d}\sum\mu_{\alpha+1+\beta}\Bigl(
      L_{p_1}(y_1,\ldots,y_{p_1}),\ldots,
      L_{p_\alpha}(y_{v-p_\alpha},\ldots,y_{v-1}),\\
    &\qquad\Psi_1(y_v),
      R_{t_1}(y_{v+1},\ldots,y_{v+t_1}),\ldots,
      R_{t_\beta}(y_{d-t_\beta+1},\ldots,y_d)\Bigr).
\end{align*}
The inner sum is over the same partitions as before, with $p_1+\cdots+p_\alpha=v-1$ and $t_1+\cdots+t_\beta=d-v$. Terms with $y_v$ an individual arrow vanish since $\Psi_1(y_v)=0$.

It remains to check the equation
\[
    \mu_1\Psi_1(x)+\Psi_1(\mu_1x)=D_1(x).
\]

To do so, we again induct on the number of vertices of $T$. To begin, consider a tree with one internal vertex,
\[
    x=\mu_k(x_1,\dots,x_k),
\]
where $k\geq2$ and each $x_i$ is an individual $a$-, $b$-, or
$c$-arrow. We have already set $\Psi_1(x)=0$, and know that $D_1(x)=0$, so it suffices to show
that $\Psi_1(\mu_1x)=0$. The $A_\infty$-relations give
\begin{align*}
    &\Psi_1\bigl(\mu_1\mu_k(x_1,\dots,x_k)\bigr)\\
    &\quad=
    \sum_{r=2}^{k-1}\sum_{s=0}^{k-r}
    \Psi_1\Bigl(\mu_{k-r+1}\bigl(
      x_1,\dots,x_s,
      \mu_r(x_{s+1},\dots,x_{s+r}),
      x_{s+r+1},\dots,x_k\bigr)\Bigr)\\
    &\quad=
    \sum_{r=2}^{k-1}\sum_{s=0}^{k-r}
    D_{k-r+1}\bigl(
      x_1,\dots,x_s,
      \mu_r(x_{s+1},\dots,x_{s+r}),
      x_{s+r+1},\dots,x_k\bigr)
\end{align*}
Now, since $D=R+L$ is defined componentwise, adding the two functor equations gives
\begin{align*}
    &\sum_{r=2}^{k-1}\sum_{s=0}^{k-r}
    D_{k-r+1}\bigl(x_1,\dots,x_s,\mu_r(x_{s+1},\dots,x_{s+r}),
    x_{s+r+1},\dots,x_k\bigr)\\
    &\quad=\sum_{t=2}^{k}\sum_{\substack{q_1+\cdots+q_t=k\\q_1,\dots,q_t\geq1}}
    \Bigl(\mu_t\bigl(R_{q_1}(x_1,\dots,x_{q_1}),\dots,
    R_{q_t}(x_{k-q_t+1},\dots,x_k)\bigr)\\
    &\qquad\qquad+\mu_t\bigl(L_{q_1}(x_1,\dots,x_{q_1}),\dots,
    L_{q_t}(x_{k-q_t+1},\dots,x_k)\bigr)\Bigr)=0.
\end{align*}
The two terms agree on each partition because $L_{q_j}=R_{q_j}$ on each block of individual arrows.
The combined sum is precisely $\Psi_1(\mu_1x)$, so this vanishes. Since $\Psi_1(x)=0$ as well, we obtain
\[
    \mu_1\Psi_1(x)+\Psi_1(\mu_1x)=0=D_1(x).
\]
Now we induct again. Assume the equation holds for all trees with at most $k$ internal vertices, and put
\[
    z=\mu_d(T_1,\ldots,T_d).
\]
Here $d\geq2$, each $T_i$ has $k_i\leq k$ internal vertices, and $\sum_{i=1}^{d}k_i=k$. As before, a subtree with no internal vertices is an individual arrow. Then $z$ has $k+1$ internal vertices, and for each $i$ by assumption we have
\[
    \mu_1\Psi_1(T_i)+\Psi_1(\mu_1T_i)=D_1(T_i).
\]
Set $y_i=T_i$. Then we find
\begin{align*}
    \mu_1\Psi_1(z)
    &=\mu_1D_d(y_1,\ldots,y_d)\\
    &\quad+\sum_{v=1}^{d}\sum\mu_1\mu_{\alpha+1+\beta}\Bigl(
      L_{p_1}(y_1,\ldots,y_{p_1}),\ldots,
      L_{p_\alpha}(y_{v-p_\alpha},\ldots,y_{v-1}),\\
    &\qquad\Psi_1(y_v),
      R_{t_1}(y_{v+1},\ldots,y_{v+t_1}),\ldots,
      R_{t_\beta}(y_{d-t_\beta+1},\ldots,y_d)\Bigr).
\end{align*}
We can expand the $\mu_1D_d$ term via the functor relations to obtain
\begin{align*}
    \mu_1D_d(y_1,\ldots,y_d)
    &=D_1(z)+\sum_{h=1}^{d}
    D_d(y_1,\ldots,y_{h-1},\mu_1y_h,y_{h+1},\ldots,y_d)\\
    &\quad+\sum_{\gamma=2}^{d-1}\sum_{o=0}^{d-\gamma}
    D_{d-\gamma+1}\Bigl(y_1,\ldots,y_o,
    \mu_\gamma(y_{o+1},\ldots,y_{o+\gamma}),y_{o+\gamma+1},\ldots,y_d\Bigr)\\
    &\quad+\sum_{q=2}^{d}\sum_{\substack{p_1+\cdots+p_q=d\\p_1,\ldots,p_q\geq1}}
    \sum_{v=1}^{q}\mu_q\Bigl(
    L_{p_1}(y_1,\ldots,y_{p_1}),\ldots,
    L_{p_{v-1}}(y_{o-p_{v-1}+1},\ldots,y_o),\\
    &\qquad D_{p_v}(y_{o+1},\ldots,y_{o+p_v}),\\
    &\qquad R_{p_{v+1}}(y_{o+p_v+1},\ldots,y_{o+p_v+p_{v+1}}),\ldots,
    R_{p_q}(y_{d-p_q+1},\ldots,y_d)\Bigr).
\end{align*}
Here $o=p_1+\cdots+p_{v-1}$ in the last sum. Next we expand the remaining sum, omitting the inputs for $L_{p_i}$ and $R_{t_j}$ to avoid too much notation:
\begin{align*}
    &\mu_1\mu_{\alpha+1+\beta}
    \bigl(L_{p_1},\ldots,L_{p_\alpha},\Psi_1(y_v),R_{t_1},\ldots,R_{t_\beta}\bigr)\\
    &\quad=\sum_{1\leq i\leq j\leq\alpha}
    \mu_{\alpha+1+\beta+i-j}\Bigl(L_{p_1},\ldots,L_{p_{i-1}},\\
    &\qquad\mu_{j-i+1}(L_{p_i},\ldots,L_{p_j}),
    L_{p_{j+1}},\ldots,L_{p_\alpha},\Psi_1(y_v),R_{t_1},\ldots,R_{t_\beta}\Bigr)\\
    &\qquad+\sum_{\substack{0\leq\gamma\leq\alpha,\;0\leq\delta\leq\beta\\
    \gamma+\delta<\alpha+\beta}}
    \mu_{\alpha+1+\beta-\gamma-\delta}\Bigl(L_{p_1},\ldots,L_{p_{\alpha-\gamma}},\\
    &\qquad\quad\mu_{\gamma+1+\delta}\bigl(
    L_{p_{\alpha-\gamma+1}},\ldots,L_{p_\alpha},\Psi_1(y_v),R_{t_1},\ldots,R_{t_\delta}\bigr),
    R_{t_{\delta+1}},\ldots,R_{t_\beta}\Bigr)\\
    &\qquad+\sum_{1\leq i\leq j\leq\beta}
    \mu_{\alpha+1+\beta+i-j}\Bigl(
    L_{p_1},\ldots,L_{p_\alpha},\Psi_1(y_v),R_{t_1},\ldots,R_{t_{i-1}},\\
    &\qquad\quad\mu_{j-i+1}(R_{t_i},\ldots,R_{t_j}),R_{t_{j+1}},\ldots,R_{t_\beta}\Bigr).
\end{align*}
The cases $i=j$ include the differentials of the $L$- and $R$-terms, while $\gamma=\delta=0$ gives the term containing $\mu_1\Psi_1(y_v)$.
Putting these together, we have
\begin{align*}
    \mu_1\Psi_1(z)
    &=D_1(z)+\sum_{h=1}^{d}
    D_d(y_1,\ldots,y_{h-1},\mu_1y_h,y_{h+1},\ldots,y_d)\\
    &\quad+\sum_{\gamma=2}^{d-1}\sum_{o=0}^{d-\gamma}
    D_{d-\gamma+1}\Bigl(y_1,\ldots,y_o,
    \mu_\gamma(y_{o+1},\ldots,y_{o+\gamma}),y_{o+\gamma+1},\ldots,y_d\Bigr)\\
    &\quad+\sum_{q=2}^{d}\sum_{\substack{p_1+\cdots+p_q=d\\p_1,\ldots,p_q\geq1}}
    \sum_{v=1}^{q}\mu_q\Bigl(
    L_{p_1}(y_1,\ldots,y_{p_1}),\ldots,
    L_{p_{v-1}}(y_{o-p_{v-1}+1},\ldots,y_o),\\
    &\qquad D_{p_v}(y_{o+1},\ldots,y_{o+p_v}),\\
    &\qquad R_{p_{v+1}}(y_{o+p_v+1},\ldots,y_{o+p_v+p_{v+1}}),\ldots,
    R_{p_q}(y_{d-p_q+1},\ldots,y_d)\Bigr)\\
    &\quad+\sum_{v=1}^{d}
    \sum_{\substack{\alpha,\beta\geq0\\p_1+\cdots+p_\alpha=v-1\\t_1+\cdots+t_\beta=d-v}}
    \sum_{1\leq i\leq j\leq\alpha}
    \mu_{\alpha+1+\beta+i-j}\Bigl(L_{p_1},\ldots,L_{p_{i-1}},\\
    &\qquad\mu_{j-i+1}(L_{p_i},\ldots,L_{p_j}),
    L_{p_{j+1}},\ldots,L_{p_\alpha},\Psi_1(y_v),R_{t_1},\ldots,R_{t_\beta}\Bigr)\\
    &\quad+\sum_{v=1}^{d}
    \sum_{\substack{\alpha,\beta\geq0\\p_1+\cdots+p_\alpha=v-1\\t_1+\cdots+t_\beta=d-v}}
    \sum_{\substack{0\leq\gamma\leq\alpha,\;0\leq\delta\leq\beta\\
    \gamma+\delta<\alpha+\beta}}
    \mu_{\alpha+1+\beta-\gamma-\delta}\Bigl(L_{p_1},\ldots,L_{p_{\alpha-\gamma}},\\
    &\qquad\mu_{\gamma+1+\delta}\bigl(
    L_{p_{\alpha-\gamma+1}},\ldots,L_{p_\alpha},\Psi_1(y_v),R_{t_1},\ldots,R_{t_\delta}\bigr),
    R_{t_{\delta+1}},\ldots,R_{t_\beta}\Bigr)\\
    &\quad+\sum_{v=1}^{d}
    \sum_{\substack{\alpha,\beta\geq0\\p_1+\cdots+p_\alpha=v-1\\t_1+\cdots+t_\beta=d-v}}
    \sum_{1\leq i\leq j\leq\beta}
    \mu_{\alpha+1+\beta+i-j}\Bigl(
    L_{p_1},\ldots,L_{p_\alpha},\Psi_1(y_v),R_{t_1},\ldots,R_{t_{i-1}},\\
    &\qquad\mu_{j-i+1}(R_{t_i},\ldots,R_{t_j}),R_{t_{j+1}},\ldots,R_{t_\beta}\Bigr).
\end{align*}
Notice there are seven types of summations. 
Next we consider $\Psi_1(\mu_1z)$:
\begin{align*}
    \Psi_1(\mu_1z)
    &=\sum_{h=1}^{d}\Psi_1\Bigl(
    \mu_d(y_1,\ldots,y_{h-1},\mu_1y_h,y_{h+1},\ldots,y_d)\Bigr)\\
    &\quad+\sum_{e=2}^{d-1}\sum_{j=0}^{d-e}
    \Psi_1\Bigl(\mu_{d-e+1}\bigl(y_1,\ldots,y_j,
    \mu_e(y_{j+1},\ldots,y_{j+e}),y_{j+e+1},\ldots,y_d\bigr)\Bigr).
\end{align*}
Here we keep the terms containing $\mu_1y_h$ separate. For each $e,j$, write
\[
    (y'_1,\ldots,y'_{d-e+1})
    =\bigl(y_1,\ldots,y_j,\mu_e(y_{j+1},\ldots,y_{j+e}),y_{j+e+1},\ldots,y_d\bigr).
\]
Then we can expand each term of the second summation:
\begin{align*}
    &\Psi_1\bigl(\mu_{d-e+1}(y'_1,\ldots,y'_{d-e+1})\bigr)\\
    &\quad=D_{d-e+1}(y'_1,\ldots,y'_{d-e+1})\\
    &\qquad+\sum\mu_{\alpha+1+\beta}\Bigl(
    L_{p_1}(y_1,\ldots,y_{p_1}),\ldots,
    L_{p_\alpha}(y_{j-p_\alpha+1},\ldots,y_j),\\
    &\qquad\quad\Psi_1\bigl(\mu_e(y_{j+1},\ldots,y_{j+e})\bigr),
    R_{t_1}(y_{j+e+1},\ldots,y_{j+e+t_1}),\ldots,
    R_{t_\beta}(y_{d-t_\beta+1},\ldots,y_d)\Bigr)\\
    &\qquad+\sum_{\substack{1\leq v\leq d-e+1\\v\neq j+1}}\sum
    \mu_{\alpha+1+\beta}\Bigl(
    L_{p_1}(y'_1,\ldots,y'_{p_1}),\ldots,
    L_{p_\alpha}(y'_{v-p_\alpha},\ldots,y'_{v-1}),\\
    &\qquad\quad\Psi_1(y'_v),
    R_{t_1}(y'_{v+1},\ldots,y'_{v+t_1}),\ldots,
    R_{t_\beta}(y'_{d-e-t_\beta+2},\ldots,y'_{d-e+1})\Bigr).
\end{align*}

For the first remaining terms, we have
\begin{align*}
    &\Psi_1\bigl(\mu_e(y_{j+1},\ldots,y_{j+e})\bigr)\\
    &\quad=D_e(y_{j+1},\ldots,y_{j+e})\\
    &\qquad+\sum_{h=j+1}^{j+e}\sum\mu_{\gamma+1+\delta}\Bigl(
    L_{p'_1}(y_{j+1},\ldots,y_{j+p'_1}),\ldots,
    L_{p'_\gamma}(y_{h-p'_\gamma},\ldots,y_{h-1}),\\
    &\qquad\quad\Psi_1(y_h),
    R_{t'_1}(y_{h+1},\ldots,y_{h+t'_1}),\ldots,
    R_{t'_\delta}(y_{j+e-t'_\delta+1},\ldots,y_{j+e})\Bigr).
\end{align*}
And so putting all together we have
\begin{align*}
    \Psi_1(\mu_1z)
    &=\sum_{h=1}^{d}\Psi_1\Bigl(
    \mu_d(y_1,\ldots,y_{h-1},\mu_1y_h,y_{h+1},\ldots,y_d)\Bigr)\\
    &\quad+\sum_{e=2}^{d-1}\sum_{j=0}^{d-e}
    D_{d-e+1}\Bigl(y_1,\ldots,y_j,
    \mu_e(y_{j+1},\ldots,y_{j+e}),y_{j+e+1},\ldots,y_d\Bigr)\\
    &\quad+\sum_{e=2}^{d-1}\sum_{j=0}^{d-e}\sum
    \mu_{\alpha+1+\beta}\Biggl(
    L_{p_1}(y_1,\ldots,y_{p_1}),\ldots,
    L_{p_\alpha}(y_{j-p_\alpha+1},\ldots,y_j),\\
    &\qquad\Bigl[D_e(y_{j+1},\ldots,y_{j+e})\\
    &\qquad\quad+\sum_{h=j+1}^{j+e}\sum\mu_{\gamma+1+\delta}\Bigl(
    L_{p'_1}(y_{j+1},\ldots,y_{j+p'_1}),\ldots,
    L_{p'_\gamma}(y_{h-p'_\gamma},\ldots,y_{h-1}),\\
    &\qquad\qquad\Psi_1(y_h),
    R_{t'_1}(y_{h+1},\ldots,y_{h+t'_1}),\ldots,
    R_{t'_\delta}(y_{j+e-t'_\delta+1},\ldots,y_{j+e})\Bigr)\Bigr],\\
    &\qquad R_{t_1}(y_{j+e+1},\ldots,y_{j+e+t_1}),\ldots,
    R_{t_\beta}(y_{d-t_\beta+1},\ldots,y_d)\Biggr)\\
    &\quad+\sum_{e=2}^{d-1}\sum_{j=0}^{d-e}
    \sum_{\substack{1\leq v\leq d-e+1\\v\neq j+1}}\sum
    \mu_{\alpha+1+\beta}\Bigl(
    L_{p_1}(y'_1,\ldots,y'_{p_1}),\ldots,
    L_{p_\alpha}(y'_{v-p_\alpha},\ldots,y'_{v-1}),\\
    &\qquad\Psi_1(y'_v),
    R_{t_1}(y'_{v+1},\ldots,y'_{v+t_1}),\ldots,
    R_{t_\beta}(y'_{d-e-t_\beta+2},\ldots,y'_{d-e+1})\Bigr).
\end{align*}
In the last sum, $y'_v$ is an input outside the inner operation, and the $L$- and $R$-blocks are taken before and after it in the tuple $(y'_1,\ldots,y'_{d-e+1})$ defined above.

For $1\leq h\leq d$, write $y_h^h=\mu_1y_h$ and $y_i^h=y_i$ for $i\neq h$. We then have
\begin{align*}
    &\Psi_1\Bigl(\mu_d(y_1,\ldots,y_{h-1},\mu_1y_h,y_{h+1},\ldots,y_d)\Bigr)\\
    &\quad=D_d(y_1,\ldots,y_{h-1},\mu_1y_h,y_{h+1},\ldots,y_d)\\
    &\qquad+\sum\mu_{\alpha+1+\beta}\Bigl(
    L_{p_1}(y_1,\ldots,y_{p_1}),\ldots,
    L_{p_\alpha}(y_{h-p_\alpha},\ldots,y_{h-1}),\\
    &\qquad\quad\Psi_1(\mu_1y_h),
    R_{t_1}(y_{h+1},\ldots,y_{h+t_1}),\ldots,
    R_{t_\beta}(y_{d-t_\beta+1},\ldots,y_d)\Bigr)\\
    &\qquad+\sum_{\substack{1\leq v\leq d\\v\neq h}}\sum
    \mu_{\alpha+1+\beta}\Bigl(
    L_{p_1}(y_1^h,\ldots,y_{p_1}^h),\ldots,
    L_{p_\alpha}(y_{v-p_\alpha}^h,\ldots,y_{v-1}^h),\\
    &\qquad\quad\Psi_1(y_v),
    R_{t_1}(y_{v+1}^h,\ldots,y_{v+t_1}^h),\ldots,
    R_{t_\beta}(y_{d-t_\beta+1}^h,\ldots,y_d^h)\Bigr).
\end{align*}
One can carry out the explicit calculations and verify that
\begin{align}
    \mu_1\Psi_1(z)+\Psi_1(\mu_1z)
    &=D_1(z)\notag\\
    &\quad+\sum_{v=1}^{d}\sum
    \mu_{\alpha+1+\beta}\Bigl(L_{p_1},\ldots,L_{p_\alpha},\notag\\
    &\qquad\bigl[D_1(y_v)+\mu_1\Psi_1(y_v)+\Psi_1(\mu_1y_v)\bigr],
    R_{t_1},\ldots,R_{t_\beta}\Bigr).
\end{align}
But we know that $D_1(y_v)+\mu_1\Psi_1(y_v)+\Psi_1(\mu_1y_v)=0$ by the inductive assumption, so we recover
\begin{equation}
    \mu_1\Psi_1(z)+\Psi_1(\mu_1z)=D_1(z).
\end{equation}

\end{proof}

\subsection{Parallel Gluings}
Okay, so the two orders of taking the pushout give two different but homotopic functors, and we can build an explicit homotopy $\Psi$ controlling the difference. Next we consider the gluing functors. Consider $Q_n\bigvee_{ij}Q_m\bigvee_{i'j'}Q_\ell$. Recall that the category was built in two stages: one order starts by pushing out over $Q_n$ and $Q_m$, and the other by pushing out over $Q_m$ and $Q_\ell$. Here we assume $n,m,\ell\geq1$ and $j\neq i'$, and use the same final object labels and shifts in both orders.
After the first gluing, write $\widetilde{i'}$ and $\widetilde j$ for the labels of the remaining attachment points, so that
\[
    (\Phi_{ij})_0(\widetilde{i'})=i'_{Q_m},
    \qquad
    (\Phi_{i'j'})_0(\widetilde j)=j_{Q_m}.
\]
From this we can consider the following two compositions, with common source $Q_{n+m+\ell-2}=Q_{(n+m-1)+\ell-1}=Q_{n+(m+\ell-1)-1}$:
\begin{equation}
    \begin{tikzcd}[column sep=large]
     Q_{n+m+\ell-2} \ar[r,"\Phi_{\widetilde{i'}j'}"]
     & Q_{n+m-1}\bigvee_{\widetilde{i'}j'}Q_\ell
       \ar[r,"\Phi_{ij}\vee 1_{Q_\ell}"]
     & Q_n\bigvee_{ij}Q_m\bigvee_{i'j'}Q_\ell
    \end{tikzcd}
\end{equation}
\begin{equation}
    \begin{tikzcd}[column sep=large]
     Q_{n+m+\ell-2} \ar[r,"\Phi_{i\widetilde j}"]
     & Q_n\bigvee_{i\widetilde j}Q_{m+\ell-1}
       \ar[r,"1_{Q_n}\vee\Phi_{i'j'}"]
     & Q_n\bigvee_{ij}Q_m\bigvee_{i'j'}Q_\ell
    \end{tikzcd}
\end{equation}
For simplicity, let's write
\begin{align}
    \Phi^L&:=(\Phi_{ij}\vee 1_{Q_\ell})\circ\Phi_{\widetilde{i'}j'},
    &\Phi^R&:=(1_{Q_n}\vee\Phi_{i'j'})\circ\Phi_{i\widetilde j}.
\end{align}
\begin{prop}
    The strictly unital $A_\infty$-functors $\Phi^L$ and $\Phi^R$ are homotopic. In particular, there exists a pre-natural transformation
    \begin{equation}
        \Xi\in \operatorname{Fun}_{A_\infty}\bigl(Q_{n+m+\ell-2},Q_n\bigvee_{ij}Q_m\bigvee_{i'j'}Q_\ell\bigr)(\Phi^L,\Phi^R)_0
    \end{equation}
    such that 
    \begin{equation}
        \mu_1^{\text{Fun}}\Xi=\Phi^L+\Phi^R,
    \end{equation}
    with $\Xi_0=0$.
\end{prop}
\begin{proof}
Similarly to the case of $\Psi$, we find that $\Phi^R$ and $\Phi^L$ agree on sequences whose image lies entirely in $Q_n$, $Q_m$ or $Q_\ell$. They also agree on paths whose image makes at most one jump (this is easy to verify). Thus, we take $\Xi_k=0$ on all such sequences.

Now consider a sequence in $Q_{n+m+\ell-2}$ whose image passes from $Q_n$ through $Q_m$ to $Q_\ell$. Write its expanded arrow sequence as
\[
    (a_{i-p+1},\ldots,a_i,\ b_{j+1},\ldots,b_{j+t},\ c_{j'+1},\ldots,c_{j'+q}).
\]
Write $d_r$ for the arrows in the domain, we find
   \begin{align*}
 &\Phi^L_k(d_{i-p+1},\ldots,d_{\widetilde{i'}+q-1})\\
 &\quad=\sum_{r=0}^{p+t-2}\sum_{s=0}^{q-1}
   (\Phi_{ij}\vee 1_{Q_\ell})_{r+s+1}\Bigl(
   \overline{ab}_{i-p+1},\ldots,\overline{ab}_{i-p+r},\\
 &\qquad\qquad
   \mu_{k+1-r-s}(\overline{ab}_{i-p+r+1},\ldots,\overline{ab}_{\widetilde{i'}},
                        c_{j'+1},\ldots,c_{j'+q-s}),
   c_{j'+q-s+1},\ldots,c_{j'+q}\Bigr).
\end{align*}
where $\overline{ab}_r$ denotes an arrow of $Q_{n+m-1}$ before applying $\Phi_{ij}$. We then use $\Phi_{ij}$ and our definition of $(-)\vee(-)$ to find 
\begin{equation}\label{philefttwocrossings}
\begin{aligned}
 &\Phi^L_k(d_{i-p+1},\ldots,d_{\widetilde{i'}+q-1})\\
 &\quad=\sum_{u=1}^{p}\sum_{v=1}^{t}
   \mu_{p+t+q-u-v+1}\Bigl(a_{i-p+1},\ldots,a_{i-u},\\
 &\qquad\qquad
   \mu_{u+v}(a_{i-u+1},\ldots,a_i,b_{j+1},\ldots,b_{j+v}),
   b_{j+v+1},\ldots,b_{j+t},c_{j'+1},\ldots,c_{j'+q}\Bigr).
\end{aligned}
\end{equation}
One can repeat the same calculation for $\Phi^R$ and find that 
\begin{equation}\label{phirighttwocrossings}
\begin{aligned}
 &\Phi^R_k(d_{i-p+1},\ldots,d_{\widetilde{i'}+q-1})\\
 &\quad=\sum_{v=1}^{t}\sum_{w=1}^{q}
 \mu_{p+t+q-v-w+1}\Bigl(
 a_{i-p+1},\ldots,a_i,b_{j+1},\ldots,b_{j+t-v},\\
 &\qquad\qquad
 \mu_{v+w}(b_{j+t-v+1},\ldots,b_{j+t},c_{j'+1},\ldots,c_{j'+w}),
 c_{j'+w+1},\ldots,c_{j'+q}\Bigr).
\end{aligned}
\end{equation}
One will hopefully notice that these two sums fit into the homotopy equation below. Indeed, this is a direct consequence of the $A_\infty$-relations, with the additional realisation that 
\begin{equation*}
    \mu_d(x_1,\dots,x_d)=0
\end{equation*}
for a sequence $(x_1,\dots,x_d)$ with the image we are considering. Thus we set
\begin{equation}\label{xitwocrossings}
\begin{aligned}
 &\Xi_k(d_{i-p+1},\ldots,d_{\widetilde{i'}+q-1})\\
 &\quad:=\mu_{p+t+q}(a_{i-p+1},\ldots,a_i,\,
       b_{j+1},\ldots,b_{j+t},\,c_{j'+1},\ldots,c_{j'+q}).
\end{aligned}
\end{equation}
Indeed, summing equations~\eqref{philefttwocrossings} and~\eqref{phirighttwocrossings} we find
\[
\begin{aligned}
 &(\Phi^L_k+\Phi^R_k)(d_{i-p+1},\ldots,d_{\widetilde{i'}+q-1})\\
 &\quad=\mu_1\Xi_k(d_{i-p+1},\ldots,d_{\widetilde{i'}+q-1})\\
 &\qquad+\sum_{\substack{1\le u\le p,\;1\le w\le q\\(u,w)\ne(p,q)}}
    \mu_{p-u+1+q-w}\Bigl(a_{i-p+1},\ldots,a_{i-u},\\
 &\qquad\qquad
    \Xi_{u+t+w-2}(d_{i-u+1},\ldots,d_{\widetilde{i'}+w-1}),
    c_{j'+w+1},\ldots,c_{j'+q}\Bigr).
\end{aligned}
\]
Now one iterates this construction on the length of the input sequence. Assume that the components $\Xi_d$ have been defined and satisfy the homotopy equations for all $d<k$. To define $\Xi_k$, we rearrange the above equation to give
\[
\begin{aligned}
 &\mu_1\Xi_k(d_{i-p+1},\ldots,d_{\widetilde{i'}+q-1})\\
 &\quad=(\Phi^L_k+\Phi^R_k)(d_{i-p+1},\ldots,d_{\widetilde{i'}+q-1})\\
 &\qquad+\sum_{\substack{1\le u\le p,\;1\le w\le q\\(u,w)\ne(p,q)}}
    \mu_{p-u+1+q-w}\Bigl(a_{i-p+1},\ldots,a_{i-u},\\
 &\qquad\qquad
    \Xi_{u+t+w-2}(d_{i-u+1},\ldots,d_{\widetilde{i'}+w-1}),
    c_{j'+w+1},\ldots,c_{j'+q}\Bigr).
\end{aligned}
\]
On the right-hand side, every value is known, and the inductive assumption implies that the sum is a cycle. By the extension of Lemma~\ref{lemphihomology} to three factors, the only possible non-trivial homology classes are those of a unit or one of the $f_i$. For the sequence above, when $k\geq2$, the degree does not match that of any non-trivial homology representative. Thus the right-hand side is a boundary, and we can choose $\Xi_k$ to be a primitive. We treat $k=1$ below.

In general, we must consider the full homotopy equation. {We take $\Xi_s$ to vanish whenever one of its inputs is a unit. On $(d_{r+1},\ldots,d_{r+k})$, the equation is
\begin{align}\label{generalxiequation}
    &\mu_1\Xi_k(d_{r+1},\ldots,d_{r+k})\notag\\
    &\quad=(\Phi^L_k+\Phi^R_k)(d_{r+1},\ldots,d_{r+k})\notag\\
    &\qquad+\sum
    \mu_{\alpha+1+\beta}\bigl(
    \Phi^L_{p_1},\ldots,\Phi^L_{p_\alpha},\Xi_s,
    \Phi^R_{t_1},\ldots,\Phi^R_{t_\beta}\bigr).
\end{align}
noting that the RHS is already known, and also that we set $\Xi_1(1)=0$. }

The degree of the right-hand side on $(d_{r+1},\ldots,d_{r+k})$ is
\[
    k-1+\sum_{s=1}^{k}|d_{r+s}|.
\]
For a unit, the input sequence must have length $k=\alpha(n+m+\ell-1)$ for some $\alpha\geq1$. Since the arrows in a full loop in the domain have total degree $2-(n+m+\ell-1)$, we find
\[
    k-1+\sum_{s=1}^{k}|d_{r+s}|=2\alpha-1\neq0.
\]
Thus the degrees cannot match, and the cycle must be a boundary. In particular, for a single full loop, $k=n+m+\ell-1$, the right-hand side of~\eqref{generalxiequation} has degree $1$, while the unit has degree $0$. We can therefore choose $\Xi_{n+m+\ell-1}$ to be a primitive in this case as well. For a class represented by a morphism in $\bigl(Q_n\bigvee_{ij}Q_m\bigvee_{i'j'}Q_\ell\bigr)(\Phi^L_0(r),\Phi^L_0(r+1))$, the sequence must have length $k=\alpha(n+m+\ell-1)+1$, and
\[
    k-1+\sum_{s=1}^{k}|d_{r+s}|=2\alpha+|d_{r+1}|.
\]
Again, the degrees do not match when $\alpha\geq1$.
We are left with $\alpha=0$, so we just need to consider $\Xi_1(d_{r+1})$. The only time that $\Phi^R$ and $\Phi^L$ can differ on $d_{r+1}$ is when the copies of $Q_n$ and $Q_\ell$ are attached to $Q_m$ at consecutive points. In this case, writing $a,b,c$ for the three arrows in path order, we have
\begin{align*}
    (\Phi^L_1+\Phi^R_1)(d_{r+1})
    &=\mu_2(\mu_2(a,b),c)+\mu_2(a,\mu_2(b,c))\\
    &=\mu_1\mu_3(a,b,c).
\end{align*}
Thus we can set $\Xi_1(d_{r+1})=\mu_3(a,b,c)$. This completes the induction.
\end{proof}

\begin{cor}\label{cortwogluedfunctors}
    Given three functors $\ff:Q_n\to\cc$, $\gg:Q_m\to\cc$ and $\hh:Q_\ell\to\cc$, with $\ff_0(i)=\gg_0(j)$ and $\gg_0(i')=\hh_0(j')$, the two glued functors
    \begin{align}
        (\ff\#_{ij}\gg)\#_{\widetilde{i'}j'}\hh&:=(\Phi^L)^*L,\\
        \ff\#_{i\widetilde j}(\gg\#_{i'j'}\hh)&:=(\Phi^R)^*R
    \end{align}
    are homotopic.
\end{cor}
\begin{proof}
    We find
    \[
        (\Phi^L)^*L\simeq(\Phi^R)^*L\simeq(\Phi^R)^*R.
    \]
\end{proof}

The remainder of this section will be devoted to proving the following theorem. 
\begin{theorem}\label{thmparallelgluings}
    Let $\gg:Q_m\to\cc$ and consider a family of functors $\ff^I:=\{\ff^i:Q_{n_i}\to\cc\}_{i\in I\subseteq\{0,\dots,m-1\}}$, with $\ff^i_0(n_i)=\gg_0(i)$ and $m,n_i\geq1$. Then there exists a family of functors 
    \begin{equation}\Big\{\#^p(\ff^I;\gg):Q_{\sum_{i\in I}n_i +m-|I|}\to \cc\Big\}_{p\in \text{Ord}(I)}\end{equation}
   where $\text{Ord}(I)$ denotes the set of all ordered sequences $(i_1,\dots,i_{|I|})$ with $i_j\in I$ appearing exactly once. These, together with corresponding higher homotopies, satisfy the cubical homotopy relations. 

\end{theorem}

For each order $p$, we will choose an iterated pushout functor built using $(-)\vee(-)$, together with a functor from $Q_{\sum_{i\in I}n_i+m-|I|}$ into the same pushout, generalising $\Phi^L$ and $\Phi^R$. We define $\#^p(\ff^I;\gg)$ by pulling back the iterated pushout functor along this chosen functor.

The vertices of $[0,1]^I$ are the subsets $S\subseteq I$ of inputs
already attached. An edge $S\to S\cup\{i\}$ attaches $\ff^i$, and an
order $p$ gives a path from $\varnothing$ to $I$ (i.e., from $(0,\dots,0)$ to $(1,\dots,1)$).

Consider the case of two functors $\ff$ and $\hh$ attached to $\gg$ at
the distinct vertices $j$ and $i'$ of $Q_m$, as in
Corollary~\ref{cortwogluedfunctors}. The two coordinates of $[0,1]^2$
record whether $\ff$ and $\hh$, respectively, have been attached.
Starting at $(0,0)$ with $\gg$, we can follow either path in the
directed square
\begin{equation}\label{parallelgluingsquare}
\begin{tikzcd}[row sep=large,column sep=huge]
    {(0,1)=}\quad \gg\#_{i'j'}\hh
        \arrow[r,"\text{attach }\ff"]
    & {(1,1)=}\quad \text{both attached}\\
    {(0,0)=}\quad \gg
        \arrow[u,"\text{attach }\hh"]
        \arrow[r,"\text{attach }\ff"']
    & {(1,0)=}\quad \ff\#_{ij}\gg
        \arrow[u,"\text{attach }\hh"']
\end{tikzcd}
\end{equation}
The path through $(1,0)$ attaches $\ff$ first and corresponds to
$(\Phi^L)^*L=(\ff\#_{ij}\gg)\#_{\widetilde{i'}j'}\hh$. The path through $(0,1)$ attaches $\hh$ first and gives
 $(\Phi^R)^*R=\ff\#_{i\widetilde j}(\gg\#_{i'j'}\hh)$. Thus the
vertex $(1,1)$ records the same pair of attachments in both cases,
while the functor obtained there may depend on the path. A filling of
this square is a homotopy between the glued functors obtained along
its two paths. When one goes one order higher and considers $I=\{i_0,i_1,i_2\}$, we obtain a cube of attachments:  
\begin{equation}\label{parallelgluingcube}
\begin{tikzcd}[row sep=large,column sep=large]
 & {(0,1,1)} \arrow[rr,"\ff^{i_0}"] & & {(1,1,1)} \\
 {(0,1,0)} \arrow[rr,"\ff^{i_0}" near end] \arrow[ur,"\ff^{i_2}"]
 & & {(1,1,0)} \arrow[ur,"\ff^{i_2}"'] & \\
 & {(0,0,1)} \arrow[uu,dashed,"\ff^{i_1}" near start]
             \arrow[rr,dashed,"\ff^{i_0}" near start]
 & & {(1,0,1)} \arrow[uu,"\ff^{i_1}"'] \\
 {(0,0,0)=}\;\gg \arrow[rr,"\ff^{i_0}"'] \arrow[uu,"\ff^{i_1}"]
                \arrow[ur,dashed,"\ff^{i_2}"']
 & & {(1,0,0)} \arrow[uu,"\ff^{i_1}"' near start] \arrow[ur,"\ff^{i_2}"'] &
\end{tikzcd}
\end{equation}
Every face of the cube will have a filling by a homotopy between the two paths around said face. Then the whole cube will have a filling corresponding to a higher homotopy {whose boundary is the sum of the six face homotopies, each composed with the remaining attachment before or after the face}. Note that the faces carry canonical directions induced by the edges.

\begin{rmk}
    We will give a formal proof below, but to understand why one expects this to be true, it helps to realise that there is a correspondence between the constructed functors and homotopies and the $A_\infty$-structure relations. Indeed, recall that when considering $\Phi^R$, $\Phi^L$ and $\Xi$, we calculated
    \begin{equation}
    \begin{aligned}
 &\Phi^L_k(d_{i-p+1},\ldots,d_{\widetilde{i'}+q-1})\\
 &\quad=\sum_{u=1}^{p}\sum_{v=1}^{t}
   \mu_{p+t+q-u-v+1}\Bigl(a_{i-p+1},\ldots,a_{i-u},\\
 &\qquad\qquad
   \mu_{u+v}(a_{i-u+1},\ldots,a_i,b_{j+1},\ldots,b_{j+v}),
   b_{j+v+1},\ldots,b_{j+t},c_{j'+1},\ldots,c_{j'+q}\Bigr).
\end{aligned}
\end{equation}
\begin{equation}
\begin{aligned}
 &\Phi^R_k(d_{i-p+1},\ldots,d_{\widetilde{i'}+q-1})\\
 &\quad=\sum_{v=1}^{t}\sum_{w=1}^{q}
 \mu_{p+t+q-v-w+1}\Bigl(
 a_{i-p+1},\ldots,a_i,b_{j+1},\ldots,b_{j+t-v},\\
 &\qquad\qquad
 \mu_{v+w}(b_{j+t-v+1},\ldots,b_{j+t},c_{j'+1},\ldots,c_{j'+w}),
 c_{j'+w+1},\ldots,c_{j'+q}\Bigr).
\end{aligned}
\end{equation}
The terms of $\Phi^R$ group together the right jump, and the terms of $\Phi^L$ the left jump. Their sum will remind one of \begin{equation*}
    \mu_2(\mu_2(a,b),c)+\mu_2(a,\mu_2(b,c)).
\end{equation*}
Our homotopy $\Xi$ then groups both jumps, which reminds one of $\mu_3(a,b,c)$. The functor homotopy differential can then be interpreted through the relation
\[
    \mu_1\mu_3(a,b,c)=\mu_2(\mu_2(a,b),c)+\mu_2(a,\mu_2(b,c))
\]
when $a$, $b$ and $c$ are cycles. Though below we give a more direct proof, one would expect to be able to construct a proof that interprets the functors in this manner.
\end{rmk}
\subsubsection{The homologies of parallel pushouts}
Let's set some notation to begin:
\[
 \gg:Q_m\longrightarrow\cc,\qquad
 \ff^i:Q_{n_i}\longrightarrow\cc,\qquad
 \ff^i_0(n_i)=\gg_0(i),
 \quad i\in I\subseteq\{0,\ldots,m-1\},
\]
with $m,n_i>0$. Write
\[
 d(I)=m+\sum_{i\in I}(n_i-1).
\]
We denote the iterated pushout by $Q_m\bigvee_{i\in I}Q_{n_i}$ and fix compatible object shifts. Note that the order of forming the pushout does not matter here.
For every attachment order $p\in\text{Ord}(I)$ there is a pair $(\Phi^p,E_p)$
\[
 \Phi^p:Q_{d(I)}\longrightarrow Q_m\bigvee_{i\in I}Q_{n_i},
 \qquad
 E_p:Q_m\bigvee_{i\in I}Q_{n_i}\longrightarrow\cc,
\]
where both functors use the same order and fixed choices. We then set
\[
 \#^p(\ff^I;\gg):=(\Phi^p)^*(E_p).
\]
In particular, changing the order changes both entries of the pair. Given a predetermined $(\Phi,E)$ we can attach to it a new $\ff^i$ by the operation
\[
 (\Phi,E)\longmapsto
 \bigl((\Phi\vee1)\circ\phi,\ E\vee\ff^i\bigr),
\]
when the new factor is on the right. Here $\phi$ denotes the general functor $Q_{\alpha+\beta-1}\to Q_{\alpha}\bigvee Q_{\beta}$ with the indexing understood. We will denote by $a^i$ a generic (non-identity) morphism in $Q_{n_i}$ and $b$ a generic (non-identity) morphism in $Q_m$. We will omit the indices on these morphisms to avoid an overload of notation. For example, in the already worked out case of $I=\{i,j\}$ and a gluing of $\gg:Q_m \to \cc$, $\ff^i:Q_{n_i}\to \cc$, and $\ff^j:Q_{n_j}\to \cc $ with $\ff^i_0(n_i)=\gg_0(i)$ and $\ff^j_0(n_j)=\gg_0(j)$ we could consider sequences in $Q_m\bigvee_{l\in I}Q_{n_l}$ of the form
\begin{equation*}
(a^i,\dots,a^i,b,\dots,b,a^j,\dots,a^j).
\end{equation*}
More generally, we will need to consider sequences of the form
\begin{equation*}
    (a^1,\dots,a^1,b,\dots,b,a^2,\dots,a^2,b,\dots,b,\dots,b,\dots,b,a^k,\dots,a^k).
\end{equation*}
\begin{rmk}
    For $|I|>2$ we also need to account for more complex sequences, for example
    \begin{equation*}
        (a^1,\dots,a^1,b,\dots,b,a^3,\dots,a^3,b,\dots,b,a^2,\dots,a^2),
    \end{equation*}
    as we can `use $Q_m$' to jump from any $Q_{n_i}$ to any $Q_{n_j}$ where $i$ need not be before $j$.
\end{rmk}
We therefore also need to consider sequences of the form
\begin{equation*}
(a^{i_1},\dots,a^{i_1},b,\dots,b,a^{i_2},\dots,a^{i_2},b,\dots,b,\dots,b,\dots,b,a^{i_l},\dots,a^{i_l}).
\end{equation*}
Sequences can also begin and/or terminate with $b$'s, and $l$ can be arbitrarily large, since we allow repeats, i.e., $i_s=i_{s'}$ for distinct $s$ and $s'$.

The proof will use the following corollary.
\begin{cor}\label{corhom}
 The homology groups of the morphism complexes $\bigl(Q_m\bigvee_{i\in I}Q_{n_i}\bigr)(x,y)$ have the following generators:
\begin{itemize}
    \item $[1_x]$ when $x=y$.

    \item $[a^i_y]$ when $x=y-1$ inside the copy of $Q_{n_i}$,
    with indices modulo $n_i+1$.

    \item $[b_y]$ when $x=y-1$ inside $Q_m$,
    with indices modulo $m+1$.

    \item $[\mu_2(a^i_{n_i},b_{i+1})]$ when
    $x=(n_i-1)_{Q_{n_i}},\,
        y=(i+1)_{Q_m},\,
         i\in I$.

    \item $[\mu_2(b_i,a^i_0)]$ when $x=(i-1)_{Q_m},\,
         y=0_{Q_{n_i}},\,
        i\in I$.

    \item
    $[\mu_2(\mu_2(a^i_{n_i},b_{i+1}),a^{i+1}_0)]$ when $x=(n_i-1)_{Q_{n_i}},\,y=0_{Q_{n_{i+1}}},\,i,i+1\in I$.
    \item Trivial otherwise.
\end{itemize}

\end{cor}
In other words, the homology can only occur for endomorphism complexes, generated by the unit, or for morphism complexes between consecutive elements in each copy of $Q_{n_i}$ or $Q_m$. It can also occur between the final (non-glued) element of $Q_{n_i}$ and the next (non-glued) element of $Q_m$, or between the element before a gluing in $Q_m$ and the first element in the copy of $Q_{n_i}$ which was glued. Finally, it can occur between the final (non-glued) element of $Q_{n_i}$ and the first non-glued element of $Q_{n_{i+1}}$ \textit{when} they are glued at successive elements of $Q_m$.

\begin{proof}[Proof of \ref{corhom}]
We induct on $|I|$, starting with Lemma~\ref{lemphihomology}. For $p\in I$, we can write
\[
    Q_m\bigvee_{i\in I}Q_{n_i}
    \cong
    \left(Q_m\bigvee_{i\in I\setminus\{p\}}Q_{n_i}\right)
    \bigvee_{p,n_p}Q_{n_p}.
\]
We then re-run the same argument of the lemma at this last attachment. The classes already present in either factor persist, and the only additional classes are represented by
\[
    [\mu_2(A,a^p_0)]
    \qquad\text{and}\qquad
    [\mu_2(a^p_{n_p},B)],
\]
where $[A]$ and $[B]$ are non-identity homology classes in the first factor ending and starting at the glued object $p$, respectively. Here $a^p_0$ is the first arrow after the glued object in $Q_{n_p}$. By induction, $A$ is represented by $b_p$, or by $\mu_2(a^{p-1}_{n_{p-1}},b_p)$ if $p-1$ is already attached; similarly for $B$ in the other direction. These give exactly the mixed classes listed above.

In particular, when two pushouts are taken at consecutive objects of $Q_m$, say $Q_{n_i}$ and $Q_{n_{i+1}}$ attached at $i\in Q_m$ and $i+1\in Q_m$ respectively, picture
\begin{equation}
    \begin{tikzcd}
                              & {}                                                   & 0                                          \\
                              & i+1 \arrow[u, "b_{i+2}"'] \arrow[ru, "a_{0}^{i+1}"'] &                                            \\
0                             &                                                      & n_{i+1}-1 \arrow[lu, "a^{i+1}_{n_{i+1}}"'] \\
                              & i \arrow[uu, "b_{i+1}"'] \arrow[lu, "a^i_0"]         &                                            \\
n_i-1 \arrow[ru, "a^i_{n_i}"] & {} \arrow[u, "b_i"']                                 &
\end{tikzcd}
\end{equation}
One finds that a non-trivial homology class can appear in the form
\[
    [\mu_2(\mu_2(a^i_{n_i},b_{i+1}),a^{i+1}_0)]
    =[\mu_2(a^i_{n_i},\mu_2(b_{i+1},a^{i+1}_0))].
\]
Note, though, that no non-trivial homology classes can occur as morphisms from $n_i-1$ to $0_{Q_{n_{i+2}}}$ (and in fact to $0_{Q_{n_j}}$ for $j\notin\{i,i+1\}$), since
\begin{align*}
    [a^{i+2}_0]\circ[b_{i+2}]\circ[b_{i+1}]\circ[a^i_{n_i}]
    &=[a^{i+2}_0]\circ[\mu_2(b_{i+1},b_{i+2})]\circ[a^i_{n_i}]\\
    &=0.
\end{align*}
\end{proof}

We now fix more notation. Let $I$
denote the chosen family of attachments. We will consider pairs $(\Psi^I,\Xi^I)$, generalisations of $\Psi$ and $\Xi$ produced earlier. For example, if we consider $I=\{i,j\}$ then we have homotopies
\begin{equation*}
    \Psi^{ij}:E_{ij}\implies E_{ji}\qquad
    \Xi^{ij}:\Phi^{ij}\implies \Phi^{ji}
\end{equation*}
These two homotopies allow one to fill the square with the corresponding homotopy, which we denote by $\chi^{ij}:
\#^{ij}(\ff^I;\gg)\simeq\#^{ji}(\ff^I;\gg)$:
\begin{equation}
\begin{tikzcd}[row sep=large,column sep=huge]
    {(0,1)}
        \arrow[r,"\text{attach }\ff^i"]
        \arrow[dr,phantom,"\Downarrow\,\chi^{ij}" description]
    & {(1,1)}\\
    {(0,0)}
        \arrow[u,"\text{attach }\ff^j"]
        \arrow[r,"\text{attach }\ff^i"']
    & {(1,0)}
        \arrow[u,"\text{attach }\ff^j"']
\end{tikzcd}
\end{equation}
One dimension higher, we can consider the permutations of $I=\{i,j,k\}$, to construct the following functors  \[
 \#^{ijk}(\ff^I;\gg),\quad\#^{ikj}(\ff^I;\gg),\quad\#^{jik}(\ff^I;\gg),\quad
 \#^{jki}(\ff^I;\gg),\quad\#^{kij}(\ff^I;\gg),\quad\#^{kji}(\ff^I;\gg),
\]
These correspond to the six
directed paths from $(0,0,0)$ to $(1,1,1)$ in the attachment cube:
\begin{equation}\label{eq:attachmentcube}
\begin{tikzcd}[row sep=large,column sep=large]
 & {(0,1,1)} \arrow[rr,"\text{attach }\ff^i"]
 & & {(1,1,1)}\\
 {(0,1,0)}
   \arrow[rr,"\text{attach }\ff^i" near end]
   \arrow[ur,"\text{attach }\ff^k"]
   \arrow[drrr,phantom,"\chi^{ijk}" description]
 & & {(1,1,0)}
       \arrow[ur,"\text{attach }\ff^k"']
 &\\
 & {(0,0,1)}
     \arrow[uu,dashed,"\text{attach }\ff^j" near start]
     \arrow[rr,dashed,"\text{attach }\ff^i" near start]
 & & {(1,0,1)}
       \arrow[uu,"\text{attach }\ff^j"']\\
 {(0,0,0)}
   \arrow[rr,"\text{attach }\ff^i"']
   \arrow[uu,"\text{attach }\ff^j"]
   \arrow[ur,dashed,"\text{attach }\ff^k"']
 & & {(1,0,0)}
       \arrow[uu,"\text{attach }\ff^j"' near start]
       \arrow[ur,"\text{attach }\ff^k"']
 &
\end{tikzcd}
\end{equation}
Each square face carries the corresponding two-attachment homotopy $\chi$. One should view this filling $\chi^{ijk}$ as follows: there are six paths that use a single (directed) face. For example, one can take the path $(0,0,0)\to(1,0,0)\to(1,1,1)$, where the first arrow attaches $\ff^i$, and the second represents the homotopy between the two orders of attaching $\ff^j$ and $\ff^k$ to the result. We denote this homotopy by $\chi^{i|jk}$. More generally, we use
\[\chi^{i|jk} \qquad \chi^{ij|k}\]
to denote applying the homotopy $\chi^{jk}$ \textit{after} attaching $\ff^i$, and attaching $\ff^k$ \textit{after} applying the homotopy $\chi^{ij}$, respectively.
{Explicitly, writing $\chi^i$ for the attachment of $\ff^i$, we use
\[
    \chi^{i|jk}:=\chi^{jk}\circ\chi^i,
    \qquad
    \chi^{ij|k}:=\chi^k\circ\chi^{ij}.
\]
Here $\circ$ denotes the chosen attachment operation, with each subsequent homotopy or attachment applied to the result of the preceding one. This is constructed by the obvious extension of the gluing of two functors.}
We will also use the same notation for higher analogues of $\Psi$ and $\Xi$. In general, if our indexing set is $I$, we will need to consider all such compositions
\[\chi^{A_1|A_2|\dots|A_k}:=\chi^{A_k}\circ\chi^{A_{k-1}}\circ\dots\circ\chi^{A_1},\]
where the $A_i$ are nonempty and pairwise disjoint, with $\bigcup_{i=1}^k A_i=I$.

{Returning to the $3$-cube, one might expect the equation
\[
\begin{aligned}
    \mu_1^{\mathrm{Fun}}\chi^{ijk}
    &\stackrel{?}{=}\chi^{i|jk}+\chi^{jk|i}+\chi^{j|ik}
      +\chi^{ik|j}+\chi^{k|ij}+\chi^{ij|k}.
\end{aligned}
\]
However, there is an issue with the endpoints. Each face is a $1$-homotopy between two of our output functors, but these endpoints differ for different faces. For example,
\[
    \chi^{i|jk}:\#^{ijk}\implies\#^{ikj},
    \qquad
    \chi^{ij|k}:\#^{ijk}\implies\#^{jik}.
\]
We want $\chi^{ijk}$ to compare two paths of face homotopies with initial functor $\#^{ijk}$ and terminal functor $\#^{kji}$. Thus its equation must use the differential with these fixed left and right endpoint functors, and we need
\[
    \bigl(\mu_{1;\#^{ijk},\#^{kji}}^{\mathrm{Fun}}\bigr)^2\chi^{ijk}=0.
\]
In particular, the right-hand side must be closed under this same differential. The individual face equations use their own endpoint functors, when we change to the chosen common endpoints, additional terms appear in the differential. These need to be accounted for.}

{To see how to do this, consider two successive $1$-homotopies $H_1:f_1\implies f_2$ and $H_2:f_2\implies f_3$, where all $f_i$ agree on objects and the homotopies have trivial zeroth components. Comparing the endpoint terms in the functor category differential (see \cite{Se,Or}) and using the $A_\infty$-relations gives
\[
\begin{aligned}
    \mu_{1;f_1,f_3}^{\mathrm{Fun}}(H_1+H_2)
    &=f_1+f_3
      +\mu_{1;f_1,f_3}^{\mathrm{Fun}}\mu_2^{\mathrm{Fun}}(H_1,H_2).
\end{aligned}
\]
The last term accounts for changing the right endpoint of $H_1$ from $f_2$ to $f_3$ and the left endpoint of $H_2$ from $f_2$ to $f_1$. Since we work modulo $2$, we therefore take
\[
    H_1\star H_2:=H_1+H_2+\mu_2^{\mathrm{Fun}}(H_1,H_2),
    \qquad
    \mu_{1;f_1,f_3}^{\mathrm{Fun}}(H_1\star H_2)=f_1+f_3.
\]
This composes two homotopies into one with the required endpoints.}

{Continuing with the three-cube example, where $I=\{i,j,k\}$, we now concentrate on the homotopies $\Psi^{(-|-)}$. These have different initial and final endpoints. For example,
\[
\begin{aligned}
    \Psi^{i|jk}:E_{ijk}\implies E_{ikj},\\
   \Psi^{jk|i}:E_{jki}\implies E_{kji}.
\end{aligned}
\]
More explicitly,
\[
\begin{aligned}
\Psi^{i|jk}:((\gg\vee\ff^i)\vee\ff^{j})\vee\ff^k \implies ((\gg\vee\ff^i)\vee\ff^k)\vee\ff^j,\\
\Psi^{jk|i}:((\gg\vee\ff^j)\vee\ff^k)\vee\ff^i \implies ((\gg\vee\ff^k)\vee\ff^j)\vee\ff^i.
\end{aligned}
\]

The six face homotopies give two routes. Taking left-associated $\star$ `compositions', we find
\begin{equation}\label{eq:cubicalroutes}
\begin{gathered}
(\Psi^{ij|k}\star\Psi^{j|ik})\star\Psi^{jk|i}:
 E_{ijk}\xRightarrow{\Psi^{ij|k}}E_{jik}
 \xRightarrow{\Psi^{j|ik}}E_{jki}
 \xRightarrow{\Psi^{jk|i}}E_{kji},\\
(\Psi^{i|jk}\star\Psi^{ik|j})\star\Psi^{k|ij}:
 E_{ijk}\xRightarrow{\Psi^{i|jk}}E_{ikj}
 \xRightarrow{\Psi^{ik|j}}E_{kij}
 \xRightarrow{\Psi^{k|ij}}E_{kji}.
\end{gathered}
\end{equation}
Thus the six paths through the cube using only edges become the vertices of a hexagon, and the paths using a square face become its edges.

The $2$-homotopy $\Psi^{ijk}$ fills the difference between the paths. Its equation is
\begin{equation}\label{eq:threecubeboundary}
\begin{aligned}
    \mu_1^{\mathrm{Fun}}\Psi^{ijk}
    &=(\Psi^{ij|k}\star\Psi^{j|ik})\star\Psi^{jk|i}
      +(\Psi^{i|jk}\star\Psi^{ik|j})\star\Psi^{k|ij}\\
    &=\Psi^{i|jk}+\Psi^{jk|i}+\Psi^{j|ik}
      +\Psi^{ik|j}+\Psi^{k|ij}+\Psi^{ij|k}\\
 &\quad+\mu_2^{\mathrm{Fun}}(\Psi^{ij|k},\Psi^{j|ik})\\
 &\quad+\mu_2^{\mathrm{Fun}}\Bigl(
    \bigl[\Psi^{ij|k}+\Psi^{j|ik}+\mu_2^{\mathrm{Fun}}(\Psi^{ij|k},\Psi^{j|ik})\bigr],\Psi^{jk|i}\Bigr)\\
 &\quad+\mu_2^{\mathrm{Fun}}(\Psi^{i|jk},\Psi^{ik|j})\\
 &\quad+\mu_2^{\mathrm{Fun}}\Bigl(
    \bigl[\Psi^{i|jk}+\Psi^{ik|j}+\mu_2^{\mathrm{Fun}}(\Psi^{i|jk},\Psi^{ik|j})\bigr],\Psi^{k|ij}\Bigr).
\end{aligned}
\end{equation}

More generally, for a fixed indexing set $I=\{0,1,\dots,k\}$, we require
\begin{equation}\label{eq:cubicalarity}
    \mu_1^{\mathrm{Fun}}\Psi^I
    =\sum_{\substack{A\sqcup B=I\\A,B\neq\emptyset}}\Psi^{A|B}+C^I.
\end{equation}
Here $C^I$ contains the composition terms, as in~\eqref{eq:threecubeboundary}. To specify them in higher dimensions, we must also choose how to compose the higher face homotopies. We aim to make these choices inductively to match the following constructions.

{
\begin{rmk}
If one expands the arity-$r$ equation for $\Psi^{ijk}$ in~\eqref{eq:threecubeboundary}, one finds, schematically,
\[
\begin{aligned}
    \sum_{\substack{A\sqcup B=\{i,j,k\}\\A,B\neq\emptyset}}
    \Psi_r^{A|B}(x_1,\dots,x_r)
    &=\sum\mu\bigl(\dots,\Psi^{ijk}(\dots),\dots\bigr)
      +\sum\Psi^{ijk}\bigl(\dots,\mu(\dots),\dots\bigr)\\
    &\quad+\sum\mu\bigl(\dots,\Psi^{i|jk}(\dots),\dots,
      \Psi^{ik|j}(\dots),\dots\bigr)+\cdots.
\end{aligned}
\]
i.e., the equation contains both terms involving the $2$-homotopy and terms combining the face homotopies along the chosen paths.
With this, one could include the correction terms directly in the defining equation for a $2$-homotopy between the chosen paths. We use the notation $\Psi^{A|B}$ for the individual faces and write the correction terms separately so that the combinatorics of the cube remain clearly visible.
\end{rmk}
}

   To understand these choices, we can go back to the jump from three functors, including the root $\gg$, and the $2$-cube, to four functors and the $3$-cube. Each square face of the $3$-cube carries the homotopy with the direction chosen in the two-attachment case. Once these lower homotopies are chosen, each term $\Psi^{A|B}$ is determined by a face of the $3$-cube, i.e., which pair of indices we take and whether the remaining attachment comes before or after the homotopy.

With this data and the additional choice of endpoints corresponding to $ijk$ and $kji$, the two routes around the hexagon are determined as above. Indeed, every face homotopy can be interpreted as switching two adjacent entries in an attachment order. For example, starting at $ijk$, the homotopy $\Psi^{i|jk}$ gives the switch $ijk\xrightarrow{i(jk)}ikj$. There are exactly $3!=6$ permutations. For our fixed endpoints, the minimum number of adjacent switches is $3$, and there are exactly two such minimal sequences, displayed in~\eqref{eq:cubicalroutes}. These correspond precisely to the two paths around the hexagon. Furthermore, we can interpret these two choices as two connected unions of square faces of the $3$-cube. Their intersection is precisely the two edge paths corresponding to the initial and final choices $ijk$ and $kji$, and their union recovers the entire boundary of the $3$-cube.
The three faces for the first composition are shown below, together with their flattened picture. The double arrows indicate the directions of the face homotopies.
\begin{equation*}
\begin{tikzpicture}[
  >={Stealth[length=1.8mm]},
  line cap=round,line join=round,
  edge/.style={->,draw=black!70,line width=.65pt},
  initial/.style={->,draw=blue!75!black,line width=1.45pt},
  terminal/.style={->,draw=magenta!75!black,line width=1.45pt},
  facearrow/.style={-{Implies[length=2mm,width=2.4mm]},
    double,double distance=.65pt,line width=.45pt},
  vertex/.style={circle,fill=black,inner sep=1.25pt},
  coordlabel/.style={font=\scriptsize,text=black,inner sep=1.5pt},
  edgelabel/.style={font=\small,fill=white,inner sep=1pt},
  facelabel/.style={font=\small,text=black,inner sep=1pt}
]
\color{black}

\begin{scope}[xshift=-3.2cm,yshift=-.3cm,scale=1.1]
  \coordinate (c000) at (-1.2,-1.3);
  \coordinate (c100) at (1.5,-1.3);
  \coordinate (c010) at (-1.2,1.1);
  \coordinate (c110) at (1.5,1.1);
  \coordinate (c001) at (-2.2,-.3);
  \coordinate (c101) at (.5,-.3);
  \coordinate (c011) at (-2.2,2.1);
  \coordinate (c111) at (.5,2.1);

  \fill[teal!14] (c000)--(c100)--(c110)--(c010)--cycle;
  \fill[orange!18] (c010)--(c110)--(c111)--(c011)--cycle;
  \fill[violet!13] (c000)--(c010)--(c011)--(c001)--cycle;
  \foreach \a/\b in {100/101,001/101,101/111}
    \draw[edge,dashed,draw=black!30,line width=.5pt]
      (c\a)--(c\b);
  \node[circle,fill=black!35,inner sep=1pt] at (c101) {};

  \draw[edge] (c000)--(c010);
  \draw[edge] (c010)--(c110);
  \draw[edge] (c010)--(c011);
  \draw[initial] (c000)--node[edgelabel,below,text=blue!75!black] {$i$} (c100);
  \draw[initial] (c100)--node[edgelabel,right,text=blue!75!black] {$j$} (c110);
  \draw[initial] (c110)--node[edgelabel,above right,text=blue!75!black] {$k$} (c111);
  \draw[terminal] (c000)--node[edgelabel,below left,text=magenta!75!black] {$k$} (c001);
  \draw[terminal] (c001)--node[edgelabel,left,text=magenta!75!black] {$j$} (c011);
  \draw[terminal] (c011)--node[edgelabel,above,text=magenta!75!black] {$i$} (c111);

  \draw[facearrow,draw=teal!70!black] (.45,-.85)--(-.70,.45);
  \node[facelabel] at (.80,.05) {$ij$};
  \draw[facearrow,draw=orange!80!black] (.35,1.35)--(-.95,1.85);
  \node[facelabel] at (-.10,1.93) {$ik$};
  \draw[facearrow,draw=violet!75!black] (-1.43,.70)--(-1.97,.02);
  \node[facelabel] at (-1.92,1.12) {$jk$};

  \foreach \v in {000,100,010,110,001,011,111}
    \node[vertex] at (c\v) {};
  \node[coordlabel,below left] at (c000) {$000$};
  \node[coordlabel,below right] at (c100) {$100$};
  \node[coordlabel,below left] at (c010) {$010$};
  \node[coordlabel,right] at (c110) {$110$};
  \node[coordlabel,left] at (c001) {$001$};
  \node[coordlabel,above left] at (c011) {$011$};
  \node[coordlabel,above right] at (c111) {$111$};
\end{scope}

\begin{scope}[xshift=3.5cm,yshift=.1cm,scale=1.1]
  \coordinate (h000) at (-2.05,0);
  \coordinate (h100) at (-1.025,-1.775);
  \coordinate (h110) at (1.025,-1.775);
  \coordinate (h111) at (2.05,0);
  \coordinate (h011) at (1.025,1.775);
  \coordinate (h001) at (-1.025,1.775);
  \coordinate (h010) at (0,0);

  \fill[teal!14] (h000)--(h100)--(h110)--(h010)--cycle;
  \fill[orange!18] (h010)--(h110)--(h111)--(h011)--cycle;
  \fill[violet!13] (h000)--(h010)--(h011)--(h001)--cycle;
  \draw[edge] (h000)--(h010);
  \draw[edge] (h010)--(h110);
  \draw[edge] (h010)--(h011);
  \draw[initial] (h000)--node[edgelabel,below left,text=blue!75!black] {$i$} (h100);
  \draw[initial] (h100)--node[edgelabel,below,text=blue!75!black] {$j$} (h110);
  \draw[initial] (h110)--node[edgelabel,below right,text=blue!75!black] {$k$} (h111);
  \draw[terminal] (h000)--node[edgelabel,above left,text=magenta!75!black] {$k$} (h001);
  \draw[terminal] (h001)--node[edgelabel,above,text=magenta!75!black] {$j$} (h011);
  \draw[terminal] (h011)--node[edgelabel,above right,text=magenta!75!black] {$i$} (h111);

  \draw[facearrow,draw=teal!70!black] (-.90,-1.25)--(-.30,-.40);
  \node[facelabel] at (-1.12,-.58) {$ij$};
  \draw[facearrow,draw=orange!80!black] (1.25,-.65)--(1.25,.65);
  \node[facelabel] at (1.60,0) {$ik$};
  \draw[facearrow,draw=violet!75!black] (-.30,.40)--(-.90,1.25);
  \node[facelabel] at (-1.12,.58) {$jk$};

  \foreach \v in {000,100,010,110,001,011,111}
    \node[vertex] at (h\v) {};
  \node[coordlabel,left] at (h000) {$000$};
  \node[coordlabel,below left] at (h100) {$100$};
  \node[coordlabel,below right] at (h110) {$110$};
  \node[coordlabel,right] at (h111) {$111$};
  \node[coordlabel,above right] at (h011) {$011$};
  \node[coordlabel,above left] at (h001) {$001$};
  \node[coordlabel,above left] at (h010) {$010$};
\end{scope}

\end{tikzpicture}
\end{equation*}

{
\begin{rmk}
One should view $\Psi^{ijk}$ as filling the permutahedron (at this dimension, the hexagon) of attachment orders. Its six vertices correspond to the six complete directed edge paths through the $3$-cube, and its edges correspond to the face homotopies. The filling compares the two composite paths of $1$-homotopies from $E_{ijk}$ to $E_{kji}$. Writing $H_L$ and $H_R$ for the first and second composites in~\eqref{eq:cubicalroutes}, respectively, the equation is
\[
    \mu_1^{\mathrm{Fun}}\Psi^{ijk}
    =H_L+H_R
    =\sum_{\substack{A\sqcup B=\{i,j,k\}\\A,B\neq\emptyset}}
      \Psi^{A|B}+C^{ijk}.
\]
Here $C^{ijk}$ contains the additional terms arising when we compose the individual face homotopies, accounting for their changing endpoints. The cube tells us which face homotopies occur, while their composition in the functor category determines these correction terms.
\end{rmk}
}

If we now go up to five functors, including $\gg$, and thus the $4$-cube, we retain the choices on each $3$-cube face and choose $ijkl$ and $lkji$ as initial and terminal attachment orders. There are $4!$ edge paths corresponding to orders of gluings. We now need to compare paths of $1$-homotopies from $ijkl$ to $lkji$ using the $2$-homotopies, for example $\Psi^{ijk|l}$, but now also $\Psi^{ij|kl}$. These have endpoint orders $ijkl\xrightarrow{(ijk)l}kjil$ and $ijkl\xrightarrow{(ij)(kl)}jilk$, respectively; here the arrows record only the endpoint orders, while each $2$-homotopy compares two paths between them. We can compose a pair of $2$-homotopies along a shared $1$-homotopy, represented by a sub-$2$-cube or a sequence of such squares. Each of these basic $2$-homotopies is represented by either a sub-$3$-cube or a pair of sub-$2$-cubes meeting at a vertex.
\begin{equation*}
\begin{tikzpicture}[
  >={Stealth[length=1.8mm]},
  line cap=round,line join=round,
  edge/.style={->,draw=black!65,line width=.65pt},
  hidden/.style={edge,dashed,draw=black!30,line width=.5pt},
  initial/.style={->,draw=blue!75!black,line width=1.35pt},
  terminal/.style={->,draw=magenta!75!black,line width=1.35pt},
  shared/.style={->,draw=orange!80!black,line width=1.15pt},
  facearrow/.style={-{Implies[length=2mm,width=2.4mm]},
    double,double distance=.65pt,line width=.45pt},
  vertex/.style={circle,fill=black,inner sep=1.2pt},
  coordlabel/.style={font=\scriptsize,text=black,inner sep=1.5pt},
  edgelabel/.style={font=\small,fill=white,inner sep=1pt},
  homlabel/.style={font=\small,text=black,inner sep=2pt},
  homdirection/.style={font=\Large,inner sep=0pt}
]
\color{black}

\begin{scope}[xshift=-2.4cm]
  \coordinate (u0000) at (-2.8,-1.25);
  \coordinate (u1000) at (-1.2,-1.25);
  \coordinate (u1100) at (1.25,-1.25);
  \coordinate (u1010) at (-1.2,1.05);
  \coordinate (u1110) at (1.25,1.05);
  \coordinate (u1001) at (-2.15,-.30);
  \coordinate (u1101) at (.30,-.30);
  \coordinate (u1011) at (-2.15,2);
  \coordinate (u1111) at (.30,2);

  \fill[teal!14] (u1000)--(u1100)--(u1110)--(u1010)--cycle;
  \fill[teal!8] (u1010)--(u1110)--(u1111)--(u1011)--cycle;
  \fill[teal!21] (u1000)--(u1010)--(u1011)--(u1001)--cycle;
  \foreach \a/\b in {1100/1101,1001/1101,1101/1111}
    \draw[hidden] (u\a)--(u\b);
  \draw[edge,line width=1pt] (u0000)--node[edgelabel,above] {$i$} (u1000);
  \foreach \a/\b in {1000/1010,1010/1110,1010/1011}
    \draw[edge] (u\a)--(u\b);
  \draw[initial] (u1000)--node[edgelabel,below,text=blue!75!black] {$j$} (u1100);
  \draw[initial] (u1100)--node[edgelabel,right,text=blue!75!black] {$k$} (u1110);
  \draw[initial] (u1110)--node[edgelabel,above right,text=blue!75!black] {$l$} (u1111);
  \draw[terminal] (u1000)--node[edgelabel,below left,text=magenta!75!black] {$l$} (u1001);
  \draw[terminal] (u1001)--node[edgelabel,left,text=magenta!75!black] {$k$} (u1011);
  \draw[terminal] (u1011)--node[edgelabel,above,text=magenta!75!black] {$j$} (u1111);
  \node[homlabel,fill=teal!14] at (-.10,.62) {$\Psi^{i|jkl}$};
  \node[homdirection,rotate=-42,text=teal!75!black] at (-.25,-.05) {$\Rrightarrow$};
  \foreach \v in {0000,1000,1100,1010,1110,1001,1011,1111}
    \node[vertex] at (u\v) {};
  \node[circle,fill=black!35,inner sep=1pt] at (u1101) {};
  \node[coordlabel,below left] at (u0000) {$0000$};
  \node[coordlabel,below] at (u1000) {$1000$};
  \node[coordlabel,above right] at (u1111) {$1111$};
\end{scope}

\begin{scope}[xshift=3.1cm]
  \coordinate (s0000) at (-1.8,-1.8);
  \coordinate (s1000) at (0,-1.8);
  \coordinate (s0100) at (-1.8,0);
  \coordinate (s1100) at (0,0);
  \coordinate (s1110) at (1.8,0);
  \coordinate (s1101) at (0,1.8);
  \coordinate (s1111) at (1.8,1.8);
  \fill[teal!14] (s0000)--(s1000)--(s1100)--(s0100)--cycle;
  \fill[violet!13] (s1100)--(s1110)--(s1111)--(s1101)--cycle;
  \draw[initial] (s0000)--node[edgelabel,below,text=blue!75!black] {$i$} (s1000);
  \draw[initial] (s1000)--node[edgelabel,right,text=blue!75!black] {$j$} (s1100);
  \draw[initial] (s1100)--node[edgelabel,below,text=blue!75!black] {$k$} (s1110);
  \draw[initial] (s1110)--node[edgelabel,right,text=blue!75!black] {$l$} (s1111);
  \draw[terminal] (s0000)--node[edgelabel,left,text=magenta!75!black] {$j$} (s0100);
  \draw[terminal] (s0100)--node[edgelabel,above,text=magenta!75!black] {$i$} (s1100);
  \draw[terminal] (s1100)--node[edgelabel,left,text=magenta!75!black] {$l$} (s1101);
  \draw[terminal] (s1101)--node[edgelabel,above,text=magenta!75!black] {$k$} (s1111);
  \draw[facearrow,draw=teal!70!black] (-.35,-1.45)--(-1.45,-.35);
  \node[homlabel] at (-1.35,-1.25) {$ij$};
  \draw[facearrow,draw=violet!75!black] (1.45,.35)--(.35,1.45);
  \node[homlabel] at (1.25,1.25) {$kl$};
  \node[homlabel] at (1,-1) {$\Psi^{ij|kl}$};
  \foreach \v in {0000,1000,0100,1100,1110,1101,1111}
    \node[vertex] at (s\v) {};
  \node[coordlabel,below left] at (s0000) {$0000$};
  \node[coordlabel,above left,fill=white] at (s1100) {$1100$};
  \node[coordlabel,above right] at (s1111) {$1111$};
\end{scope}

\begin{scope}[xshift=-2.8cm,yshift=-5cm,scale=.78]
  \coordinate (g0000) at (-3.5,-1.3);
  \coordinate (g0100) at (-3.5,1.1);
  \coordinate (g0010) at (-4.45,-.35);
  \coordinate (g0110) at (-4.45,2.05);
  \coordinate (g1000) at (-.5,-1.3);
  \coordinate (g1100) at (-.5,1.1);
  \coordinate (g1010) at (-1.45,-.35);
  \coordinate (g1110) at (-1.45,2.05);
  \coordinate (g1001) at (2.5,-1.3);
  \coordinate (g1101) at (2.5,1.1);
  \coordinate (g1011) at (1.55,-.35);
  \coordinate (g1111) at (1.55,2.05);

  \fill[teal!12] (g0000)--(g1000)--(g1100)--(g0100)--cycle;
  \fill[teal!7] (g0100)--(g1100)--(g1110)--(g0110)--cycle;
  \fill[teal!18] (g0000)--(g0100)--(g0110)--(g0010)--cycle;
  \fill[violet!12] (g1000)--(g1001)--(g1101)--(g1100)--cycle;
  \fill[violet!7] (g1100)--(g1101)--(g1111)--(g1110)--cycle;
  \fill[orange!25] (g1000)--(g1100)--(g1110)--(g1010)--cycle;

  \foreach \a/\b in {0010/1010,1010/1011,1001/1011,1011/1111}
    \draw[hidden] (g\a)--(g\b);
  \foreach \a/\b in {0000/0100,0000/0010,0010/0110,0100/0110,
    0000/1000,0100/1100,0110/1110,1000/1001,1100/1101,
    1110/1111,1001/1101,1101/1111}
    \draw[edge] (g\a)--(g\b);
  \foreach \a/\b in {1000/1100,1000/1010,1100/1110,1010/1110}
    \draw[shared] (g\a)--(g\b);
  \draw[facearrow,draw=orange!80!black] (-.74,.74)--(-1.21,.01);
  \node[homlabel,text=orange!80!black] at (-1.17,1.05) {$jk$};
  \node[homlabel,fill=teal!12] at (-2.65,.60) {$\Psi^{ijk|l}$};
  \node[homdirection,rotate=-35,text=teal!75!black] at (-2.65,-.03) {$\Rrightarrow$};
  \node[homlabel,fill=violet!12] at (.85,.60) {$\Psi^{i|jkl}$};
  \node[homdirection,rotate=20,text=violet!75!black] at (.85,-.03) {$\Rrightarrow$};

  \foreach \v in {0000,0100,0010,0110,1000,1100,1010,1110,1001,1101,1111}
    \node[vertex] at (g\v) {};
  \node[circle,fill=black!35,inner sep=1pt] at (g1011) {};
  \node[coordlabel,below left] at (g0000) {$0000$};
  \node[coordlabel,below right] at (g1000) {$1000$};
  \node[coordlabel,above right] at (g1100) {$1100$};
  \node[coordlabel,left] at (g1010) {$1010$};
  \node[coordlabel,above left] at (g1110) {$1110$};
  \node[coordlabel,above right] at (g1111) {$1111$};
\end{scope}

\begin{scope}[xshift=1.6cm,yshift=-5cm]
  \coordinate (v0000) at (0,-1.2);
  \coordinate (v1000) at (2,-1.2);
  \coordinate (v0100) at (0,.8);
  \coordinate (v1100) at (2,.8);
  \coordinate (v0010) at (-.75,-.45);
  \coordinate (v1010) at (1.25,-.45);
  \coordinate (v0110) at (-.75,1.55);
  \coordinate (v1110) at (1.25,1.55);
  \coordinate (v1101) at (3.6,.8);
  \coordinate (v1111) at (2.85,1.55);

  \fill[teal!7] (v0100)--(v1100)--(v1110)--(v0110)--cycle;
  \fill[teal!18] (v0000)--(v0100)--(v0110)--(v0010)--cycle;
  \fill[orange!25] (v0000)--(v1000)--(v1100)--(v0100)--cycle;
  \fill[violet!13] (v1100)--(v1101)--(v1111)--(v1110)--cycle;
  \foreach \a/\b in {1000/1010,0010/1010,1010/1110}
    \draw[hidden] (v\a)--(v\b);
  \foreach \a/\b in {0000/0010,0100/0110,0010/0110,0110/1110,
    1100/1101,1101/1111,1110/1111}
    \draw[edge] (v\a)--(v\b);
  \foreach \a/\b in {0000/1000,0000/0100,1000/1100,0100/1100}
    \draw[shared] (v\a)--(v\b);
  \draw[edge,draw=violet!75!black,line width=1.15pt]
    (v1100)--(v1110);

  \node[homlabel,text=orange!80!black] at (.90,.51) {$ij$};
  \node[homlabel,fill=orange!25] at (.85,-.08) {$\Psi^{ijk|l}$};
  \node[homdirection,rotate=-45,text=teal!75!black] at (.90,-.65) {$\Rrightarrow$};
  \draw[facearrow,draw=violet!75!black] (1.85,1.36)--(2.98,1);
  \node[homlabel] at (2.80,1.40) {$kl$};

  \foreach \v in {0000,1000,0100,1100,0010,0110,1110,1101,1111}
    \node[vertex] at (v\v) {};
  \node[circle,fill=black!35,inner sep=1pt] at (v1010) {};
  \node[coordlabel,below left] at (v0000) {$0000$};
  \node[coordlabel,below right] at (v1100) {$1100$};
  \node[coordlabel,above left] at (v1110) {$1110$};
  \node[coordlabel,above right] at (v1111) {$1111$};
\end{scope}
\end{tikzpicture}
\end{equation*}
We have fixed the two attachment orders $ijkl$ and $lkji$, each represented by an edge path through the $4$-cube. We now choose an initial and a terminal $1$-homotopy between them, given by sequences of directed square faces. For example, we can take
\[
    (ij)(ik)(il)(jk)(jl)(kl)
    \qquad\text{and}\qquad
    (kl)(jl)(jk)(il)(ik)(ij),
\]
where we read from left to right, with each pair indicating the two adjacent entries being switched in the current attachment order. These give two edge paths in the three-dimensional permutohedron (see {\cite{Po}} for background on permutohedra). We can pass from one path to the other using the $2$-homotopies represented by its hexagonal and square faces. Here, hexagonal faces correspond to sub-$3$-cubes and square faces to pairs of sub-$2$-cubes joined at a vertex. We can then choose two compatible composites of these $2$-homotopies, both running between the chosen $1$-homotopies. The filling between these composites is the $3$-homotopy $\Psi^{ijkl}$, which is what we mean by filling the $4$-cube. This explains the face terms and the additional composition terms in~\eqref{eq:cubicalarity}.
}

\subsubsection{The higher $\Psi$s in parallel}
We first show that we can construct our higher $\Psi^I$'s. Note that we only seek the existence of a compatible family. We write $L^I,R^I$ for the endpoint functors chosen for the filling $\Psi^I$ (see below). These may vary between fillings, but remain fixed throughout the recursion for a given $\Psi^I$. To begin, assume we have already constructed all lower-order $\Psi^A$'s for nonempty $A\subsetneq I$, together with the face families $\Psi^{A|B}$ for $A,B\neq\emptyset$, $A\cap B=\emptyset$ and $A\cup B=I$.
{We use the coherent boundary just described and write
\[
    \widehat D^I:=D^I+C^I,
\]
where $D^I$ is the face sum below and $C^I$ contains the composition terms, as in~\eqref{eq:threecubeboundary}. {These use only previously constructed face homotopies and the operations $\mu_r^{\mathrm{Fun}}$.}}
Similarly to the construction of $\Psi^{ij}$, we take $\Psi^I_k=0$ for all $k\neq1$ and recursively define $\Psi^I_1$ on expressions $\mu_k(x_1,\dots,x_k)$ for composable morphisms in $Q_m\bigvee_{i\in I}Q_{n_i}$. To begin, we set $\Psi^I_1(\mu_k(x_1,\dots,x_k))=0$ whenever every $x_j$ belongs to the same $Q_{n_i}$ or to $Q_m$. We also set $\Psi^I_1=0$ on individual arrows and units. For $k\geq2$, the equation we need to solve is
\begin{equation}\label{higherPsieq}
\begin{aligned}
 &\Psi^I_1\mu_k(x_1,\dots,x_k)={\widehat D^I_k(x_1,\dots,x_k)}\\
 &\quad+\sum_{v=1}^{k}\sum\mu_{\alpha+1+\beta}\Bigl(
 L^I_{p_1}(x_1,\ldots,x_{p_1}),\ldots,
 L^I_{p_\alpha}(x_{v-p_\alpha},\ldots,x_{v-1}),\\
 &\qquad\Psi^I_1(x_v),
 R^I_{t_1}(x_{v+1},\ldots,x_{v+t_1}),\ldots,
 R^I_{t_\beta}(x_{k-t_\beta+1},\ldots,x_k)\Bigr),
\end{aligned}
\end{equation}
where we retain the notation
\begin{equation}
    D^I_k(x_1,\dots,x_k):=
    \sum_{\substack{A\cup B=I,\;A\cap B=\emptyset\\A,B\neq\emptyset}}
    \Psi^{A|B}_k(x_1,\dots,x_k).
\end{equation}
All block lengths are positive; either list of blocks may be empty. Note that equation~\eqref{higherPsieq} is the $k$-th component of the homotopy equation, assuming $\Psi^I_l=0$ for $l\neq1$. With the chosen coherent face compositions, we find that
\[
    \mu_1^{\mathrm{Fun}}\widehat D^I=0.
\]
Geometrically, each codimension-two face occurs in the boundary of exactly two codimension-one faces. Thus in taking the boundary they will pairwise cancel.

The key point is that in equation~\eqref{higherPsieq}, every term on the right-hand side is already known. For example, if every $x_v$ is an individual $a^{i_j}$ or $b$, then we obtain
\begin{equation*}
    \Psi^I_1\mu_k(x_1,\dots,x_k)={\widehat D^{I}_k(x_1,\dots,x_k)}.
\end{equation*}
Furthermore, the unary equation follows by the same induction as in the $\Psi^{ij}$ case, now using {$\mu_1^{\mathrm{Fun}}\widehat D^I=0$} for the chosen endpoints $L^I,R^I$.
{Indeed, set $Z=\mu_1^{\mathrm{Fun}}\Psi^I+\widehat D^I$. By construction, $Z_0=0$ and $Z_k=0$ for $k\geq2$, while $\mu_1^{\mathrm{Fun}}Z=0$. Also, $Z_1$ vanishes on the original morphisms. The $k$-th component of $\mu_1^{\mathrm{Fun}}Z=0$ (one can check) will result in
\[
    Z_1\mu_k(x_1,\dots,x_k)
    =\sum_{v=1}^k\sum\mu_{\alpha+1+\beta}
    \bigl(L^I_{p_1},\dots,L^I_{p_\alpha},Z_1(x_v),
    R^I_{t_1},\dots,R^I_{t_\beta}\bigr).
\]
Just as in the $\Psi^{ij}$ case, with the grouping of terms and cancellations now given by the combinatorial structure of the permutahedra. 
Each $x_v$ has fewer vertices, so induction gives $Z_1=0$ on every tree.} Thus, We obtain
\[
    \mu_1\Psi^I_1(x)+\Psi^I_1(\mu_1x)={\widehat D^I_1(x)}.
\]

\subsubsection{The higher $\Xi$s in parallel}

Now we construct the homotopies corresponding to the various functors of the form $\Phi^p:Q_{d(I)}\to  Q_m \bigvee_{i\in I}Q_{n_i}$. For this, assume we have built all such $\Xi^A$'s for $|A|\leq h$. For $|I|=h+1$, we construct the components in increasing input length $k$. We set $\Xi^I_0=0$, and $\Xi^I_k(\dots,1,\dots)=0$ for all $k\geq 1$.

Here $L$ and $R$ denote the chosen $\Phi$-functors at the two ends. {Using the same choices of face compositions as above, now for the lower $\Xi$'s, we write
\[
    \widehat B^I:=B^I+C_\Xi^I,
    \qquad \mu_1^{\mathrm{Fun}}\Xi^I=\widehat B^I.
\]
Here $C_\Xi^I$ contains the composition terms for the $\Xi$-family; all its terms are already known from the lower homotopies. We use normalized face compositions, so $\widehat B^I_0=0$ and $\widehat B^I_k(\dots,1,\dots)=0$.} For composable non-identity source arrows $x_1,\dots,x_k$, the equation we need to satisfy is
\begin{equation}\label{eq:xirecursion}
\begin{aligned}
 &\mu_1\Xi^I_k(x_1,\ldots,x_k)=\widehat B^I_k(x_1,\ldots,x_k)\\
 &\quad+\sum_{e=1}^{k-1}\sum_{v=1}^{k-e+1}\sum\mu_{\alpha+1+\beta}\Bigl(
 L_{p_1}(x_1,\ldots,x_{p_1}),\ldots,
 L_{p_\alpha}(x_{v-p_\alpha},\ldots,x_{v-1}),\\
 &\qquad\Xi^I_e(x_v,\ldots,x_{v+e-1}),
 R_{t_1}(x_{v+e},\ldots,x_{v+e+t_1-1}),\ldots,
 R_{t_\beta}(x_{k-t_\beta+1},\ldots,x_k)\Bigr),
\end{aligned}
\end{equation}
where $B^I_k$ is given by
\begin{equation}
  B^I_k(x_1,\dots,x_k):=  \sum_{\substack{A\cup B=I,\;A\cap B=\emptyset\\A,B\neq\emptyset}} \Xi^{A|B}_k(x_1,\dots,x_k).
\end{equation}
{For example, with the same compositions as in~\eqref{eq:threecubeboundary},
\[
\begin{aligned}
    \widehat B^{ijk}
    &=(\Xi^{ij|k}\star\Xi^{j|ik})\star\Xi^{jk|i}\\
    &\quad+(\Xi^{i|jk}\star\Xi^{ik|j})\star\Xi^{k|ij}.
\end{aligned}
\]
Thus $C_\Xi^{ijk}$ is given by the composition terms in~\eqref{eq:threecubeboundary}, with each $\Psi$ replaced by $\Xi$.}
The argument then runs in the same manner as that for $\Xi^{ij}$: we show that the right-hand side is a cycle, then argue that it is a boundary, now using Corollary~\ref{corhom}.

To see that the RHS is a cycle, we recall that the equation is secretly the rearrangement of (as in the $\Xi^{ij}$ case)
\begin{equation*}
    (\mu_1^{\mathrm{Fun}}\Xi^I)_k=\widehat B^I_k.
\end{equation*}
Indeed, we have
\[
\begin{aligned}
 &(\mu_1^{\mathrm{Fun}}\Xi^I)_k(x_1,\ldots,x_k)
 =\mu_1\Xi^I_k(x_1,\ldots,x_k)\\
 &\quad+\sum_{e=1}^{k}\sum_{j=0}^{k-e}
 \Xi^I_{k-e+1}\bigl(x_1,\ldots,x_j,
 \mu_e(x_{j+1},\ldots,x_{j+e}),x_{j+e+1},\ldots,x_k\bigr)\\
 &\quad+\sum_{e=1}^{k-1}\sum_{v=1}^{k-e+1}\sum\mu_{\alpha+1+\beta}\Bigl(
 L_{p_1}(x_1,\ldots,x_{p_1}),\ldots,
 L_{p_\alpha}(x_{v-p_\alpha},\ldots,x_{v-1}),\\
 &\qquad\Xi^I_e(x_v,\ldots,x_{v+e-1}),
 R_{t_1}(x_{v+e},\ldots,x_{v+e+t_1-1}),\ldots,
 R_{t_\beta}(x_{k-t_\beta+1},\ldots,x_k)\Bigr).
\end{aligned}
\]
Thus we can write the RHS of equation~\eqref{eq:xirecursion} as
\begin{equation}
    \mathrm{RHS}=\bigl(\widehat B^I+\mu_1^{\mathrm{Fun}}\Xi^{I,<k}\bigr)_k,
\end{equation}
where $\Xi^{I,<k}_l=\Xi^I_l$ for $l<k$ and is zero for $l\geq k$. The point is that $(\mu_1^{\text{Fun}}\Xi^I)_k=(\mu_1^\text{Fun}\Xi^{I,<k})_k+\mu_1\Xi^I_k$. One then finds
\begin{equation}
    \mu_1^{\text{Fun}}\left[\widehat B^I+\mu_1^{\text{Fun}}\Xi^{I,<k}\right]=\mu_1^{\text{Fun}}\widehat B^I,
\end{equation}
i.e., it remains to show $\mu_1^{\text{Fun}}\widehat B^I=0$. But this again is just the combinatorial structure of permutahedra, indeed, for example, if $I-=\{i,j,k\}$ then 
\begin{equation*}
\begin{aligned}
    (\mu_1^\text{Fun}\widehat B^{ijk})_1
    &=\Bigl(\mu_1^{\text{Fun}}\bigl[
    \Xi^{i|jk}+\Xi^{jk|i}+\Xi^{j|ik}\\
    &\qquad+\Xi^{ik|j}+\Xi^{k|ij}+\Xi^{ij|k}\bigr]\Bigr)_1.
\end{aligned}
\end{equation*}
One finds that
 \[
\begin{aligned}
(\mu_1^{\text{Fun}}\Xi^{i|jk})_1
&=(\Phi^{ijk}_1+\Phi^{ikj}_1),\\
(\mu_1^{\text{Fun}}\Xi^{jk|i})_1
&=(\Phi^{jki}_1+\Phi^{kji}_1)
\end{aligned}
\]
and similarly for the other permutations, giving
\[
\begin{aligned}
&\Bigl(\mu_1^{\text{Fun}}\bigl[
\Xi^{i|jk}+\Xi^{j|ik}+\Xi^{k|ij}
+\Xi^{ij|k}+\Xi^{ik|j}+\Xi^{jk|i}
\bigr]\Bigr)_1\\
&=
(\Phi^{ijk}_1+\Phi^{ikj}_1)
+(\Phi^{jik}_1+\Phi^{jki}_1)
+(\Phi^{kij}_1+\Phi^{kji}_1)\\
&\qquad
+(\Phi^{ijk}_1+\Phi^{jik}_1)
+(\Phi^{ikj}_1+\Phi^{kij}_1)
+(\Phi^{jki}_1+\Phi^{kji}_1)
\\
&=0,
\end{aligned}
\]
as everything cancels in pairs. It then follows that
\begin{equation}
    \mu_1^{\text{Fun}}\left[\widehat B^I+\mu_1^{\text{Fun}}\Xi^{I,<k}\right]=\mu_1^{\text{Fun}}\widehat B^I=0.
\end{equation}

It follows then that the RHS of equation~\eqref{eq:xirecursion} is a cycle. Now we use Corollary~\ref{corhom} to see that, for any sequence $(x_1,\dots,x_k)$, the RHS applied to the sequence must be a boundary for degree reasons, as in the construction of $\Xi^{ij}$. Thus we can recursively choose every $\Xi^I_k$ satisfying equation~\eqref{eq:xirecursion}. One should note that, since $|I|>2$, the degree of the right-hand side satisfies
\begin{equation}
    \sum_{i=1}^k|x_i|+k+|I|-3\geq\sum_{i=1}^k|x_i|+k-1.
\end{equation}

\begin{rmk}
For each order $p$, we take $\#^p(\ff^I;\gg):=(\Phi^p)^*E_p$, using the choices above. Applying $A_\infty$-composition to the families $\Xi$ and $\Psi$, including its higher composition terms, gives a family $\Gamma^I$ of homotopies between these pullbacks. This proves Theorem~\ref{thmparallelgluings}.
\end{rmk}

\subsection{Ordered Trees}\label{subsec:ordered-trees}
We can reinterpret the parallel gluing combinatorics in terms of discretely labelled trees. As mentioned in the introduction, this idea is similar to the Boardman--Vogt construction \cite{BV} and its extension to coloured operads \cite{BM}, though we work in a discrete setting, with labels representing attachment orders. Indeed, we can picture the fixed functor $\gg:Q_m\to\cc$ as a vertex, the root of our tree. Every attachment of an $\ff^i:Q_{n_i}\to\cc$ will then be viewed as attaching a branch with a single vertex corresponding to $\ff^i$ to the root vertex. For example, the following diagram
\begin{equation}
\begin{tikzcd}[row sep=large,column sep=large]
\ff^i & \ff^j & \ff^k\\
& \gg \arrow[ul,no head] \arrow[u,no head] \arrow[ur,no head] &
\end{tikzcd}
\end{equation}
is a tree representing the attachment of $\ff^i,\ff^j,\ff^k$ to some $\gg$. The diagram, however, does not tell us in what order the branches were attached. One can keep track of this by associating to every branch $b$ a number $\text{ord}(b)\in\N$. For example,
\begin{equation}
\begin{tikzcd}[row sep=large,column sep=large]
\ff^i & \ff^j & \ff^k\\
& \gg \arrow[ul,no head,"2"] \arrow[u,no head,"0"] \arrow[ur,no head,"1"'] &
\end{tikzcd}
\end{equation}
now represents the functor $\#^{jki}(\ff^I;\gg)$, where we attach in the order $j,k,i$. We can also keep track of homotopies. For example,
\begin{equation}\label{homtreeeee}
\begin{tikzcd}[row sep=large,column sep=large]
\ff^i & \ff^j & \ff^k\\
& \gg \arrow[ul,no head,"1"] \arrow[u,no head,"0"] \arrow[ur,no head,"1"'] &
\end{tikzcd}
\end{equation}
represents the homotopy $\Gamma^{j|ik}$, i.e., we attach $\ff^j$, then apply the homotopy $\Gamma^{ik}$. We view this as simultaneously attaching $\ff^i$ and $\ff^k$.

Consider $\text{Tree}^{\gg}(I)$, the collection of formal sums of such trees, with branches corresponding to $i\in I\subseteq\{0,\dots,m-1\}$, where $\gg:Q_m\to\cc$, and with numbering given by $\text{ord}(i)\in\{0,\dots,|I|-1\}$. Furthermore, we require the labels that occur to be consecutive, starting at $0$, i.e., if an edge has label $e$, then every label $0,\dots,e$ must occur. There is a map
\[\mu:\text{Tree}^{\gg}(I)\to\text{Tree}^{\gg}(I)\]
given as follows. Write $k=|I|$ and $I=\{i_0<i_1<\dots<i_{k-1}\}$. We denote a single tree $T$ by its sequence of order labels
\[
    T=(\text{ord}(i_0),\dots,\text{ord}(i_{k-1})).
\]
If label $q$ occurs $m_q$ times, it contributes $m_q-1$ to the degree. Thus, writing $0,\dots,s-1$ for the distinct edge labels, we have
\[
    \deg(T)=\sum_{q=0}^{s-1}(m_q-1)=k-s.
\]
We define $\mu$ by
\begin{equation}\label{treemu}
\mu(a_0,\ldots,a_{k-1})
=
\sum_{q=0}^{\max\{a_0,\dots,a_{k-1}\}}
\ \sum_{\varnothing\neq S\subsetneq\{j:a_j=q\}}
(a^{q,S}_0,\ldots,a^{q,S}_{k-1}),
\end{equation}
where we set
\begin{equation}
a^{q,S}_j=
\begin{cases}
a_j+1,&a_j>q\ \text{or}\ j\in S,\\
a_j,&\text{otherwise}.
\end{cases}
\end{equation}
The map $\mu$ lowers degree by one.
\begin{rmk}
We will see in Proposition~\ref{tree-differential} that this $\mu$ is a differential, i.e., $\mu^2=0$.
\end{rmk}

Unpacking this, the outer sum runs over the order labels $q$. For each $q$, we choose a nonempty proper subset $S$ of the branches labelled $q$, increase their labels by $1$, and also increase every label greater than $q$ by $1$. An example will make it clear. Consider the tree in \eqref{homtreeeee}, which we view as $(a_0,a_1,a_2)=(1,0,1)$. Here, $\{l:a_l=0\}=\{1\}$ and $\{l:a_l=1\}=\{0,2\}$, giving
\[
\begin{aligned}
\mu(1,0,1)
={}&
\sum_{\varnothing\neq S\subsetneq\{1\}}
(a^{0,S}_0,a^{0,S}_1,a^{0,S}_2)
+
\sum_{\varnothing\neq S\subsetneq\{0,2\}}
(a^{1,S}_0,a^{1,S}_1,a^{1,S}_2)\\
={}&0+
(a^{1,\{0\}}_0,a^{1,\{0\}}_1,a^{1,\{0\}}_2)
+
(a^{1,\{2\}}_0,a^{1,\{2\}}_1,a^{1,\{2\}}_2)\\
={}&(2,0,1)+(1,0,2).
\end{aligned}
\]
We can picture this as
\begin{equation}
\begin{tikzpicture}[baseline=(current bounding box.center)]
  \foreach \xx/\aa/\bb/\cc in {0/1/0/1,4/2/0/1,8/1/0/2}{
    \begin{scope}[xshift=\xx cm]
      \node (G) at (0,0) {$\gg$};
      \node (Fi) at (-1.1,1.5) {$\ff^i$};
      \node (Fj) at (0,1.5) {$\ff^j$};
      \node (Fk) at (1.1,1.5) {$\ff^k$};
      \draw (G)--node[pos=.55,left] {$\aa$} (Fi);
      \draw (G)--node[pos=.55,left] {$\bb$} (Fj);
      \draw (G)--node[pos=.55,right] {$\cc$} (Fk);
    \end{scope}
  }
  \node at (2,.75) {$\overset{\mu}{\longmapsto}$};
  \node at (6,.75) {$+$};
\end{tikzpicture}
\end{equation}
A more complicated example would be $I=\{i,j,k,l\}$, with $\ff^i$ and $\ff^k$ attached simultaneously first, then $\ff^j$ and $\ff^l$ attached simultaneously second, i.e., a tree of the form $(0,1,0,1)$. We find that
\[
    \mu(0,1,0,1)=(1,2,0,2)+(0,2,1,2)+(0,2,0,1)+(0,1,0,2),
\]
i.e.,
\begin{equation}
\begin{tikzpicture}[baseline=(current bounding box.center)]
  \foreach \xx/\yy/\aa/\bb/\cc/\dd in {
    0/0/0/1/0/1,
    4.5/1.3/0/2/1/2,
    8.5/1.3/1/2/0/2,
    4.5/-1.3/0/1/0/2,
    8.5/-1.3/0/2/0/1}{
    \begin{scope}[shift={(\xx,\yy)}]
      \node (G) at (0,0) {$\gg$};
      \node (Fi) at (-1.35,1.5) {$\ff^i$};
      \node (Fj) at (-0.45,1.5) {$\ff^j$};
      \node (Fk) at (0.45,1.5) {$\ff^k$};
      \node (Fl) at (1.35,1.5) {$\ff^l$};
      \draw (G)--node[pos=.6,left] {$\aa$} (Fi);
      \draw (G)--node[pos=.6,left] {$\bb$} (Fj);
      \draw (G)--node[pos=.6,right] {$\cc$} (Fk);
      \draw (G)--node[pos=.6,right] {$\dd$} (Fl);
    \end{scope}
  }
  \node at (2.25,0.75) {$\overset{\mu}{\longmapsto}$};
  \node at (6.5,2.05) {$+$};
  \node at (2.4,-0.55) {$+$};
  \node at (6.5,-0.55) {$+$};
\end{tikzpicture}
\end{equation}
The reason for viewing our cubical homotopies in these terms will become apparent shortly.

\subsection{Series Gluings}\label{subsec:series-gluings}
The previous section described the combinatorics of parallel gluings, i.e., several functors attached to a single fixed functor. We realised this can be viewed in terms of attaching several vertices to a single common root. We can also glue functors sequentially, i.e., in series. That is, we consider gluing a family $\ff^0,\dots,\ff^k$ together by attaching each $\ff^i$ to $\ff^{i+1}$. Fortunately, the combinatorics here are simpler and resemble those of the Stasheff associahedra {\cite{St1,St2}}. We now make this formal. Let's fix some more notation. The idea is to break up the orders of gluing via parentheses. For example, \[
  (\ff^0\#\ff^1)\#\ff^2,
  \qquad
  \ff^0\#(\ff^1\#\ff^2)
\]
represent the two orders of attaching three functors in series. The structure is that of the associahedron $K(3)$, i.e., a line segment: there are two vertices corresponding to the functors and a single edge corresponding to the homotopy. So assume that we have a family of functors \[
  \ff^i:Q_{n_i}\longrightarrow\cc,
  \qquad 0\leq i\leq k,
\]
with $n_i>0$, such that 
\[
  \ff^i_0(n_i)=\ff^{i+1}_0(j_{i+1}),
  \qquad 0\leq i<k, \quad 0\leq j_{i+1}<n_{i+1}.
\]
We fix compatible object shifts throughout and suppress the
attachment indices in the notation $\#$. Now consider the Stasheff polyhedron $K(k+1)$. It has 
\begin{equation*}
    C_k=\frac{1}{k+1}{2k \choose k}
\end{equation*}
vertices, i.e., the $k$-th Catalan number. These vertices correspond to full binary bracketings of the sequence $0,1,\dots,k$ (see {\cite{MTTV}}). For example, there is a vertex in $K(4)$ corresponding to $((0,1),(2,3))$. In terms of functors $\ff^0,\dots ,\ff^3$, this would correspond to $(\ff^0\#\ff^1)\#(\ff^2\#\ff^3)$. So denote the set of vertices of $K(k+1)$ by $V(k)$. For each $v\in V(k)$, we obtain a functor  
\begin{equation}
  \#^v(\ff^0,\dots,\ff^k):
  Q_{\sum_{i=0}^{k}n_i-k}\longrightarrow\cc,
\end{equation}
given by gluing in the order corresponding to $v\in V(k)$. We write $d(I)=\sum_{i=0}^k n_i-k$, where $I=\{0,\dots,k\}$.  Now we return to our tree interpretation. We can represent any $v\in V(k)$ by an ordered path tree, with every branch $b$ numbered by $\text{ord}(b)\in \{0,\dots,k-1\}$. For example, we can consider $\ff^0,\ff^1$ and $\ff^2$ with $V(2)=\{((0,1),2),(0,(1,2))\}$, and obtain trees
\begin{equation}
    \begin{tikzcd}
\ff^2                         &  & \ff^2                         \\
\ff^1 \arrow[u, "1", no head] &  & \ff^1 \arrow[u, "0", no head] \\
\ff^0 \arrow[u, "0", no head] &  & \ff^0 \arrow[u, "1", no head]         
\end{tikzcd}
\end{equation}
Our desired homotopy from $\#^{v_1}(\ff^0,\ff^1,\ff^2)$ to $\#^{v_2}(\ff^0,\ff^1,\ff^2)$ should then correspond to the tree \begin{equation}
    \begin{tikzcd}
\ff^2                         &                       \\
\ff^1 \arrow[u, "0", no head] &   \\
\ff^0 \arrow[u, "0", no head] &          
\end{tikzcd}
\end{equation}
where one applies (a slight variant of) $\mu$ to this tree to obtain the sum of the two desired vertices.

\begin{prop}
    There is a homotopy $\Theta:\#^{v_1}(\ff^0,\ff^1,\ff^2) \implies \#^{v_2}(\ff^0,\ff^1,\ff^2)$.  I.e., 
    \begin{equation}
        \mu_1^\text{Fun}\Theta=\#^{v_1}(\ff^0,\ff^1,\ff^2)+\#^{v_2}(\ff^0,\ff^1,\ff^2)
    \end{equation}
\end{prop}
\begin{proof}
    The homotopy has in fact already been constructed. Indeed, the trees
    \begin{equation}
    \begin{tikzcd}
\ff^2                         &   &   & &                 \\
\ff^1 \arrow[u, "0", no head] &  &  \ff^2& &\ff^0\\
\ff^0 \arrow[u, "0", no head] &   &  & \ff^1\ar[ul,no head,"0"]\ar[ur,no head,swap,"0"]&     
\end{tikzcd}
\end{equation}
represent the same attachments after changing the root, so our homotopy is just that of $\Gamma$ constructed for parallel gluings. More formally, when $n_1\geq2$, we choose an object of $\ff^1$ away from the two attachment points as its output. Using the cyclic structure, we relabel the objects and apply compatible shifts if needed. We can then view $Q_{n_0}$ and $Q_{n_2}$ as glued in parallel to $Q_{n_1}$.
When $n_1=1$, we view $\ff^1$ as having an empty output and apply the same construction.
\end{proof}

So we want to generalise this to arbitrary $k$ to provide the following statement.
\begin{theorem}\label{Main theroem for series}
   For every $k\geq 2$, the family of functors \[
\#^v(\ff^0,\ldots,\ff^k):
Q_{\sum_{i=0}^k n_i-k}\longrightarrow\cc,
\qquad v\in V(k),
\]
admits a coherent family of $l$-homotopies, for $1\leq l\leq k-1$, that correspond to the fillings of the cells of the Stasheff associahedron $K(k+1)$.
\end{theorem}

We should fix some notation and explain what this statement means. Given $I=\{0,1,2\}$, there is a homotopy $\Theta^{012}:=\Theta$ which gives 
\begin{equation*}
    \mu_1^\text{Fun}\Theta^{012}=\#^{v_1}(\ff^0,\ff^1,\ff^2)+\#^{v_2}(\ff^0,\ff^1,\ff^2)
\end{equation*}
Here, $v_1$ and $v_2$ correspond to $((0,1),2)$ and $(0,(1,2))$, so we can write
\begin{equation*}
    \Theta^{(01)2}:=\#^{v_1}(\ff^0,\ff^1,\ff^2),
    \qquad
    \Theta^{0(12)}:=\#^{v_2}(\ff^0,\ff^1,\ff^2).
\end{equation*}
Thus we have
\begin{equation*}
    \mu_1^\text{Fun}\Theta^{012}=\Theta^{(01)2}+\Theta^{0(12)}.
\end{equation*}
If we go one order higher, with $I=\{0,1,2,3\}$, one will find that the corresponding $2$-homotopy should satisfy 
\begin{equation}\label{01234brackettheta}
\begin{aligned}
    \mu_1^\text{Fun}\Theta^{0123}
    &=\Theta^{(01)23}+\Theta^{0(12)3}+\Theta^{01(23)}\\
    &\quad+\Theta^{(012)3}+\Theta^{0(123)}+C^{0123}.
\end{aligned}
\end{equation}
Here, $C^{0123}$ consists of the correction terms for composing the face homotopies, as in the parallel setting.
Indeed, if we consider the orders of gluing together $4$ functors, there are $3$ gluings required, giving $3!=6$ orders. However, two of these orders result in the same functor. Consider the two trees
\begin{equation*}
    \begin{tikzcd}
 \ff^3   & \ff^3\\
\ff^2     \ar[u,no head,"1"] & \ff^2     \arrow[u, "0", no head]                                    \\
\ff^1 \arrow[u, "2", no head]  &\ff^1   \arrow[u, "2", no head]\\
\ff^0 \arrow[u, "0", no head]     &\ff^0      \arrow[u, "1", no head]
\end{tikzcd}
\end{equation*}
These are equal to one another: the gluings $\ff^2\#\ff^3$ and $\ff^0\# \ff^1$ happen away from one another, so the resulting functors will be identical. This can be seen in terms of bracketing, as $(01)(23)=(01)(23)$ will not depend on the order you place the brackets. Each of the five $\Theta$-terms on the RHS of equation~\eqref{01234brackettheta} corresponds to a homotopy between possible output functors. In particular, we have
\begin{align*}
\mu_1^{\mathrm{Fun}}\Theta^{(01)23}
    &=\Theta^{((01)2)3}+\Theta^{(01)(23)},\\
\mu_1^{\mathrm{Fun}}\Theta^{0(12)3}
    &=\Theta^{(0(12))3}+\Theta^{0((12)3)},\\
\mu_1^{\mathrm{Fun}}\Theta^{01(23)}
    &=\Theta^{(01)(23)}+\Theta^{0(1(23))},\\
\mu_1^{\mathrm{Fun}}\Theta^{(012)3}
    &=\Theta^{((01)2)3}+\Theta^{(0(12))3},\\
\mu_1^{\mathrm{Fun}}\Theta^{0(123)}
    &=\Theta^{0((12)3)}+\Theta^{0(1(23))}.
\end{align*}
Notice that these cancel in pairs. In general, a homotopy $\Theta^{0\dots k}$ should have boundary given by a sum over all ways of placing a single bracket inside the sequence $01\dots k$ with length of the sequence inside the bracket $l$ such that $1<l<k+1$, together with the corresponding correction terms $C^{0\dots k}$, i.e.,
\begin{equation*}
    \mu_1^{\text{Fun}}\Theta^{0\dots k}=\sum_{l=2}^{k}\sum_{r=0}^{k+1-l}\Theta^{0\dots(r\dots r+l-1)\dots k}+C^{0\dots k}.
\end{equation*}

\begin{rmk}
    We can alternatively view $\Theta$ as an equivalence class of an ordered tree as before. Given a tree containing a sub-tree of the form 
    \begin{equation}
        \begin{tikzcd}[column sep=small]
{} \arrow[r, no head, dashed] & \ff^{i-1} \arrow[rr, "x", no head] &  & \ff^i \arrow[rr, "y", no head] &  & \ff^{i+1} \arrow[rr, "z", no head] &  & \ff^{i+2} \arrow[r, no head, dashed] & {}
\end{tikzcd}
    \end{equation}
    with $y>x>z$, it is equivalent to the tree given by swapping the orders of the entire component with $x$ and the entire component with $z$. These swaps relabel the whole components while keeping the order, and any equal labels, within each component fixed. For example, let $\mathbf a$ and $\mathbf b$ record the labels in two independent components whose attachments occupy consecutive blocks of stages. Numbering these stages from $0$, suppose $\mathbf a$ uses $0,\dots,p-1$ and $\mathbf b$ uses $p,\dots,p+q-1$, with no other attachments at these stages. If the components are joined at a later stage $y\geq p+q$, the swap is
\[
    (\mathbf a,y,\mathbf b)\sim(\mathbf a+q,y,\mathbf b-p),
\]
where the additions apply to every label in the corresponding component; all other labels stay fixed. For instance,
\[
    ((0,1),3,(2))\sim((1,2),3,(0)).
\]
Thus, trees are equivalent if there is a sequence of such swaps. Furthermore, we need to impose that trees with a subtree as above with $y>x=z$ are equivalent to $0$. In general, if a tree has a subsequence of edge labels, $(a_{r-1},a_r,\dots,a_{r+s},a_{r+s+1})$ with $a_{r-1}=a_{r+s+1}$ and $a_{r+i}>a_{r-1}$ for $0\leq i\leq s$, then it is also declared equivalent to zero. We apply this zero relation recursively as we build the tree, contracting edges with smaller labels before considering the next attachments. For example, take the tree with edge labels $(1,2,0,1)$. We first glue $\ff^2$ to $\ff^3$ along the edge labelled $0$. The remaining tree has vertices $\ff^0,\ff^1,\ff^2\#\ff^3,\ff^4$ and edge labels $(1,2,1)$. The zero relation now applies, so the original tree is equivalent to $0$. This will account for trees whose boundary is a sum of equivalent trees.
The boundary equations can then be interpreted in terms of $\mu$ as in the parallel case, though we work with equivalence classes of trees. For example, the following trees are equivalent 
    \begin{equation}\begin{tikzpicture}[line cap=round,line join=round]
  \colorlet{seriesfinal}{blue!65!black}
  \foreach \yy/\aa/\bb/\cc/\dd/\name in {
    0/0/1/3/2/1,
    -1.6/0/2/3/1/2,
    -3.2/1/2/3/0/3}{
    \begin{scope}[yshift=\yy cm]
      \node at (-1.1,0) {$T_{\name}:$};
      \node (Fzero) at (0,0) {$\mathcal F^0$};
      \node (Fone) at (2,0) {$\mathcal F^1$};
      \node (Ftwo) at (4,0) {$\mathcal F^2$};
      \node (Fthree) at (6,0) {$\mathcal F^3$};
      \node (Ffour) at (8,0) {$\mathcal F^4$};
      \draw (Fzero)--node[above=3pt] {$\aa$} (Fone);
      \draw (Fone)--node[above=3pt] {$\bb$} (Ftwo);
      \draw[seriesfinal,line width=1.2pt]
        (Ftwo)--node[above=3pt] {$\cc$} (Fthree);
      \draw (Fthree)--node[above=3pt] {$\dd$} (Ffour);
    \end{scope}
  }
\end{tikzpicture}
    \end{equation}
To see $T_1\sim T_2$, first attach $\ff^0$ to $\ff^1$ at time $0$. We can then attach $\ff^0\#\ff^1$ to $\ff^2$ and attach $\ff^3$ to $\ff^4$ in either order, since both happen before the two components are joined at time $3$. Thus we can swap the labels $1$ and $2$.

To see $T_1\sim T_3$, we can instead attach $\ff^3$ to $\ff^4$ first, then attach $\ff^0$ to $\ff^1$, and finally attach $\ff^0\#\ff^1$ to $\ff^2$. The two components are still joined at time $3$, and the order of attachments within each component is unchanged.
\end{rmk}
\begin{prop}
    The definition of $\mu:\text{Tree}^{\ff^0}(I)\to \text{Tree}^{\ff^0}(I) $ in equation~\eqref{treemu} descends to equivalence classes of trees.
\end{prop}
\begin{proof}
    First, assume a tree $T$ is equivalent to $0$, i.e., $T\sim0$. We must show $\mu(T)\sim 0$. Indeed, assume $I=\{0,\dots,k\}$ and $T=(a_1,\dots,a_{r-1},a_r,a_{r+1},\dots,a_k)$ with $a_{r-1}=a_{r+1}$ and $a_r>a_{r+1}$. Then 
\begin{equation*}
\mu(a_1,\dots,a_{r-1},a_r,a_{r+1},\dots, a_{k})
=
\sum_{q=0}^{\max{\{a_1,\dots,a_k\}}}
\ \sum_{\varnothing\neq S\subsetneq\{j:a_j=q\}}
(a^{q,S}_1,\ldots,a_{r-1}^{q,S},a_r^{q,S},a_{r+1}^{q,S},\ldots,a^{q,S}_k),
\end{equation*}
recalling that 
\begin{equation*}
a^{q,S}_j=
\begin{cases}
a_j+1,&a_j>q\ \text{or}\ j\in S,\\
a_j,&\text{otherwise}.
\end{cases}   
\end{equation*}

The terms on the right-hand side either remain equivalent to zero or cancel in pairs of equivalent trees. Thus $\mu T\sim0$. The same argument applies to the longer subsequences in the zero relation.

Now, assume $T\sim T'$ by a single swapping and relabeling of blocks, so that
\[
    T=(\mathbf a,y,\mathbf b),
    \qquad T'=(\mathbf a+q,y,\mathbf b-p).
\]
If $\mu$ splits a label in $\mathbf a$, write $\mathbf a'$ for the resulting block. It now has $p+1$ distinct labels, and the corresponding terms in $\mu T$ and $\mu T'$ are related by
\[
    (\mathbf a',y+1,\mathbf b+1)
    \sim(\mathbf a'+q,y+1,\mathbf b-p).
\]
Similarly, if $\mu$ splits a label in $\mathbf b$, the resulting block $\mathbf b'$ has $q+1$ distinct labels, and we have
\[
    (\mathbf a,y+1,\mathbf b')
    \sim(\mathbf a+q+1,y+1,\mathbf b'-p).
\]
Splitting a label outside the two blocks leaves their numbers of stages unchanged, so the corresponding terms are again equivalent. Thus the terms of $\mu T$ and $\mu T'$ match in equivalent pairs, giving $\mu T\sim\mu T'$. The same holds if we take a sequence of such equivlences, thus in general.
\end{proof}

\begin{rmk}
    The above statement implies we can reinterpret the statement of Theorem~\ref{Main theroem for series} as saying that there exists a family of coherent homotopies that carry the combinatorial structure of $(\text{Tree}^{\ff^0}(I)/\sim,\mu)$.
\end{rmk}

Again, we will give the proof in parts, though we note that it will follow almost identically to the proof of parallel gluings. Indeed, we show a similar statement for the corresponding versions of $\Psi$ and $\Xi$ and use these to provide the desired homotopies. To begin, we show that the homology of the pushout categories behaves similarly to that of the parallel setting. Set $I=\{0,\dots,k\}$ and denote this category by 
\begin{equation}
    \bigvee_{i\in I}Q_{n_i}.
\end{equation}
\begin{prop}\label{homseriesQs}
  The homology classes of the morphism complexes $\bigvee_{i\in I}Q_{n_i}(x,y)$ are given by
\begin{itemize}
    \item $[1_x]$ when $x=y$.

    \item $[a^i_y]$ when $x=y-1$ inside the copy of $Q_{n_i}$,
    with indices modulo $n_i+1$.

    \item $[\mu_2(a^i_{n_i},a^{i+1}_{j_{i+1}+1})]$ when
    $x=(n_i-1)_{Q_{n_i}},\,
        y=(j_{i+1}+1)_{Q_{n_{i+1}}},\,
         0\leq i<k.$

    \item $[\mu_2(a^{i+1}_{j_{i+1}},a^i_0)]$ when $x=(j_{i+1}-1)_{Q_{n_{i+1}}},\,
         y=0_{Q_{n_i}},\,
        0\leq i<k.$
   \item
$\bigl[\mu_2\bigl(\cdots\mu_2(a^i_{n_i},
a^{i+1}_{n_{i+1}})\cdots,
a^{i+l}_{j_{i+l}+1}\bigr)\bigr]$
when
\[
    x=(n_i-1)_{Q_{n_i}},\qquad
    y=(j_{i+l}+1)_{Q_{n_{i+l}}},
\]
where $l\geq2$, $i,i+1,\ldots,i+l\in I$, and
\[
    j_r=n_r-1
    \qquad\text{for every }i+1\leq r\leq i+l-1.
\]

\item
$\bigl[\mu_2\bigl(\cdots\mu_2(a^{i+l}_{j_{i+l}},
a^{i+l-1}_0)\cdots,a^i_0\bigr)\bigr]$
when
\[
    x=(j_{i+l}-1)_{Q_{n_{i+l}}},\qquad
    y=0_{Q_{n_i}},
\]
where $l\geq2$, $i,i+1,\ldots,i+l\in I$, and
\[
    j_r=0
    \qquad\text{for every }i+1\leq r\leq i+l-1.
\]
    \item Trivial otherwise.
\end{itemize}
All local indices are taken modulo $n_r+1$ in the corresponding factor $Q_{n_r}$.
\end{prop}
\begin{proof}
    The cases where the homology is trivial follow by applying the argument from the proof of Lemma~\ref{lemphihomology} at each successive attachment, keeping track of where the gluings take place. Similarly, the proof of points 1 through 4 follows in the same way. What is left are the terms that appear as iterated $\mu_2$'s. Here one might already note the difference with parallel gluings: we can get non-zero homology between objects that live multiple gluings away. One might further notice that the associahedra are appearing here. We can check that
    \[
        \mu_1\mu_2(x,y)=\mu_2(\mu_1x,y)+\mu_2(x,\mu_1y).
    \]
Since the original arrows are cycles, iterating this identity gives
    \[
        \mu_1\mu_2\bigl(\cdots\mu_2(a^i_{n_i},
        a^{i+1}_{n_{i+1}})\cdots,
        a^{i+l}_{j_{i+l}+1}\bigr)=0.
    \]
    Moreover, as in the proof of Lemma~\ref{lemphihomology}, one can check that these iterated composites are not boundaries for degree reasons, so they represent non-trivial homology classes. A similar check is made for the $\bigl[\mu_2\bigl(\cdots\mu_2(a^{i+l}_{j_{i+l}},
a^{i+l-1}_0)\cdots,a^i_0\bigr)\bigr]$ case.
\end{proof}

\subsubsection{The higher $\Psi$s in series}
We consider the various orders in which one can attach our $Q_{n_i}$'s in series. {For each $v\in V(k)$, we have the corresponding pair
\[
    \Phi^v:Q_{d(I)}\longrightarrow\bigvee_{i\in I}Q_{n_i},
    \qquad
    E_v:\bigvee_{i\in I}Q_{n_i}\longrightarrow\cc.
\]
Here the $\Psi$'s compare the extensions $E_v$, while the $\Xi$'s below compare the $\Phi^v$'s.}
The claim is that there exists $\Psi^{0\dots k}$ satisfying the corresponding higher homotopy equation. We retain the face sum
\[
    D^{0\dots k}:=\sum_{l=2}^{k}\sum_{r=0}^{k+1-l}
    \Psi^{0\dots(r\dots r+l-1)\dots k}.
\]
{As in the parallel case, we use the completed boundary
\[
    \widehat D^{0\dots k}=D^{0\dots k}+C^{0\dots k},
    \qquad
    \mu_1^{\mathrm{Fun}}\Psi^{0\dots k}=\widehat D^{0\dots k},
\]
where $C^{0\dots k}$ contains the higher composition terms for the chosen face composites. We fix the endpoint extensions $L,R$ for this filling.}

We have already (informally) shown that there is a $\Psi^{012}$ with
\[
    \mu_1^{\mathrm{Fun}}\Psi^{012}=\Psi^{(01)2}+\Psi^{0(12)}.
\]
So assume that we have constructed the lower homotopies for all shorter consecutive subsequences. We show that we can construct $\Psi^{0\dots k}$. Indeed, we can again set $\Psi^{0\dots k}_l=0$ for all $l\neq1$. Thus we need to build a suitable $\Psi^{0\dots k}_1\mu_s(x_1,\dots,x_s)$ for composable inputs in $\bigvee_{i\in I}Q_{n_i}$, with $s\geq2$.

The proof follows similarly to that for $\Psi^I$. {We set $\Psi^{0\dots k}_1=0$ on the original arrows and units, and induct on the total number of operation vertices in the trees representing the inputs $x_1,\dots,x_s$.} To begin, if every $x_v$ is an original arrow, we have
\[
    \Psi^{0\dots k}_1\mu_s(x_1,\dots,x_s)
    ={\widehat D^{0\dots k}_s(x_1,\dots,x_s)}.
\]
More generally, once $\Psi^{0\dots k}_1$ is defined on the input trees, we set
\begin{equation}\label{eq:seriespsirecursion}
\begin{aligned}
    &\Psi^{0\dots k}_1\mu_s(x_1,\dots,x_s)
    ={\widehat D^{0\dots k}_s(x_1,\dots,x_s)}\\
    &\quad+\sum_{v=1}^{s}\sum
    \mu_{\alpha+1+\beta}\Bigl(
    L_{p_1}(x_1,\dots,x_{p_1}),\dots,
    L_{p_\alpha}(x_{v-p_\alpha},\dots,x_{v-1}),\\
    &\qquad \Psi^{0\dots k}_1(x_v),
    R_{t_1}(x_{v+1},\dots,x_{v+t_1}),\dots,
    R_{t_\beta}(x_{s-t_\beta+1},\dots,x_s)\Bigr).
\end{aligned}
\end{equation}
The remaining check is to verify
\begin{equation}
\begin{aligned}
    &\mu_1\Psi^{0\dots k}_1\mu_s(x_1,\dots,x_s)
    +\Psi^{0\dots k}_1\mu_1\mu_s(x_1,\dots,x_s)\\
    &\qquad={\widehat D^{0\dots k}_1}\bigl(\mu_s(x_1,\dots,x_s)\bigr).
\end{aligned}
\end{equation}
But again, this follows in the same way as for $\Psi^I$, using
{$\mu_1^{\mathrm{Fun}}\widehat D^{0\dots k}=0$}.
The combinatorics are now those of associahedra: every codimension-two face lies in exactly two codimension-one faces (see, for example, \cite{MTTV}). {As before, the chosen coherent face compositions and their correction terms are included in this cancellation.}

\subsubsection{The higher $\Xi$s in series}
Now we show the corresponding statement for $\Xi^{0\dots k}$. Again, we assume that we have constructed the lower homotopies on all shorter consecutive subsequences, with $\Xi^{012}$ as the base case. We want to build $\Xi^{0\dots k}_s(x_1,\dots,x_s)$ for composable sequences in $Q_{d(I)}$. {Here $L,R$ denote the chosen endpoint $\Phi$-functors. We use the completed boundary
\[
    \widehat B^{0\dots k}=B^{0\dots k}+C_\Xi^{0\dots k},
    \qquad
    \mu_1^{\mathrm{Fun}}\Xi^{0\dots k}=\widehat B^{0\dots k},
\]
where $C_\Xi^{0\dots k}$ contains the corresponding higher composition terms.}
We retain the notation
\[
    B^{0\dots k}:=\sum_{l=2}^{k}\sum_{r=0}^{k+1-l}
    \Xi^{0\dots(r\dots r+l-1)\dots k}.
\]
{We take $\Xi^{0\dots k}_0=0$ and $\Xi^{0\dots k}_s(\dots,1,\dots)=0$.} In general, we can reorganise the homotopy equation as
\begin{equation}\label{eq:seriesxirecursion}
\begin{aligned}
    &\mu_1\Xi^{0\dots k}_s(x_1,\dots,x_s)
    ={\widehat B^{0\dots k}_s(x_1,\dots,x_s)}\\
    &\quad+\sum_{e=1}^{s-1}\sum_{v=1}^{s-e+1}\sum
    \mu_{\alpha+1+\beta}\Bigl(
    L_{p_1}(x_1,\dots,x_{p_1}),\dots,
    L_{p_\alpha}(x_{v-p_\alpha},\dots,x_{v-1}),\\
    &\qquad \Xi^{0\dots k}_e(x_v,\dots,x_{v+e-1}),
    R_{t_1}(x_{v+e},\dots,x_{v+e+t_1-1}),\dots,
    R_{t_\beta}(x_{s-t_\beta+1},\dots,x_s)\Bigr).
\end{aligned}
\end{equation}
Again, $\mu_l(x_{r+1},\dots,x_{r+l})=0$ inside $Q_{d(I)}$, except when $l=d(I)+1$, in which case it is a unit. We then argue in the same manner as in the parallel case. For $s$ large enough, and for now excluding the cases of non-trivial homology classes from Proposition~\ref{homseriesQs}, the RHS must be a boundary as one can verify it is a cycle. 

For unit inputs, we have already taken $\Xi^{0\dots k}_s(\dots,1,\dots)=0$. We also need to consider the homology classes
\[
\bigl[\mu_2\bigl(\cdots\mu_2(a^i_{n_i},
a^{i+1}_{n_{i+1}})\cdots,a^{i+l}_{j_{i+l}+1}\bigr)\bigr]
\]
and
\[
\bigl[\mu_2\bigl(\cdots\mu_2(a^{i+l}_{j_{i+l}},
a^{i+l-1}_0)\cdots,a^i_0\bigr)\bigr].
\]
For these, we note that for $k\geq3$ the degree argument also gives the result, as in the parallel case. Indeed, one finds that the degree of the RHS of~\eqref{eq:seriesxirecursion} is 
\[\sum_{j=1}^s|x_j|+s+k-3\].

\begin{rmk}
Similarly to the parallel case, we can choose our corresponding output functors to satisfy Theorem~\ref{Main theroem for series}.
\end{rmk}

\subsection{Extending to mixed gluings}\label{subsec:all-gluings}

So we now want to put both series and parallel gluings together. Indeed, we can consider gluings of functors that correspond to a mix of the two, for example
\begin{equation}
    \begin{tikzcd}
        &                                                                & \ff^3                    \\
\ff^{0} & \ff^{1}                                                        & \ff^2 \arrow[u, no head] \\
        & \gg \arrow[lu, no head] \arrow[u, no head] \arrow[ru, no head] &
\end{tikzcd}
\end{equation}
This is where viewing the combinatorics in terms of ordered trees comes in useful. Indeed, we can label all the branches as we did before and use the same $\mu$ on equivalence classes of such trees. We use the equivalence from the series setting, as no additional relations are needed for parallel trees. We should make this more precise.

We consider a root $\ff^\emptyset:Q_{n_\emptyset}\to\cc$ and a set $I^\emptyset\subseteq\{0,\dots,n_\emptyset-1\}$. For every $\alpha^0\in I^\emptyset$, we consider a functor $\ff^{(\alpha^0)}:Q_{n_{(\alpha^0)}}\to\cc$ with $\ff^{(\alpha^0)}_0(n_{(\alpha^0)})=\ff^\emptyset_0(\alpha^0)$. We then set $I^{(\alpha^0)}\subseteq\{0,\dots,n_{(\alpha^0)}-1\}$, so that for every $\alpha^1\in I^{(\alpha^0)}$ there is a functor $\ff^{(\alpha^0,\alpha^1)}:Q_{n_{(\alpha^0,\alpha^1)}}\to\cc$ with $\ff^{(\alpha^0,\alpha^1)}_0(n_{(\alpha^0,\alpha^1)})=\ff^{(\alpha^0)}_0(\alpha^1)$. Repeating this process gives a family of functors
\begin{equation}
    \ff^{(\alpha^0,\dots,\alpha^s)}:Q_{n_{(\alpha^0,\dots,\alpha^s)}}\to\cc,
\end{equation}
where $\alpha^i\in I^{(\alpha^0,\dots,\alpha^{i-1})}\subseteq\{0,\dots,n_{(\alpha^0,\dots,\alpha^{i-1})}-1\}$, with
\begin{equation}
    \ff^{(\alpha^0,\dots,\alpha^j)}_0(n_{(\alpha^0,\dots,\alpha^j)})
    =\ff^{(\alpha^0,\dots,\alpha^{j-1})}_0(\alpha^j).
\end{equation}
Here, the empty sequence denotes the root. 

In terms of ordered trees, one views the functor $\ff^{(\alpha^0,\dots,\alpha^s)}$ as the vertex we reach by starting at the root and taking the branch at input $\alpha^0$, then the branch at input $\alpha^1$ of the next vertex, and so on. For example, consider the following tree:
\begin{equation*}
    \begin{tikzcd}
. &                                           & X                     &                                                              &                                          &   \\
  & . \arrow[ru, no head] \arrow[lu, no head] &                       &                                                              & Y                                        & . \\
  &                                           & . \arrow[lu, no head] & .                                                            & . \arrow[u, no head] \arrow[ru, no head] &   \\
  &                                           &                       & . \arrow[lu, no head] \arrow[u, no head] \arrow[ru, no head] &                                          &
\end{tikzcd}
\end{equation*}
The vertex $X$ corresponds to a functor $\ff^{(0,0,1)}$, and $Y$ to a functor $\ff^{(2,0)}$. We can encode the order in which we attach the vertices by labelling the edges. We also allow multiple edges to carry the same number; these repetitions correspond to homotopies. For example,

\begin{equation*}
    \begin{tikzcd}
. &                                           & .                    &                                                              &                                          &   \\
  & . \arrow[ru, no head,swap,"1"] \arrow[lu, no head,"3"] &                       &                                                              & .                                       & . \\
  &                                           & . \arrow[lu, no head,"0"] & .                                                            & . \arrow[u, no head,"5"] \arrow[ru, no head,swap,"6"] &   \\
  &                                           &                       & . \arrow[lu, no head,"2"] \arrow[u, no head,"7"] \arrow[ru, no head,swap,"4"] &                                          &
\end{tikzcd}
\end{equation*}
will correspond to a specific output functor $\#^\tau(\ff^{\mathcal{I}}):Q_{d(\mathcal{I})}\to \cc$ with $\mathcal I$ denoting the set of vertex addresses, including $\emptyset$, and
\begin{equation}
    d(\mathcal{I})=\sum_{\mathbf{\alpha}\in\mathcal{I}}n_{\mathbf{\alpha}} - |\mathcal{I}|+1.
\end{equation}
Here, $\tau\in \overline{\text{Tree}}(\mathcal{I})$ with $\overline{\text{Tree}}(\mathcal{I})=\text{Tree}^{\ff^{\emptyset}}(\mathcal{I})/\sim$.
\begin{theorem}\label{thm:all-gluings}
   Given a collection of $\ff^{\mathbf{\alpha}}:Q_{n_\mathbf{\alpha}}\to \cc$ with $\mathbf{\alpha}\in \mathcal{I}$ for some fixed $\mathcal{I}$ as above, there is a correspondence between all orders of gluings of functors, along with a coherent family of higher homotopies, and $\overline{\text{Tree}}(\mathcal{I})$, with boundary relations described by $\mu:\overline{\text{Tree}}(\mathcal{I})\to \overline{\text{Tree}}(\mathcal{I})$ together with the composition terms below. More precisely, to each tree class $\tau$
of degree zero we assign the functor
\[
    \Theta^\tau:=\#^\tau(\ff^{\mathcal I}):
    Q_{d(\mathcal I)}\longrightarrow\cc,
\]
and to each tree class $\tau$ of degree $q>0$ we assign a
$q$-dimensional homotopy $\Theta^\tau$. These assignments respect
the equivalence and zero relations and, for $q>0$, satisfy
\[
    \mu_1^{\mathrm{Fun}}\Theta^\tau
    =\Theta^{\mu\tau}+C^\tau.
\]
Here, $\Theta^{\mu\tau}=\sum_i\Theta^{\tau_i}$ when $\mu\tau=\sum_i\tau_i$, and $C^\tau$ consists of the correction terms for the chosen compositions of lower-degree face homotopies, as in the parallel and series cases.
\end{theorem}

\begin{rmk}
    We have already provided the proofs for when $\Theta$ corresponds to a single parallel tree, or for when it is a single series tree. Any underlying tree can be built by a sequence of gluings of such trees. We first describe this gluing and its counterpart on the functor side; the higher homotopies for arbitrary edge labellings are then constructed by induction below. In particular, we need a map
    \begin{equation*}
        \text{Glue}:\overline{\text{Tree}}(\mathcal{I})\times \overline{\text{Tree}}(\mathcal{J})\to \overline{\text{Tree}}(\mathcal I*_{\mathbf\alpha}\mathcal J)
    \end{equation*}
    such that $\mu \text{Glue}(-,-)=\text{Glue}(\mu(-),-)+\text{Glue}(-,\mu(-))$ where $\mathbf{\alpha}\in \mathcal{I}$ is a chosen vertex and  $\mathcal{I}*_{\mathbf{\alpha}}\mathcal{J}$ is the resulting set of all sequences when you glue the root of the generic tree corresponding to $\mathcal{J}$ to $\mathcal{I}$ at $\mathbf{\alpha}$.
\end{rmk}

\begin{prop}\label{tree-differential}
    The map $\mu$ is a differential on ordered trees, i.e., $\mu^2=0$, and hence also on $\overline{\text{Tree}}(\mathcal I)$.
\end{prop}
\begin{proof}
    We work with labelled representatives before quotienting. If every branch of $\tau$ has a different order, then we take $\mu\tau=0$, i.e., every degree-zero tree is a cycle. The degree of a general tree is the sum of the multiplicities minus one for each occupied branch order. It is clear that $\mu$ lowers the degree by one, so for any degree-$1$ tree we automatically obtain $\mu^2\tau=0$.

    Next, let's assume $\text{deg}(\tau)=q>0$. We can write this as
\begin{equation}
\deg\tau=\sum_{m=0}^{\max\{\text{ord}(a):a\in \tau\}}\bigl(|m_\tau|-1\bigr),
\end{equation}
where we have written
\begin{equation}
    |m_\tau|:=\left|\{a\in \tau : \text{ord}(a)=m\}\right|.
\end{equation}
From the formula for $\mu$, we find that the number of terms in $\mu\tau$, before quotienting by the equivalence, is given by
\begin{equation}
   \sum_{m=0}^{\max\{\text{ord}(a):a\in \tau\}}\sum_{i=1}^{|m_\tau|-1}{|m_\tau|\choose i}.
\end{equation}
It is a standard result that
\begin{equation*}
    \sum_{i=0}^n{n\choose i}=2^n,
\end{equation*}
so we find
\begin{align}
    \sum_m\sum_{i=1}^{|m_\tau|-1}{|m_\tau|\choose i}
    &=\sum_m\left(2^{|m_\tau|}-{|m_\tau|\choose0}-{|m_\tau|\choose|m_\tau|}\right)\notag\\
    &=2\sum_m\left(2^{|m_\tau|-1}-1\right),
\end{align}
where $m$ runs from $0$ to $\max\{\text{ord}(a):a\in\tau\}$. This is even. In particular, the number of terms in $\mu^2\tau$, counted with multiplicity, is even. We claim more precisely that each resulting tree appears exactly twice. Explicitly, if $\mu\tau=\sum\tau_i$, each tree in $\mu^2\tau$ appears in $\mu\tau_{s}$ and $\mu\tau_{s'}$ for two distinct terms $\tau_s,\tau_{s'}$.

For example, consider $\deg(\tau)=2$ with $|m_\tau|=3$ for some $m$. Away from these three edges, every resulting tree in $\mu^2\tau$ looks the same. There are a total of $2\cdot(2\cdot(2^2-1))=12$ terms. The three edges that originally had order $m$ now have orders $m,m+1,m+2$, in some order. Thus there are only $3!=6$ distinct resulting trees, and each occurs twice. Indeed, label the three edges $(a,b,c)$ and write $(0,2,1)$ to mean that we have raised the order of $a$ by zero, $b$ by two and $c$ by one. This occurs in the following two ways:
\begin{align*}
    (0,0,0)&\mapsto(0,1,1)\mapsto(0,2,1),\\
    (0,0,0)&\mapsto(0,1,0)\mapsto(0,2,1).
\end{align*}

More generally, label the edges at every order $m$ with $|m_\tau|>1$, and list them as
\begin{equation}
    (a^{m_0}_0,\dots,a^{m_0}_{k_0},
    a^{m_1}_0,\dots,a^{m_1}_{k_1},\dots,
    a^{m_l}_0,\dots,a^{m_l}_{k_l}).
\end{equation}
Applying $\mu$ to $\tau$, we choose an order $m^1=m_i\in\{m_0,\dots,m_l\}$ and a subset
\[
    \emptyset\neq R^1\subsetneq\{a^{m_i}_0,\dots,a^{m_i}_{k_i}\}.
\]
We raise the order of every edge in $R^1$ by one, and also raise the order of every edge with order greater than $m^1$ by one.

If we repeat this, there are three cases. Either the second order we pick is $m^2=m^1+1$, in which case $\emptyset\neq R^2\subsetneq R^1$; or $m^2=m^1$, in which case $R^2$ is a nonempty proper subset of the edges left at order $m^1$; or $m^2\notin\{m^1,m^1+1\}$.

Assume first that $m^2=m^1+1$. Every edge originally above order $m^1$ has been raised by two, every edge in $R^2$ has been raised by two, and every edge in $R^1\setminus R^2$ has been raised by one.

We realise that this exact tree appears again by choosing the orders $\bar m^1=\bar m^2=m^1$, with
\[
    \bar R^1=R^2,
    \qquad
    \bar R^2=R^1\setminus R^2.
\]
Indeed, every edge in $\bar R^1=R^2$ is raised by two and every edge in $\bar R^2=R^1\setminus R^2$ is raised by one.

Next, assume $m^2\notin\{m^1,m^1+1\}$. This case is more straightforward. If $m^2<m^1$, the same tree appears by first choosing $\bar m^1=m^2$ with $\bar R^1=R^2$, and then $\bar m^2=m^1+1$ with $\bar R^2=R^1$. Similarly, if $m^2>m^1+1$, we choose $\bar m^1=m^2-1$ and $\bar m^2=m^1$, again with $\bar R^1=R^2$ and $\bar R^2=R^1$. Thus every resulting tree in $\mu^2\tau$ appears twice. Since we work modulo two, the sum is zero.
\end{proof}
\begin{rmk}
The idea of composing trees by substitution and by what they call `grafting' also appears in the work of Loday and Vallette \cite{LV}.
\end{rmk}

We can define $\text{Glue}(-,-)$ in the obvious way. If we have $(\tau,\tau')\in\overline{\text{Tree}}(\mathcal{I})\times \overline{\text{Tree}}(\mathcal{J}) $ we construct a new tree $\text{Glue}_{\mathbf\alpha}(\tau,\tau')\in \overline{\text{Tree}}(\mathcal I*_{\mathbf\alpha}\mathcal J)$ by identifying the root of $\tau'$ with $\mathbf\alpha\in\tau$. We then correct the branch ordering by setting
\begin{equation}
    \text{ord}_{\text{Glue}_{\mathbf\alpha}(\tau,\tau')}(b):=\begin{cases}
        \text{ord}_{\tau}(b) & b\in \tau\\
        \text{ord}_{\tau'}(b)+\max_{e\in \tau}\{\text{ord}(e)\}+1& b\in \tau'.
    \end{cases}
\end{equation}
For example, we have
\begin{equation*}
    \begin{tikzcd}[column sep=small,row sep=small]
             &   &                                                                              &        &     &                           &   &   &                             &                           & .                                                  \\
             &   &                                                                              &        &     & .                         &   &   &                             &                           & . \arrow[u, "2", no head]                          \\
\text{Glue}_{\mathbf\alpha}( & . & .                                                                            & \alpha & {,} & . \arrow[u, "0", no head] & ) & = & . \arrow[rd, "0"', no head] & .                         & . \arrow[ld, "1", no head] \arrow[u, "3", no head] \\
             &   & . \arrow[lu, "0", no head] \arrow[u, "1", no head] \arrow[ru, "1"', no head] &        &     & . \arrow[u, "1", no head] &   &   &                             & . \arrow[u, "1", no head] &
\end{tikzcd}
\end{equation*}
Here, $\mu(\tau')=0$, as $\tau'$ has degree zero. We find that $\mu(\tau)$ is
\begin{equation*}
    \begin{tikzcd}
. & .                                                                            & \alpha & + & . \arrow[rd, "0"', no head] & .                                                   & \alpha \\
  & . \arrow[lu, "0", no head] \arrow[u, "1", no head] \arrow[ru, "2"', no head] &        &   &                             & . \arrow[u, "2", no head] \arrow[ru, "1"', no head] &
\end{tikzcd}
\end{equation*}

Then $\text{Glue}_{\mathbf\alpha}(\mu(\tau),\tau')$ gives
\begin{equation*}
    \begin{tikzcd}
  &                                                                              & .                              &   &                             &                                                     & .                              \\
  &                                                                              & . \arrow[u, "3", no head]      &   &                             &                                                     & . \arrow[u, "3", no head]      \\
. & .                                                                            & \alpha \arrow[u, "4", no head] & + & . \arrow[rd, "0"', no head] & .                                                   & \alpha \arrow[u, "4", no head] \\
  & . \arrow[lu, "0", no head] \arrow[u, "1", no head] \arrow[ru, "2"', no head] &                                &   &                             & . \arrow[u, "2", no head] \arrow[ru, "1"', no head] &
\end{tikzcd}
\end{equation*}
Then $\mu\text{Glue}_{\mathbf\alpha}(\tau,\tau')$ is given by
\begin{equation*}
\begin{tikzcd}
  &                                                                              & .                              &   &                             &                                                     & .                              \\
  &                                                                              & . \arrow[u, "3", no head]      &   &                             &                                                     & . \arrow[u, "3", no head]      \\
. & .                                                                            & \alpha \arrow[u, "4", no head] & + & . \arrow[rd, "0"', no head] & .                                                   & \alpha \arrow[u, "4", no head] \\
  & . \arrow[lu, "0", no head] \arrow[u, "1", no head] \arrow[ru, "2"', no head] &                                &   &                             & . \arrow[u, "2", no head] \arrow[ru, "1"', no head] &
\end{tikzcd}
\end{equation*}
which agrees with $\text{Glue}(\mu\tau,\tau')$.

\begin{prop}
    The map $\text{Glue}_{\mathbf\alpha}(-,-):\overline{\text{Tree}}(\mathcal{I})\times \overline{\text{Tree}}(\mathcal{J})\to \overline{\text{Tree}}(\mathcal I*_{\mathbf\alpha}\mathcal J)$ satisfies
    \begin{equation}
        \mu \text{Glue}_{\mathbf\alpha}(-,-)=\text{Glue}_{\mathbf\alpha}(\mu(-),-)+\text{Glue}_{\mathbf\alpha}(-,\mu(-))
    \end{equation}
\end{prop}
\begin{proof}
Choose any pair $(\tau,\tau')$. Both compositions give the same underlying tree shape. We need to verify that every labelled tree obtained by the clockwise path in
\begin{equation}
    \begin{tikzcd}[column sep=small]
        \overline{\text{Tree}}(\mathcal I)\otimes\overline{\text{Tree}}(\mathcal J)
        \arrow[d,"\mu\otimes1+1\otimes\mu"']
        \arrow[r,"\text{Glue}_{\mathbf\alpha}"]
        & \overline{\text{Tree}}(\mathcal I*_{\mathbf\alpha}\mathcal J)\arrow[d,"\mu"]\\
        \overline{\text{Tree}}(\mathcal I)\otimes\overline{\text{Tree}}(\mathcal J)
        \arrow[r,"\text{Glue}_{\mathbf\alpha}"]
        & \overline{\text{Tree}}(\mathcal I*_{\mathbf\alpha}\mathcal J)
    \end{tikzcd}
\end{equation}
matches a unique term obtained by the other path.
Note that we take the spans of tree classes over $\mathbb Z/2\mathbb Z$ and extend $\text{Glue}_{\mathbf\alpha}$ bilinearly.

First, the number of terms is the same. Write $m=\max\text{ord}_\tau$ and $n=\max\text{ord}_{\tau'}$, taking the maximum to be $-1$ for a tree with no edges. By definition,
\[
    \text{ord}_{\text{Glue}_{\mathbf\alpha}(\tau,\tau')}(b)=
    \begin{cases}
        \text{ord}_\tau(b),&b\in\tau,\\
        \text{ord}_{\tau'}(b)+m+1,&b\in\tau'.
    \end{cases}
\]
Thus every order value in $\tau$ has the same multiplicity in the glued tree, and every order value $q$ in $\tau'$ becomes $q+m+1$, with the same multiplicity. The number of terms before quotienting is therefore the same.

Next we check the labels. Write $\mu\tau=\sum_i\tau_i$ and $\mu\tau'=\sum_j\tau'_j$. We have
\[
\begin{aligned}
    &\text{Glue}_{\mathbf\alpha}(\mu\tau,\tau')
      +\text{Glue}_{\mathbf\alpha}(\tau,\mu\tau')\\
    &\qquad=\sum_i\text{Glue}_{\mathbf\alpha}(\tau_i,\tau')
      +\sum_j\text{Glue}_{\mathbf\alpha}(\tau,\tau'_j).
\end{aligned}
\]
We need to realise each of these as a term in $\mu\text{Glue}_{\mathbf\alpha}(\tau,\tau')$. Write $(a_1,\dots,a_k)$ for the edge labels of $\text{Glue}_{\mathbf\alpha}(\tau,\tau')$. Using the notation of~\eqref{treemu}, we can split the sum as
\begin{align*}
    \mu\text{Glue}_{\mathbf\alpha}(\tau,\tau')
    &=\sum_{q=0}^{m+n+1}
      \sum_{\emptyset\neq S\subsetneq\{j:a_j=q\}}
      (a^{q,S}_1,\dots,a^{q,S}_k)\\
    &=\sum_{q=0}^{m}
      \sum_{\emptyset\neq S\subsetneq\{j:a_j=q\}}
      (a^{q,S}_1,\dots,a^{q,S}_k)\\
    &\quad+\sum_{q=m+1}^{m+n+1}
      \sum_{\emptyset\neq S\subsetneq\{j:a_j=q\}}
      (a^{q,S}_1,\dots,a^{q,S}_k).
\end{align*}
In the first sum, we split an order value in $\tau$. This also raises every label from $\tau'$ by one, exactly as the new maximum $m+1$ in $\tau_i$ requires. These are the terms $\text{Glue}_{\mathbf\alpha}(\tau_i,\tau')$. In the second sum, we split an order value in $\tau'$ and leave $\tau$ unchanged. These are the terms $\text{Glue}_{\mathbf\alpha}(\tau,\tau'_j)$.
\end{proof}

We now show that we can choose a family of functors and higher homotopies whose gluing is controlled by $\text{Glue}$ and $\mu$. We keep the representing trees as part of the data and define
\begin{equation}
    \overline{\text{Glue}_{\mathbf\alpha}}(\Theta^\tau,\Theta^{\tau'})
    :=\Theta^{\text{Glue}_{\mathbf\alpha}(\tau,\tau')}.
\end{equation}
For a glued tree of positive degree, the boundary equation then reads
\begin{equation}
\begin{aligned}
    \mu_1^{\text{Fun}}\Theta^{\text{Glue}_{\mathbf\alpha}(\tau,\tau')}
    &=\Theta^{\text{Glue}_{\mathbf\alpha}(\mu\tau,\tau')}
      +\Theta^{\text{Glue}_{\mathbf\alpha}(\tau,\mu\tau')}\\
    &\quad+C^{\text{Glue}_{\mathbf\alpha}(\tau,\tau')}.
\end{aligned}
\end{equation}
Tree sums in the superscripts denote the corresponding sums of $\Theta$-terms, and $C$ contains the composition corrections as above. Thus we can describe the compositions and higher-homotopy relations using their representing trees.
These constructions will follow from the following statements.

\begin{rmk}
Fix a tree $\tau\in\overline{\text{Tree}}(\mathcal I)$ with $|m_\tau|=1$ for every occupied order $m$, and the corresponding family of functors $\ff^{\mathbf\alpha}:Q_{n_{\mathbf\alpha}}\to\cc$ for $\mathbf\alpha\in\mathcal I$. Then there is a degree-zero functor $\Theta^\tau:Q_{d(\mathcal I)}\to\cc$, given by a sequence of pushouts and pullbacks in the order corresponding to $\tau$. That is, the map $\tau\mapsto\Theta^\tau$ is defined in degree zero.
\end{rmk}

\begin{prop}\label{treehomcom}
    The homology of the morphism complexes
    $(\bigvee_{\mathbf\alpha\in\mathcal I}Q_{n_{\mathbf\alpha}})(x,y)$
    is generated by the classes listed in the parallel and series cases (Corollary~\ref{corhom} and Proposition~\ref{homseriesQs}), along with additional classes that are iterated compositions of such classes. 
\end{prop}
\begin{proof}
    Any general tree can be built by successively gluing trees in series and in parallel. Applying the reduction from the proof of Lemma~\ref{lemphihomology} at each attachment gives the listed representatives and their mixed composites, just as in the parallel and series proofs. The surviving representatives are non-trivial for the same degree reasons.
\end{proof}

\begin{prop}
    For every ordered tree $\tau\in\overline{\text{Tree}}(\mathcal I)$ of degree $q\geq0$, there is a $\Psi^\tau:\bigvee_{\mathbf\alpha\in\mathcal I}Q_{n_{\mathbf\alpha}}\to\cc$ of degree $q$. In degree zero this is a functor, and in positive degree it is a higher homotopy. For $q>0$, it satisfies
    \begin{equation}
        \mu_1^{\text{Fun}}\Psi^\tau
        =\Psi^{\mu\tau}+C_\Psi^\tau,
        \qquad
        \Psi^{\mu\tau}:=\sum_i\Psi^{\tau_i},
    \end{equation}
    where $\mu\tau=\sum_i\tau_i$.
\end{prop}
Here, $C_\Psi^\tau$ contains the correction terms for composing the face homotopies, as in the parallel setting.
\begin{proof}
The proof is almost identical to the parallel and series cases, but with the full combinatorics of ordered trees. We induct on the degree of $\tau$. The degree-zero extensions are obtained by successively applying the pushout construction. For positive degree, we use the same recursive formula as in~\eqref{higherPsieq}, taking $\Psi^\tau_l=0$ for $l\neq1$ and replacing the boundary term by
\[
    \sum_i\Psi^{\tau_i}+C_\Psi^\tau.
\]
We apply the same construction of face compositions and correction terms as in the parallel and series cases, now using the combinatorial structure of the trees to collect the terms. The same cancellation, using $\mu^2\tau=0$ (Proposition~\ref{tree-differential}), gives the unary equation
\[
    \mu_1\Psi^\tau_1(x)+\Psi^\tau_1(\mu_1x)
    =\sum_i\Psi^{\tau_i}_1(x)+(C_\Psi^\tau)_1(x).
\]
\end{proof}

\begin{rmk}
Every ordered series tree and every ordered parallel tree can be viewed as a general ordered tree. The combinatorics restrict to the cubical setting in parallel and the associahedral setting in series.
\end{rmk}

\begin{prop}
    For every ordered tree $\tau\in\overline{\text{Tree}}(\mathcal I)$ of degree $q\geq0$, there is a $\Xi^\tau:Q_{d(\mathcal I)}\to\bigvee_{\mathbf\alpha\in\mathcal I}Q_{n_{\mathbf\alpha}}$ of degree $q$. In degree zero this is a functor, and in positive degree it is a higher homotopy. For $q>0$, it satisfies
    \begin{equation}
        \mu_1^{\text{Fun}}\Xi^\tau
        =\Xi^{\mu\tau}+C_\Xi^\tau,
        \qquad
        \Xi^{\mu\tau}:=\sum_i\Xi^{\tau_i},
    \end{equation}
    where $\mu\tau=\sum_i\tau_i$.
\end{prop}
The term $C_\Xi^\tau$ contains the corresponding composition corrections for this family.
\begin{proof}
Again, the proof is similar to the parallel and series cases. We induct on the degree of $\tau$, starting with the degree-zero gluing functors. For positive degree, we construct the components in increasing input length using~\eqref{eq:xirecursion}, with
\begin{equation}
    B^\tau_k:=\sum_i\Xi^{\tau_i}_k,
    \qquad
    \widehat B^\tau_k:=B^\tau_k+(C_\Xi^\tau)_k.
\end{equation}
We use $\widehat B^\tau_k$ as the completed boundary in that recursion.
The right-hand side is a cycle by induction. Since $\Xi^\tau_k$ has degree $q+k-1$ and $\mu_1$ lowers degree by one, its degree on $(x_1,\dots,x_k)$ is
\[
    \sum_{j=1}^k|x_j|+k+q-2.
\]
Using Proposition~\ref{treehomcom}, we argue as in the parallel and series cases that this cycle is a boundary, treating the cases where the degrees agree as before. We can then choose a primitive and recursively define $\Xi^\tau_k$.
\end{proof}

Composing the compatible $\Psi$ and $\Xi$ families, with the same higher composition terms as in the parallel case, gives the family $\Theta^\tau$.
We can then define
\begin{equation}
    \overline{\text{Glue}_{\mathbf\alpha}}(\ff,\gg)
    :=\Theta^{\text{Glue}_{\mathbf\alpha}(\tau,\tau')},
\end{equation}
where $\Theta^\tau=\ff$ and $\Theta^{\tau'}=\gg$, for the chosen tree presentations. This gives the family required by Theorem~\ref{thm:all-gluings}.
\begin{rmk}
In several of the arguments above, we assumed that our functors $\ff:Q_n\to\cc$ had $n\geq2$. This was mostly to avoid splitting the arguments into cases. The same construction works when $n=1$. For the special case of $n=0$, note that gluing $\gg:Q_0\to\cc$ to $\ff:Q_n\to\cc$, with $n\geq1$, along $\gg_0(0)=\ff_0(i)$ requires $\ff_0(i)$ to be contractible.

The resulting operation gives a new functor
\begin{equation*}
    \gg\#\ff:Q_{n-1}\to\cc
\end{equation*}
in which we `miss out' the object $\ff_0(i)$. The construction follows the same idea as the general case, using the chosen contraction at this object.
The higher homotopy constructions follow the same argument. The only
additional terms come from $\mu_1(a_0)=1_0$; these are handled by strict
unitality, in the same way as the unit terms arising from full loops when $n>0$.

Finally, consider the special case of gluing $\ff:Q_0\to\cc$ to $\gg:Q_0\to\cc$ along their common object. The resulting gluing would formally give a functor $\ff\#\gg:Q_{-1}\to\cc$. We can either exclude this gluing or adopt the convention that $Q_{-1}$ is the empty category. 
\end{rmk}
\section{The multicategory}\label{sec:multicategory}
In this section, we consider some of the properties of the resulting structure built in the previous sections. Before we begin, we make the following remark.
\begin{rmk}
    The setup of the previous sections was over $k=\Z/2\Z$. This was really just to avoid too much notation. One should be able to extend the construction to work over $\Z$, with the signs calculated in the usual way. In general, this does not appear too difficult, other than perhaps in the proof of Lemma~\ref{lemphihomology}, where we used the assumption that $k$ is a field.
\end{rmk}
\subsection{$\text{Cob}(\cc)$}\label{subsec:cob-category}
We can consider the collection of all possible functors
\begin{equation}
    \ff:Q_n\to \cc \qquad n\geq 0
\end{equation}
for a chosen strictly unital $A_\infty$-category $\cc$. As in the introduction, we assume that $\cc$ has shifts and use strictly unital functors throughout. We define the following $A_\infty$-multicategory, denoted $\text{Cob}(\cc)$, with objects given by
\begin{equation}
    \text{Obj}(\text{Cob}(\cc)):=\text{Obj}(\cc).
\end{equation}
Given a collection of objects $A_0,\dots,A_{n-1},A_n\in \cc$, we can define the collection of multimorphisms
\begin{equation}
    \text{Cob}(A_0,\dots,A_{n-1}:A_n)
    :=\Big\{\ff:Q_n\to\cc\ \big|\ \ff_0(i)=A_i[s_i]\Big\},
\end{equation}
where, for $n\geq2$, we shift by
\begin{equation}
    s_i=\begin{cases}
        -1  &i=0\\
        0 & i=1 \text{ or }i=n\\
        1 & 2\leq i \leq n-1
    \end{cases}
\end{equation}
For $n=1$, take $s_0=s_1=0$, and for $n=0$ take $s_0=0$.
We can define a family of operations
\begin{equation}
\begin{aligned}
    \circ^p:\;&
    \left(\prod_{i\in\mathbf J}
    \text{Cob}(B^i_0,\dots,B^i_{m_i-1}:B^i_{m_i})\right)
    \times\text{Cob}(A_0,\dots,A_{n-1}:A_n)\\
    &\longrightarrow
    \text{Cob}(A_0,\dots,B^{i_j}_0,\dots,B^{i_j}_{m_{i_j}-1},
    \dots,A_{n-1}:A_n),
\end{aligned}
\end{equation}
with $\mathbf J\subseteq \{0,\dots,n-1\}$,  $B^{i}_{m_i}=A_i$ and $p$ denoting an order of composing. These operations are given by shifting the input functors appropriately, then applying the gluing dictated by the order $p$. The structure of the compositions of such operations is governed by Theorem~\ref{thm:all-gluings} and thus gives the structure of an $A_\infty$-multicategory.

\begin{rmk}
    The shifts in the definition are made to ensure that our multimorphisms respect the grading of $Q_n$.
\end{rmk}

\begin{defn}
    We define the homotopy multicategory $\text{H}_0\text{Cob}(\cc)$ to be the multicategory with the same objects, but multimorphisms given by homotopy classes of the corresponding $A_\infty$-functors. That is, given $\ff:Q_n\to \cc$ and $\ff':Q_n\to \cc$ with $\ff_0(i)=\ff'_0(i)$ they are equivalent as multimorphisms if they are homotopic as $A_\infty$-functors.

\end{defn}

One can define the notion of an algebra over a multicategory in a similar manner to that of an algebra over an operad, for example, see {\cite{Le}}. We give a quick recap in terms of $\text{H}_0\text{Cob}(\cc)$. An algebra over $\text{H}_0\text{Cob}(\cc)$ consists of a monoidal category $\dd$, an assignment on objects $\aa:\text{Obj}(\cc) \to \text{Obj}(\dd)$ and a family of maps, also denoted by $\aa$,
\begin{equation}
    \aa:\text{H}_0\text{Cob}(\cc)\left(A_0,\dots,A_{n-1}:A_n\right)\to \dd\left(\aa(A_0)\otimes\dots\otimes\aa(A_{n-1}),\aa(A_n)\right)
\end{equation}
which commute with the composition operations of $\text{H}_0\text{Cob}(\cc)$ and preserve identities. For a nullary multimorphism, the source is the monoidal unit of $\dd$. We denote such an algebra by $\aa:\text{H}_0\text{Cob}(\cc)\to\dd$.
\begin{rmk}
Heuristically, one should think of such an algebra as an invariant associated to the objects of $\cc$, which behaves compatibly with generalised cones.
\end{rmk}
\begin{theorem}\label{maintheoremalg}
    There is an algebra $\aa:\text{H}_0\text{Cob}(\cc)\to \dd$ which realises the Grothendieck group of $\text{H}_0\text{Tw}\cc$.
\end{theorem}
The proof follows from the next definition and subsequent results.
\begin{defn}
    Given such an $A_\infty$-category $\cc$, we define a group
    \begin{equation}
        K^Q_0(\cc):=\frac{\Z\cdot \text{Obj}(\cc)}{\sim}
    \end{equation}
    where $\sim$ denotes the collection of relations
    \begin{equation}
        \chi(\ff):=[\ff_0(0)]+\sum_{i=2}^{n}[\ff_0(i)]-[\ff_0(1)] \qquad [A[s]]=(-1)^s[A]
    \end{equation}
    for any $\ff:Q_n\to\cc$, with $n\geq0$. As in the introduction, we take $\chi(\ff)=[\ff_0(0)]-[\ff_0(1)]$ for $n=1$ and $\chi(\ff)=[\ff_0(0)]$ for $n=0$.
\end{defn}

\begin{prop}
    There is an algebra $\aa:\text{H}_0\text{Cob}(\cc)\to \dd$ that realises $K^Q_0(\cc)$.
\end{prop}
\begin{proof}
    We consider the category $\dd$ with objects
    \begin{equation}
        \text{Obj}(\dd):=K^Q_0(\cc)
    \end{equation}
    given by the elements of $K^Q_0(\cc)$. The monoidal structure is given by
    \begin{equation}
        X\otimes Y:=X+Y\in K^Q_0(\cc).
    \end{equation}
    We define shifts by setting $X[s]:=(-1)^sX$. The morphisms are given by
    \begin{equation}
        \dd(X,Y):=\begin{cases}
            \{1_X\} & X=Y\\\
            \emptyset & X\neq Y.
        \end{cases}
    \end{equation}

The assignment $\aa$ is defined on objects by
\begin{equation*}
    \aa(A):=[ A]
\end{equation*}
On multimorphisms, take $[\ff]\in\text{H}_0\text{Cob}(\cc)(A_0,\dots,A_{n-1}:A_n)$, with $\ff_0(i)=A_i[s_i]$. We set
\begin{equation}
    \aa([\ff]):\aa(A_0)\otimes \dots \otimes \aa(A_{n-1})\to \aa(A_n)
\end{equation}
to be the unique morphism. Indeed, for $n\geq2$ this is well-defined, since
\begin{equation}
    \aa(A_0)\otimes \dots \otimes \aa(A_{n-1})=\sum^{n-1}_{i=0} [\ff_0(i)[-s_i]]
\end{equation}
which gives
\begin{equation*}
    \sum^{n-1}_{i=0} [\ff_0(i)[-s_i]]=\sum_{i=0}^{n-1}(-1)^{s_i}[\ff_0(i)]=-[\ff_0(0)]+[\ff_0(1)]-\sum_{i=2}^{n-1}[\ff_0(i)].
\end{equation*}
But by the definition of $K^Q_0(\cc)$ we know
\begin{equation}
    -[\ff_0(0)]+[\ff_0(1)]-\sum_{i=2}^{n-1}[\ff_0(i)]=[\ff_0(n)]=\aa(A_n).
\end{equation}
For $n=1$, the relation is $[A_0]=[A_1]$, while for $n=0$ it is $[A_0]=0$.
\end{proof}

\begin{prop}
    If $A\oplus B$ exists in $\cc$, then $[A\oplus B]=[A]+[B]\in K^Q_0(\cc)$.
\end{prop}
\begin{proof}
    We build a functor $\ff:Q_2\to \cc$ with
    \begin{equation}
        \ff_0(0)= A \quad \ff_0(1)=A\oplus B \quad \ff_0(2)=B.
    \end{equation}
    We simply set
    \begin{align*}
        \ff_1(a_1)&=i:A\hookrightarrow A\oplus B\\
        \ff_1(a_2)&=p:A\oplus B \to B\\
        \ff_1(a_0)&=0\\
        \ff_2(a_1,a_2)&=0\\
        \ff_2(a_2,a_0)&=p':A\oplus B \to A\\
        \ff_2(a_0,a_1)&=i':B\hookrightarrow A\oplus B\\
        \ff_k(\dots)&=0 \text{ for }k>2.
    \end{align*}
    One can check this gives a well-defined $A_\infty$-functor, and thus imposes the relation
    \begin{equation}
        \chi(\ff)=0 =[A]-[A\oplus B]+[B].
    \end{equation}
\end{proof}
\begin{lemma}\label{lem:triangle-decomposition}
    In $\text{Tw}\cc$, up to homotopy, any functor $\ff:Q_n\to\text{Tw}\cc$ with $n\geq2$ decomposes as a sequence of gluings of functors $\gg^i:Q_2\to \text{Tw}\cc$, for $i=1,\dots,n-1$.
\end{lemma}
\begin{proof}
To begin, consider $\ff:Q_3\to\text{Tw}\cc$. We complete $\ff_1(a_1):\ff_0(0)\to \ff_0(1)$ to an exact triangle, i.e., there exists a functor $\gg:Q_2\to \text{Tw}\cc$ with $\gg_1(a_1)=\ff_1(a_1)$. Then consider the diagram
\begin{equation}
    \begin{tikzcd}
\ff_0(0) \arrow[rr, "\ff_1(a_1)"] &                                   & \ff_0(1) \arrow[dd, "\ff_1(a_2)"] \arrow[ld, "\gg_1(a_2)"] \\
                                  & \gg_0(2) \arrow[lu, "\gg_1(a_0)"] &                                                            \\
\ff_0(3) \arrow[uu, "\ff_1(a_0)"] &                                   & \ff_0(2) \arrow[ll, "\ff_1(a_3)"]
\end{tikzcd}
\end{equation}
showing the (partial) images of $\ff$ and $\gg$. We may take $\gg_0(2):=\left(\ff_0(1)\oplus \ff_0(0)[1],\ff_1(a_1)\right)$. We then consider the maps
\begin{equation}
u=
\begin{pmatrix}
\ff_2(a_0,a_1)\\
\ff_1(a_0)
\end{pmatrix}:\ff_0(3)\to \gg_0(2),
\qquad
v=
\begin{pmatrix}
\ff_1(a_2)&\ff_2(a_1,a_2)
\end{pmatrix}:\gg_0(2)\to \ff_0(2).
\end{equation}
\begin{equation}
    \begin{tikzcd}
\ff_0(0) \arrow[rr, "\ff_1(a_1)"] &                                   & \ff_0(1) \arrow[dd, "\ff_1(a_2)"] \arrow[ld, "\gg_1(a_2)"] \\
                                  & \gg_0(2) \arrow[lu, "\gg_1(a_0)"] \ar[dr,dashed,"v"]&                                                            \\
\ff_0(3) \arrow[uu, "\ff_1(a_0)"]\ar[ur,dashed,"u"] &                                   & \ff_0(2) \arrow[ll, "\ff_1(a_3)"]
\end{tikzcd}
\end{equation}
One finds that there exists a $\hh:Q_2\to \text{Tw}\cc$ with $\hh_0(0)=\ff_0(3),\, \hh_0(1)=\gg_0(2)$, and $\hh_0(2)=\ff_0(2)$ such that
\begin{align*}
    \hh_1(a_1)=&u\\
    \hh_1(a_2)=&v\\
    \hh_1(a_0)=&\ff_1(a_3).
\end{align*}
Indeed, we have
\[
\begin{aligned}
\mu_2\bigl(v,\ff_1(a_3)\bigr)
={}&
\begin{pmatrix}
\mu_2\bigl(\ff_1(a_2),\ff_1(a_3)\bigr)
&
\mu_2\bigl(\ff_2(a_1,a_2),\ff_1(a_3)\bigr)
\end{pmatrix}\\
&+
\begin{pmatrix}
0&
\mu_3\bigl(\ff_1(a_1),\ff_1(a_2),\ff_1(a_3)\bigr)
\end{pmatrix}.
\end{aligned}
\]
which we can use to realise
\[
\mu_2\bigl(v,\ff_1(a_3)\bigr)
=
\mu_1
\begin{pmatrix}
\ff_2(a_2,a_3)&\ff_3(a_1,a_2,a_3)
\end{pmatrix}.
\]
Thus we may take \[
\hh_2(a_2,a_0)=
\begin{pmatrix}
\ff_2(a_2,a_3)&\ff_3(a_1,a_2,a_3)
\end{pmatrix}.
\]
The other two order-two terms are given by
\[
\begin{aligned}
\hh_2(a_0,a_1)
&=
\begin{pmatrix}
\ff_3(a_3,a_0,a_1)\\
\ff_2(a_3,a_0)
\end{pmatrix},\\
\hh_2(a_1,a_2)
&=\ff_3(a_0,a_1,a_2).
\end{aligned}
\]
One can go on to define every $\hh_k$ in terms of the $\ff_l$'s, but for the purpose of the statement it is enough to know they exist. We can then consider $\gg\#_{21}\hh$. We find
\begin{align*}
    (\gg\#_{21}\hh)_0(0)&=\ff_0(0)\\
    (\gg\#_{21}\hh)_0(1)&=\ff_0(1)\\
    (\gg\#_{21}\hh)_0(2)&=\ff_0(2)\\
    (\gg\#_{21}\hh)_0(3)&=\ff_0(3)\\
    (\gg\#_{21}\hh)_1(a_1)&=\ff_1(a_1)\\
    (\gg\#_{21}\hh)_1(a_2)&=\mu_2(\gg_1(a_2),v)=\ff_1(a_2)\\
    (\gg\#_{21}\hh)_1(a_3)&=\ff_1(a_3)\\
    (\gg\#_{21}\hh)_1(a_0)&=\mu_2(u,\gg_1(a_0))= \ff_1(a_0)
\end{align*}
where $\gg_1(a_0)$ was the projection $(0 \quad 1_{\ff_0(0)[1]})$ and $\gg_1(a_2)$ was the inclusion $\begin{pmatrix}1_{\ff_0(1)}\\0
\end{pmatrix}$. One can go on to explicitly construct the appropriate higher order terms. Thus, at the homotopy level $\ff\simeq \gg\#\hh$. For $n>3$, the proof is done by iterating this procedure. One finds that $\ff:Q_n\to \text{Tw}\cc$ can be written as $\ff\simeq \gg^1\#\hh^1$, with $\gg^1:Q_2\to \text{Tw}\cc$ and $\hh^1:Q_{n-1}\to \text{Tw}\cc$, then repeat with $\hh^1\simeq\gg^2\#\hh^2$, and so on.
\end{proof}
\begin{prop}\label{chisplit}
    Given two functors $\ff:Q_2\to \cc$ and $\gg:Q_2\to \cc$ with $\ff_0(2)=\gg_0(j)$ (with $j\neq 2$), the resulting gluing satisfies
    \begin{equation}
        \chi(\ff\#_{2j}\gg)=\chi(\ff)+(-1)^{j-1}\chi(\gg)
    \end{equation}
\end{prop}
\begin{proof}
    This is a straightforward calculation; for $j=1$ we have
    \begin{align*}
        \chi(\ff\#_{2j}\gg)=&[(\ff\#_{2j}\gg)_0(0)]+[(\ff\#_{2j}\gg)_0(2)]+[(\ff\#_{2j}\gg)_0(3)]-[(\ff\#_{2j}\gg)_0(1)]\\
        =&[\ff_0(0)]+[\gg_0(2)]+[\gg_0(0)]-[\ff_0(1)]\\
        =&[\ff_0(0)]+[\ff_0(2)]-[\ff_0(1)]+[\gg_0(0)]-[\gg_0(1)]+[\gg_0(2)]\\
        =& \chi(\ff)+\chi(\gg).
    \end{align*}
    For $j=0$ we have
    \begin{align*}
    \chi(\ff\#_{20}\gg)
    &=[\ff_0(0)]-[\ff_0(1)]
      +[\gg_0(1)]+[\gg_0(2)[-1]]\\
    &=[\ff_0(0)]-[\ff_0(1)]
      +[\gg_0(1)]-[\gg_0(2)]\\
    &=\bigl([\ff_0(0)]-[\ff_0(1)]+[\ff_0(2)]\bigr)\\
    &\quad-\bigl([\gg_0(0)]-[\gg_0(1)]+[\gg_0(2)]\bigr)\\
    &=\chi(\ff)-\chi(\gg),
\end{align*}
where we used $[A[-1]]=-[A]$.
\end{proof}
\begin{cor}
    Given $\ff:Q_n\to \cc$, if it splits as a sequence of gluings of $\gg^i:Q_2\to \cc$ for $i=1,\dots,n-1$ then
    \begin{equation}
        \chi(\ff)=\sum_{i=1}^{n-1}(-1)^{r_i}\chi(\gg^i)
    \end{equation}
    where $r_i$ depends on the glued objects.
\end{cor}
\begin{proof}
    Iterate Proposition~\ref{chisplit}.
\end{proof}

\begin{rmk}
    We can choose to split our $\ff:Q_n\to \cc$ by a sequence of gluings of $\gg^i$ so that every gluing in the sequence has the form $\gg^i\#_{21}\hh^i$. Thus we can choose the family $\gg^i$ so that
     \begin{equation}
        \chi(\ff)=\sum_{i=1}^{n-1}\chi(\gg^i)
    \end{equation}
\end{rmk}
\begin{lemma}\label{lem:q-k0}
    There is an isomorphism of groups \begin{equation}
        \lambda:K^Q_0(\cc) \to K_0(\text{H}_0\text{Tw}\cc)
    \end{equation}
    with $\lambda([A])=[A]$.
\end{lemma}
\begin{proof}
First, we check that $\lambda$ is a well-defined group homomorphism. We have
\begin{equation*}
  \lambda([A]+[B])=[A]+[B]=\lambda([A])+\lambda([B]).
\end{equation*}
Moreover, any relation imposed by $K^Q_0$ gives a relation in $K_0(\text{H}_0\text{Tw}\cc)$. Indeed, it is enough to check that for any $\ff:Q_n\to \cc$ (i.e., any generator of the relations), the relation given by $\chi(\ff)=0$ is imposed in $K_0(\text{H}_0\text{Tw}\cc)$. This is straightforward. For $n=0$, the image object is contractible, and for $n=1$, the two image objects are homotopy equivalent. For $n\geq2$, consider the embedding $\iota:\cc\hookrightarrow\text{Tw}\cc$ and $\iota\circ \ff:Q_n\to \text{Tw}\cc$. By Lemma~\ref{lem:triangle-decomposition} it follows that $\iota\circ \ff$ splits as a sequence of gluings of $\gg^i:Q_2\to\text{Tw}\cc$. Thus, we impose the relations
\begin{equation}
[\gg^i_0(0)]-[\gg^i_0(1)]+[\gg^i_0(2)]=0
\end{equation}
in $K_0(\text{H}_0\text{Tw}\cc)$. But we realise from the above remark that we can choose our $\gg^i$ such that
\begin{equation}
    \chi(\iota\circ \ff)=\sum_{i=1}^{n-1}\chi(\gg^i)=\sum_{i=1}^{n-1}\left([\gg^i_0(0)]-[\gg^i_0(1)]+[\gg^i_0(2)]\right)
\end{equation}
so every relation in $K^Q_0(\cc)$ is also a relation in $K_0(\text{H}_0\text{Tw}\cc)$.

The map is surjective, since every twisted complex is built from finitely many shifts of objects of $\cc$. The main difficulty lies in proving injectivity. Thus we need to show that every relation generated by exact triangles in $\text{H}_0\text{Tw}\cc$ can be written in the form $\sum a_i\chi(\iota\circ\ff^i)$ for some $\ff^i:Q_{n_i}\to \cc$. We only need to check relations between elements of the form $[\iota(A)]$, as every twisted complex $X\in \text{Tw}\cc$ will split as $[X]=\sum a_i[\iota(A_i)]$. Consider a relation
\begin{equation}
    R:=\sum_{i=0}^ka_i[A_i]
\end{equation}
with $A_i\in \cc$. We can rewrite $R$ as
\begin{equation}
    R=\sum_{i=1}^lc_i \left([B^i_0]-[B^i_1]+[B^i_2]\right)
\end{equation}
since exact triangles generate the relations. Each $B_j^i$ is itself a twisted complex. For a generic $B\in \text{Tw}\cc$, we can write
\begin{equation*}
B=\left(\bigoplus_{r=1}^{l}A_r[s_r],\delta\right),
\end{equation*}
and set
\begin{equation*}
    e(B):=\sum_{r=1}^{l}(-1)^{s_r}[A_r].
\end{equation*}
We choose $e(\iota(A))=[A]$ and extend $e$ linearly to formal sums of twisted-complex objects.
This twisted complex is an iterated mapping cone of the objects $A_r[s_r]$, so there is a functor $\hh^B:Q_l\to \text{Tw}\cc$ with $\hh^B_0(1)=B$ and every $\hh^B_0(i)$ for $i\neq1$ a copy of some $A_r[s_r]$. Thus we have
\begin{equation}
    \chi(\hh^B)=e(B)-[B].
\end{equation}
For every $B^i_j$ we choose such a decomposition. Then we consider functors $\gg^i:Q_2\to \text{Tw}\cc$ which represent each triangle in the decomposition of $R$, i.e., $\chi(\gg^i)=[B^i_0]-[B^i_1]+[B^i_2]$. We can then build a functor $\ff^i:Q_{n_i}\to \cc$ which is the result of gluing each $\hh^{B^i_j}$ to the copy of $B^i_j$ in $\gg^i$. It follows that
\begin{align*}
    \chi(\ff^i)=&\chi(\gg^i)+\chi(\hh^{B^i_0})-\chi(\hh^{B^i_1})+\chi(\hh^{B^i_2})\\
    =&[B^i_0]-[B^i_1]+[B^i_2]+e(B^{i}_0)-[B^{i}_0]\\
    &\quad-(e(B^{i}_1)-[B^{i}_1])+e(B^{i}_2)-[B^{i}_2]\\
    =&e(B^i_0)-e(B^i_1)+e(B^i_2).
\end{align*}
Thus we have
\begin{equation}
    \begin{aligned}
    R&=\sum_{i=1}^l c_i\left([B^i_0]-[B^i_1]+[B^i_2]\right)\\
     &=\sum_{i=1}^l c_i\left(e(B^i_0)-e(B^i_1)+e(B^i_2)\right)
      =\sum_{i=1}^l c_i\chi(\ff^i).
\end{aligned}
\end{equation}
\end{proof}

\begin{rmk}
    The above lemma, along with the realisation of $K^Q_0(\cc)$ as an algebra over $\text{H}_0\text{Cob}(\cc)$, supplies the proof of Theorem~\ref{maintheoremalg}.
\end{rmk}

\subsection{Generators}
As we proved before, if $\cc$ is pre-triangulated then we can split any $\ff:Q_n\to \cc$ with $n\geq2$ into a sequence of gluings of $\gg:Q_2\to \cc$. However, if $\cc$ is not pre-triangulated this can fail in general. Thus, we should view taking twisted complexes as forcing these splittings. It is interesting, though, to consider what happens in the case of non-pre-triangulated categories. We have already seen that $K^Q_0(\cc)=K_0(\text{H}_0\text{Tw}\cc)$, so we also find that $K^Q_0(\cc)=K^Q_0(\text{Tw}\cc)$ since $\text{TwTw}\cc\simeq \text{Tw}\cc$. Thus we cannot distinguish $\cc$ from its pre-triangulated completion in terms of $K^Q$, i.e., the Grothendieck group.

The two categories can still be different. There can exist a functor $\ff:Q_n\to \cc$ which does not decompose as a sequence of gluings of $\gg^i:Q_2\to \cc$. Thus, we want to understand the following question: which functors $\ff:Q_n\to \cc$ generate every multimorphism in $\text{Cob}(\cc)$ or perhaps every multimorphism of $\text{H}_0\text{Cob}(\cc)$?

\begin{defn}
    We say $\cc$ is \textbf{$n$-constructable} if every multimorphism in $\text{Cob}(\cc)$ can be recovered as a finite sequence of gluings of $(\gg^i:Q_{k_i}\to \cc)_{i=0,\dots,l}$, with $k_i\leq n$ for all $i=0,\dots, l$. Similarly, $\cc$ is weakly $n$-constructable if the same is true for every multimorphism in $\text{H}_0\text{Cob}(\cc)$.
\end{defn}

    If $\cc$ is $n$-constructable, then clearly it is also $m$-constructable for every $m>n$. Thus we consider
    \begin{equation}
        \text{con}(\cc):=\min\left\{n\geq0:\cc\text{ is $n$-constructable}\right\},
    \end{equation}
    with $\text{con}(\cc)=\infty$ if there is no finite bound.

\begin{lemma}\label{lem:arity-inclusion}
    Given a strictly unital $A_\infty$-category $\cc$ with a
    chosen strict zero object, for any $n\geq0$ there exists
    a fully faithful functor
    \[
        \ii_{n,n+1}:
        \text{Fun}_{A_\infty}(Q_n,\cc)
        \hookrightarrow
        \text{Fun}_{A_\infty}(Q_{n+1},\cc).
    \]
\end{lemma}

\begin{proof}
There is in fact a family of such functors. We choose
$\ii_{n,n+1}$ as follows. First consider the category $Q_n^+$,
obtained by adjoining a strict zero object $\bar{0}$ to $Q_n$.
Extend every $\ff:Q_n\to\cc$ to
$\ff^+:Q_n^+\to\cc$ by setting
\[
    \ff^+_0(\bar{0})=0_\cc
\]
and retaining the components of $\ff$ on $Q_n$. We then define a functor
\[
    \bar{\ii}_{n+1,n}:Q_{n+1}\longrightarrow Q_n^+.
\]
We use $a_i$ for the arrows in both categories. Set
\[
(\bar{\ii}_{n+1,n})_0(i)=
\begin{cases}
    i,&0\leq i\leq n,\\
    \bar{0},&i=n+1,
\end{cases}
\qquad
(\bar{\ii}_{n+1,n})_1(a_i)=
\begin{cases}
    a_i,&1\leq i\leq n,\\
    0,&i=0\text{ or }i=n+1,
\end{cases}
\]
and
\[
    (\bar{\ii}_{n+1,n})_2(a_{n+1},a_0):=a_0.
\]
All other components on sequences of non-identity arrows
are zero. For sequences with units we take
\[
    (\bar{\ii}_{n+1,n})_1(1_i)
       =1_{(\bar{\ii}_{n+1,n})_0(i)},
    \qquad
    (\bar{\ii}_{n+1,n})_k(\ldots,1_i,\ldots)=0\quad(k\geq2).
\]
In particular,
$(\bar{\ii}_{n+1,n})_1(1_{n+1})=1_{\bar{0}}=0$. One should picture, for example,
\begin{equation}
\begin{tikzcd}[row sep=3em,column sep=3.5em]
    0 \arrow[r,"a_1"]
      & 1 \arrow[d,"a_2"]\\
    3 \arrow[u,dashed,"a_0"]
      & 2 \arrow[l,dashed,"a_3"]
\end{tikzcd}
\qquad\xrightarrow{\ \bar{\ii}_{3,2}\ }\qquad
\begin{tikzcd}[row sep=3em,column sep=3.5em]
    0 \arrow[r,"a_1"]
      & 1 \arrow[d,"a_2"]\\
    \bar{0}
      & 2 \arrow[ul,"a_0" description]
\end{tikzcd}
\end{equation}
where $(\bar{\ii}_{3,2})_2(a_3,a_0)=a_0$. We claim that $\bar{\ii}_{n+1,n}$ is well-defined. For $1\leq i\leq n$, we have
\[
    \mu_1(\bar{\ii}_{n+1,n})_1(a_i)
    =\mu_1(a_i)
    =0
    =(\bar{\ii}_{n+1,n})_1(\mu_1(a_i)).
\]
For the distinguished pair, we have
\begin{align*}
&\mu_1(\bar{\ii}_{n+1,n})_2(a_{n+1},a_0)
 +\mu_2\Bigl(
    (\bar{\ii}_{n+1,n})_1(a_{n+1}),
    (\bar{\ii}_{n+1,n})_1(a_0)
  \Bigr)\\
&\quad=\mu_1(a_0)\\
&\quad=
\begin{cases}
    1_0,&n=0,\\
    0,&n\geq1,
\end{cases}\\
&\quad=(\bar{\ii}_{n+1,n})_1
          \bigl(\mu_2(a_{n+1},a_0)\bigr).
\end{align*}

More generally, we must consider
\[
\begin{aligned}
&\sum_{k=1}^{\ell}
 \sum_{\substack{p_1+\cdots+p_k=\ell\\p_h\geq1}}
 \mu_k\Bigl(
   (\bar{\ii}_{n+1,n})_{p_1}
      (a_{i+1},\ldots,a_{i+p_1}),\ldots,\\
&\qquad
   (\bar{\ii}_{n+1,n})_{p_k}
      (a_{i+\ell-p_k+1},\ldots,a_{i+\ell})
 \Bigr).
\end{aligned}
\]
If the word starts or ends at the new object $n+1$, both sides of the functor equation vanish, since its image is the strict zero object. We therefore consider words whose endpoints are retained. Almost all terms vanish: each surviving block is either a single $a_r$ with $1\leq r\leq n$, or the pair $(a_{n+1},a_0)$. We compress each such pair to $a_0$. Indices before compression are modulo $n+2$, and those after compression are modulo $n+1$. If $t$ is the number of compressed pairs, then
\[
\begin{aligned}
&\sum_{k=1}^{\ell}
 \sum_{\substack{p_1+\cdots+p_k=\ell\\p_h\geq1}}
 \mu_k\Bigl(
   (\bar{\ii}_{n+1,n})_{p_1}
      (a_{i+1},\ldots,a_{i+p_1}),\ldots,\\
&\qquad
   (\bar{\ii}_{n+1,n})_{p_k}
      (a_{i+\ell-p_k+1},\ldots,a_{i+\ell})
 \Bigr)\\
&\quad=
 \mu_{\ell-t}(a_{i+1},\ldots,a_{i+\ell-t}),
\end{aligned}
\]
If $t>1$, the output must contain $a_0$
more than once and so
\[
    \mu_{\ell-t}(a_{i+1},\ldots,a_{i+\ell-t})=0.
\]
Similarly, if $t=0$, the output also vanishes, agreeing with the source terms. It remains to consider $t=1$. The output again
vanishes unless $\ell=n+2$. In this case we have
\begin{align*}
    \mu_{n+1}(a_{i+1},\ldots,a_{i})
    &=1_i\\
    &=(\bar{\ii}_{n+1,n})_1(1_i)\\
    &=(\bar{\ii}_{n+1,n})_1
       \bigl(\mu_{n+2}(a_{i+1},\ldots,a_{i})\bigr).
\end{align*}

We can then define
\[
    \ii_{n,n+1}(\ff)
    :=(\bar{\ii}_{n+1,n})^*(\ff^+).
\]
Now we check that this is fully faithful. For faithfulness, consider $\ff,\gg:Q_n\to \cc$. We have
\begin{equation}
    \text{Fun}_{A_\infty}(Q_n,\cc)\left(\ff,\gg\right)=\text{Nat}\left (\ff,\gg\right)
\end{equation}
Such an $f\in\text{Nat}(\ff,\gg)$ consists of a family
\begin{equation}
    f_k:Q_n(x_0,x_1)\otimes\dots\otimes Q_n(x_{k-1},x_k)
    \longrightarrow\cc(\ff_0(x_0),\gg_0(x_k)),
    \qquad k\geq0,
\end{equation}
where $x_0,\dots,x_k$ are objects. We use normalized natural transformations, so components with a unit input vanish. On non-identity arrows, we therefore only need to consider words $(a_{i+1},\dots,a_{i+k})$.

To show faithfulness, consider $\bar f:=\ii_{n,n+1}(f)=(\bar\ii_{n+1,n})^*f^+$. In arity zero,
\[
    \bar f_0(i)=f_0(i)\quad(0\leq i\leq n),
    \qquad \bar f_0(n+1)=0.
\]
For $\ell\geq1$, one finds
\[
\begin{aligned}
\bar f_\ell(x_1,\ldots,x_\ell)
={}&
\sum_{r=1}^{\ell}
\sum_{\substack{p_1+\cdots+p_r=\ell\\p_j\geq1}}
f^+_r\Bigl(
(\bar{\ii}_{n+1,n})_{p_1}(x_1,\ldots,x_{p_1}),\ldots,\\
&\qquad
(\bar{\ii}_{n+1,n})_{p_r}
(x_{\ell-p_r+1},\ldots,x_\ell)
\Bigr).
\end{aligned}
\]
Only one term can be nonzero on the RHS. Indeed, an input to $f^+_r$ will be zero unless each block contributes an $(\bar{\ii}_{n+1,n})_1(a_i)=a_i$ or a $(\bar{\ii}_{n+1,n})_2(a_{n+1},a_0)=a_0$. Thus we find the only non-zero terms are
\[
\bar f_\ell(a_{i+1},\ldots,a_{i+\ell})
=
f_{\ell-t}(a_{i+1},\ldots,a_{i+\ell-t}).
\]
Faithfulness follows by observing that every original arrow word has a unique expansion, obtained by replacing each $a_0$ with $(a_{n+1},a_0)$. Thus, if $\bar{f}=\bar{g}$, then for every $\ell$ and expanded arrow word,
\begin{align*}
 \bar{f}_\ell(a_{i+1},\ldots,a_{i+\ell})
=&
f_{\ell-t}(a_{i+1},\ldots,a_{i+\ell-t})\\
=&g_{\ell-t}(a_{i+1},\ldots,a_{i+\ell-t})\\=& \bar{g}_\ell(a_{i+1},\ldots,a_{i+\ell}),
\end{align*}
i.e., $g=f$.

For fullness, consider $g\in \text{Nat}(\ii_{n,n+1}(\ff),\ii_{n,n+1}(\gg))$. We can recover $g$ as $g=\bar{f}$ for $f\in \text{Nat}(\ff,\gg)$ by defining
\begin{equation}
   f_{\ell-t}(a_{i+1},\ldots,a_{i+\ell-t}):=g_\ell(a_{i+1},\ldots,a_{i+\ell})
\end{equation}
with $f_0(i)=g_0(i)$. This defines every component of $f$, since the expansion is unique. Components of $g$ with an endpoint at the new zero object vanish, so $\bar f=g$.
\end{proof}

\begin{cor}
    For every $m>n$ the functor
    \begin{equation}\ii_{n,m}:=\ii_{m-1,m}\circ\dots\circ\ii_{n,n+1}:\text{Fun}_{A_\infty}\left( Q_n,\cc\right)\hookrightarrow \text{Fun}_{A_\infty}\left( Q_{m},\cc\right)
    \end{equation}
    is fully faithful.
\end{cor}

As noted previously, the group $K^Q_0(\cc)$ cannot distinguish $\cc$ from its pre-triangulated closure. Thus, by itself, it does not distinguish these two categories. However, we can obtain a family of groups associated to $\cc$.

\begin{defn}
    Let $n\geq0$ and let $\cc$ be a strictly unital $A_\infty$-category with shifts. We define
    \begin{equation}
        K_0^{Q;n}(\cc):=\frac{\Z\cdot \text{Obj}(\cc)}{\sim_n}
    \end{equation}
    where $\sim_n$ denotes the collection of relations
    \begin{equation}
        \chi(\ff):=[\ff_0(0)]+\sum_{i=2}^{m}[\ff_0(i)]-[\ff_0(1)] \qquad [A[s]]=(-1)^s[A]
    \end{equation}
    for any $\ff:Q_m\to\cc$ with $0\leq m\leq n$, using the same conventions for $m=0,1$ as above.
\end{defn}

\begin{prop}
    There is a family of surjective group homomorphisms
    \begin{equation}
        \lambda_{n,n+1}:K_0^{Q;n}(\cc)\to K_0^{Q;n+1}(\cc)
    \end{equation}
    given by \begin{equation}
        \lambda_{n,n+1}([A])=[A].
    \end{equation}
\end{prop}
\begin{proof}
    Both groups are generated by the same objects, and every relation imposed in $K^{Q;n}_0(\cc)$ is also imposed in $K^{Q;n+1}_0(\cc)$. Thus the map is well-defined and surjective.
\end{proof}

\begin{prop}
  There is an isomorphism
    \begin{equation}
        K^Q_0(\cc) \cong \text{colim}_{n\geq 0}\left(\to K^{Q;n}_0(\cc) \xrightarrow{\lambda_{n,n+1}}K^{Q;n+1}_0(\cc)\to \dots \right)
    \end{equation}
\end{prop}

\begin{proof}
    Every relation defining $K^Q_0(\cc)$ occurs at a finite stage: the relation $\chi(\ff)$ for $\ff:Q_m\to\cc$ is imposed in $K^{Q;m}_0(\cc)$, and the shift relations occur at every stage. Thus the union of the relation subgroups is exactly the subgroup defining $K^Q_0(\cc)$.
\end{proof}
\begin{rmk}
    Taking the colimit gives a family of morphisms
    \begin{equation}
       \begin{tikzcd}[column sep=small]
       K^{Q;0}_0(\cc)\ar[r,"{\lambda_{0,1}}"]\ar[drr,"\lambda_0"] & K^{Q;1}_0(\cc) \ar[dr,"\lambda_1"]\ar[r,"{\lambda_{1,2}}"] &\dots\ar[r,"{\lambda_{n-1,n}}"] & K^{Q;n}_0(\cc)\ar[dl,"\lambda_n"] \ar[r,"{\lambda_{n,n+1}}"]& K^{Q;n+1}_0(\cc) \ar[dll,"\lambda_{n+1}"]\ar[r,"{\lambda_{n+1,n+2}}"]&\dots \\
    &&  K^Q_0(\cc)
       \end{tikzcd}
    \end{equation}
    defined in the obvious way.
\end{rmk}

From this sequence we can define an invariant of $\cc$ as follows.
\begin{defn}
    We define
    \begin{equation}
\kappa(\cc):=\inf_{n\geq0}\left\{n:\lambda_n \text{ is an isomorphism}\right\},
    \end{equation}
    with $\kappa(\cc)=\infty$ if no $\lambda_n:K^{Q;n}_0(\cc)\to K^Q_0(\cc)$ is an isomorphism.
\end{defn}
Equivalently, we have
\begin{equation}
    \kappa(\cc)=\min\left\{n\geq0:\ker(\lambda_n)=0\right\},
\end{equation}
with the same convention when the set is empty.

\begin{prop}
    The category $Q_n$ has $\kappa(Q_n)=n$.
\end{prop}
\begin{proof}
    For $n\geq1$, consider the identity functor $1_{Q_n}$, taking values in (the shift completion of $Q_n$). Its relation is
    \begin{equation}
        [1]=[0]+\sum_{i=2}^{n}[i].
    \end{equation}
 The free group on objects is $\Z^{n+1}$, with generators $[0],\dots,[n]$. Thus,
    \begin{equation}
        \frac{\Z^{n+1}}{\langle\chi(1_{Q_n})\rangle}\cong\Z^n.
    \end{equation}
    By \cite{Na}, we know that
    \begin{equation}
        K^Q_0(Q_n)\cong K_0(\text{H}_0\text{Tw}Q_n)\cong\Z^n.
    \end{equation}
    Thus, $\lambda_{n}$ must be an isomorphism and every relation is a multiple of $\chi(1_{Q_n})$, and we have $K^{Q;n}_0(Q_n)\cong\Z^n$.

    If $\ff:Q_l\to Q_n$ has $l<n$, then the relation imposed (i.e., $\chi(\ff)$) must be trivial, as else $\ff$ would represent a multiple of the relation imposed by $1_{Q_n}$ but this relation used $n+1$ terms and $l+1<n+1$. In particular
    \begin{equation}
        K^{Q;l}_0(Q_n)\cong\Z^{n+1}\qquad\text{for }0\leq l<n,
    \end{equation}
    which proves $\kappa(Q_n)=n$. For $n=0$, the relation $\chi(1_{Q_0})=[0]$ gives $K^{Q;0}_0(Q_0)=K^Q_0(Q_0)=0$, so $\kappa(Q_0)=0$ as well.
\end{proof}

This invariant thus distinguishes $Q_n$ from $\text{Tw}Q_n$ {(for $n>2$)}.

\begin{lemma}
    If $\cc$ is weakly $n$-constructable, then $\kappa(\cc)\leq n$.
\end{lemma}
\begin{proof}
  Indeed, the collection of all $\ff:Q_m\to \cc$ for any $m\geq 0$ generates the relations of $K^Q_0$, so any relation is of the form
  \begin{equation*}
      R=\sum_{i=1}^lc_i\chi(\ff^i).
  \end{equation*}
But under the assumptions, each $\ff^i$ can be decomposed as a sequence of gluings of $\gg^{ij}:Q_{n_{ij}}\to \cc$ with $n_{ij}\leq n$ and $1\leq j\leq k_i$. It follows that we can rewrite $R$ as
\begin{equation}
    R=\sum_{i=1}^l\sum_{j=1}^{k_i}d_{ij}\chi(\gg^{ij}).
\end{equation}
In particular, every relation defining $K^{Q}_0(\cc)$ already appears in the relations defining $K^{Q;n}_0(\cc)$, so $\lambda_n$ must be an isomorphism.
\end{proof}

\begin{cor}
    We have $\text{con}(\cc)\geq \kappa(\cc)$.
\end{cor}
\begin{proof}
    Clear.
\end{proof}

\begin{rmk}
    The converse would imply that if $\chi(\ff)=\sum_{i=1}^l c_i\chi(\gg^i)$ with $\gg^i:Q_{k_i}\to \cc$ such that $k_i\leq n=\kappa(\cc)$, then there exists a sequence of gluings of the $\gg^i$ decomposing $\ff$. This seems non-obvious.
\end{rmk}

More specifically, we can say that a collection of multimorphisms $\mathcal{S}$ of $\text{H}_0\text{Cob}(\cc)$ \textbf{generates} $\text{H}_0\text{Cob}(\cc)$ if every multimorphism factors as a sequence of compositions of elements of $\mathcal{S}$. It follows that
\begin{prop}
    If $\mathcal{S}$ generates $\text{H}_0\text{Cob}(\cc)$, then
    \begin{equation}
        K^Q_0(\cc)\cong \frac{\Z\cdot \text{Obj}(\cc)}{R_{\mathcal S}},
    \end{equation}
    where $R_\mathcal S$ denotes the subgroup generated by $\chi(\gg)$ with $\gg\in\mathcal S$, along with the relations imposed by shifts.
\end{prop}

\begin{proof}
    This is clear.
\end{proof}
\begin{rmk}
The converse would assert that every collection of generators of the relations defining $K^Q_0(\cc)$ would also generate $\text{H}_0\text{Cob}(\cc)$.
\end{rmk}
One can consider a stronger statement. We define a multimorphism $(\ff:Q_n\to \cc)\in \text{H}_0\text{Cob}(\cc)$ to be \textbf{indecomposable} if there does not exist a sequence $(\gg^i:Q_{k_i}\to \cc)_{i=1,\dots,l}$ whose gluings produce $\ff$ (up to homotopy) and such that $k_j<n$ for every $j$. We can then say a collection $\mathcal{S}$ of multimorphisms indecomposably generates $\text{H}_0\text{Cob}(\cc)$ if $\mathcal{S}$ generates it, and if every $\gg\in \mathcal{S}$ is indecomposable. Such a generating collection is, in some sense, optimal, as we cannot decompose any of the elements further.
\begin{cor}
    If $\mathcal{S}$ indecomposably generates $\text{H}_0\text{Cob}(\cc)$ then
    \begin{equation}
        K^Q_0(\cc)\cong \frac{\Z\cdot \text{Obj}(\cc)}{R_{\mathcal S}}.
    \end{equation}
\end{cor}

\section{Comments and questions}\label{sec:questions}
Here, we list a few remarks on the previous material, as well as comment on a few interesting questions that should be explored.

\subsection{Lagrangian cobordisms}
The idea for developing this theory originates from Lagrangian cobordisms. It is known, as mentioned in the introduction, that Lagrangian cobordisms induce iterated cone decompositions in the derived Fukaya category. This result is due in large part to Biran and Cornea (see {\cite{BC1,BC2,BC3}}), and was also more recently proved in the wrapped setting by Bosshard {\cite{Bo}}. In a later paper of Biran and Cornea entitled ‘A Lagrangian pictionary’ {\cite{BC4}}, they show that one can recover the entire derived Fukaya category geometrically, in terms of composing equivalence classes of Lagrangian cobordisms (working in the exact, marked, unobstructed immersed setting).

This paper thus began as a study of how to interpret this gluing algebraically. In the derived setting it is somewhat obvious: one concatenates two iterated cone decompositions to obtain a larger one. But at the chain level, and as purely geometric objects, it is a little more mysterious. Interpreting cobordisms as $A_\infty$-functors into $\text{Fuk}_{\text{emb}}(M)$ seems closer to the underlying geometry. Indeed, decomposing such a functor into gluings requires the intermediate objects to remain embedded Lagrangians, whereas an iterated cone decomposition in the derived category may involve intermediate objects for which no embedded representative is available.

More generally, when we study geometric objects through an associated $A_\infty$- or dg category, passing to twisted complexes forces us to consider objects that might not have an immediate geometric description.

Moreover, it is not clear in general whether such cobordisms can always be decomposed into cobordisms with fewer ends while keeping the intermediate ends embedded. This again is much more explicitly seen by working with $A_\infty$-functors and by considering the multicategory $\text{H}_0\text{Cob}(\text{Fuk}_{\text{emb}}(M))$.

One would expect, and it is somewhat easy to realise, a corresponding multicategory of cobordisms, that is, a multicategory $\text{H}_0\text{Cob}(M)$ with objects given by Lagrangians (with whatever flavour one prefers, perhaps embedded for the sake of this argument) and multimorphisms given by Lagrangian cobordisms (again with whatever conditions one would like). One can then interpret the surjective homomorphism
\begin{equation}
   \theta:\Omega_{\text{Lag}}(M)\to K_0\left(D\text{Fuk}(M)\right)
\end{equation}
(see, for example, \cite{Ha}; note that one requires additional Floer-theoretic assumptions for this to be well defined)
in terms of a morphism of multicategories
\begin{equation}
   \Theta:\text{H}_0\text{Cob}(M)\to\text{H}_0\text{Cob}(\text{Fuk}(M)).
\end{equation}
In particular, $\theta$ could be understood as the induced morphism of algebras over multicategories. A much stronger question is thus: when is $\Theta$ an equivalence of multicategories? A weaker question is whether the multimorphisms in the image of $\Theta$ generate $\text{H}_0\text{Cob}(\text{Fuk}(M))$. Here originates the name Lagrangian quantum field theory. Note that Biran and Cornea also previously used this name for similar reasons (see \cite{BC2}).
\begin{rmk}
Note that there is some additional work to be done here. Our construction was in the strictly unital setting, whereas Fukaya categories need not be strictly unital. However, it is known that any h-unital $A_\infty$-category over a field is quasi-equivalent to a strictly unital one; see \cite{Se}. One would expect the same construction to go through directly in the homologically unital setting by keeping track of the additional boundary terms needed.
\end{rmk}

\subsection{Quantitative LQFTs}
The work of Biran, Cornea and Zhang \cite{BCZ1,BCZ2} on triangulated persistence categories (TPCs) can be seen as an algebraic interpretation of geometric sizes that originate in the study of Lagrangian cobordisms. The work of Biran, Cornea and Shelukhin \cite{BCS} defining and studying ‘Lagrangian shadows’ gives a quantitative invariant associated to a cobordism. This can thus be interpreted as a size associated to the iterated cone decomposition it induces in the derived category. TPCs are a generalisation of the structure one obtains from this idea.

One can ask how to extend the theory of this paper to the quantitative setting, i.e., to filtered $A_\infty$-categories. This seems somewhat straightforward. One can define filtered analogues of the categories $Q_n$, which carry a size and a specific filtration, i.e., a category we can denote $Q_n^{\mathbf r}$, with $\mathbf r$ an object that keeps track of both this size ($r\in\R$) and the chosen filtration. Functors $\ff:Q_n^{\mathbf r}\to\cc$, with $\cc$ filtered, thus generalise functors $\ff:Q_n\to\cc$ and keep track of the notion of size. This again seems to fit more closely with the structure that arises from Lagrangian shadows.

Along with this, one can now consider not just weights associated to decompositions in triangulated categories, but general filtered $A_\infty$-categories. We have significant work in progress in this direction.

\subsection{Algebras over $\text{H}_0\text{Cob}(\cc)$ and iterated cones}
We showed that one can interpret the Grothendieck group as an algebra over the multicategory $\text{H}_0\text{Cob}(\cc)$. The obvious next question is: what other such algebras exist? Such an algebra should consist of an invariant associated to each object of $\cc$ that ‘behaves’ with respect to taking iterated cones.

In the symplectic setting, one would expect invariants of Lagrangians that ‘behave’ with respect to Lagrangian cobordisms to admit an interpretation as algebras over $\text{H}_0\text{Cob}(\text{Fuk}(M))$. The operations of such an algebra would describe how these invariants are related by a cobordism, with composition reflecting the gluing of cobordisms.

It's important to note the following: in general, there need not be a canonical map
\begin{equation}
   \rho:\Delta\text{H}_0\cc(\Delta,\Delta')
   \to\text{H}_0\text{Fun}_{A_\infty}(Q_2,\cc)(\ff,\ff'),
\end{equation}
where $\Delta\text{H}_0\cc$ denotes the category of exact triangles in $\text{H}_0\cc$ (assumed to be triangulated), and $\ff$, $\ff'$ represent the exact triangles $\Delta$, $\Delta'$, respectively. The issue is that there is no canonical way to lift morphisms of exact triangles to natural transformations of the functors that represent them. One can see this by noting that a morphism of exact triangles is a tuple of morphisms between the respective objects, in order, that commute with the morphisms defining the triangles. This commutativity is at the homological level. There is no canonical way to lift these morphisms to cycles that extend to (closed degree-zero) natural transformations. This issue extends to the case of general functors $\ff:Q_n\to\cc$ and their respective images in homology.

Thus it is advantageous to work directly with $\text{H}_0\text{Cob}(\cc)$, rather than attempting to build a similar operadic structure directly from a triangulated category, where one would be required to make many additional choices. It seems that this is the natural setting to build invariants of objects that respect cone decompositions, and to view them in terms of algebras over $\text{H}_0\text{Cob}(\cc)$.

\subsection{Constructability and $n$-angulated categories}
One can consider a variant of constructability for individual morphisms. Say $f\in\cc(A,B)$ is a closed degree-zero morphism. We can call $f$ $n$-constructable if there exists a functor $\ff:Q_m\to\cc$ with $1\leq m\leq n$ and $\ff_1(a_1)=f$. With this, if $\cc$ is pre-triangulated, then we automatically obtain that every $f$ is $2$-constructable. One will notice that we need to further restrict the notion of $n$-constructability for morphisms by considering the following example. Assume $\ff:Q_2\to\cc$ is well-defined. We can consider the gluing $\gg:=\ff\#_{21}(R\ff):Q_3\to\cc$, where $R\ff$ is the rotation of $\ff$ such that $(R\ff)_0(0)=\ff_0(1)$. Write
\[
A=\ff_0(0),\qquad B=\ff_0(1),\qquad C=\ff_0(2),
\qquad f=\ff_1(a_1),
\]
so that
\[
(R\ff)_0(0)=B,\qquad
(R\ff)_0(1)=C,\qquad
(R\ff)_0(2)=A[1].
\]
Thus one will find that in homology we have a diagram
\[
\begin{tikzcd}[row sep=large,column sep=large]
    A \arrow[r,"{[f]}"]
      & B \arrow[d,"0"]\\
    B \arrow[u,"0"]
      & A[1] \arrow[l,"{[\widetilde f]}"]
\end{tikzcd}
\]
where $\widetilde f:A[1]\to B$ is $f$ regarded as a morphism of degree $-1$. In general, assuming $\cc$ has shifts, one can directly construct a functor $\gg:Q_3\to\cc$ of the above form for any closed degree-zero morphism $f$, without requiring $f$ to fit into a functor from $Q_2$. Thus any such $f$ is $3$-constructable, with $\gg_1(a_1)=f$. This leads to wanting a refined notion of constructability. Perhaps we only allow functors $\ff:Q_n\to\cc$ with no $\ff_1(a_i)\simeq0$, or impose some other notion of non-degeneracy.

There is a possible connection with the theory of
$n$-angulated categories introduced by
Geiss, Keller and Oppermann \cite{GKO}. More precisely,
$Q_n$ has $n+1$ objects, and, after the appropriate object
shifts, its first-component arrows give a sequence of the form
\[
    X_0\longrightarrow X_1\longrightarrow\cdots
    \longrightarrow X_n\longrightarrow X_0[n-1].
\]
For $n\geq2$, in the terminology of \cite{GKO}, this has the form of an $(n+1)$-angle with suspension $[n-1]$. An $(n+1)$-angulated category also allows $(n+1)$-angles with some of their arrows equal to zero, but perhaps one can realise the non-degenerate functors $\ff:Q_n\to\cc$ as a well-behaved class of $(n+1)$-angles.

\subsection{Derived equivalence vs.\ equivalence}
As made clear in the previous section, we can use the filtration of $K^Q_0(\cc)$ to distinguish $\cc$ from its triangulated completion. Thus, in general, we can use this to help distinguish two categories that are derived equivalent. That is, if $\cc$ and $\cc'$ satisfy $D\cc\simeq D\cc'$ but we know, for example, that $\kappa(\cc)<\kappa(\cc')$ (here we require $\kappa(\cc),\kappa(\cc')>0$), then $\cc$ and $\cc'$ cannot be equivalent, or even quasi-equivalent. Though we note that $\kappa(\cc)=\kappa(\cc')$ does not imply that they are quasi-equivalent.

This suggests a stronger question. To what
extent does the multicategory $\mathrm H_0\mathrm{Cob}(\cc)$,
or its categories of algebras, determine $\cc$ up to
quasi-equivalence? That is, assume there are equivalences
\[
 \operatorname{Alg}_{\mathrm H_0\mathrm{Cob}(\cc)}(\dd)
 \simeq
 \operatorname{Alg}_{\mathrm H_0\mathrm{Cob}(\cc')}(\dd),
\]
with some additional compatibility requirements, where the above denotes the category of algebras in $\dd$. What does this imply about the relation between $\cc$ and $\cc'$?

\subsection{A $K$-theory spectrum associated to a multicategory}
Recent work on scissors congruence, including \cite{BGMMZ},
studies decompositions of geometric objects through categories
with covering families and their associated $K$-theory spectra.
This suggests a connection with our multicategory: its
multimorphisms describe how objects can be assembled through
iterated cones, and composition records successive assemblies.
A related construction using based multicategories appears in the work of Elmendorf and Mandell \cite{EM}.

It would be interesting to understand how this construction can be applied to a general $\text{H}_0\text{Cob}(\cc)$. It seems plausible that one will be able to construct a connective spectrum $\mathcal S(\cc)$ using their method, and moreover perhaps we can view $K^Q_0(\cc)$ as $\pi_0\mathcal S(\cc)$. This would also allow us to perhaps recover all $K$-groups associated to $\text{Tw}\cc$ in terms of a spectrum over $\text{H}_0\text{Cob}(\cc)$. A stronger question is whether this identification extends to a comparison with the $K$-theory spectrum of $\mathrm{Tw}\cc$:
\[
    \mathcal S(\cc)\longrightarrow\mathbf K(\mathrm{Tw}\cc).
\]
An equivalence would thus recover the higher groups
$K_i(\mathrm{Tw}\cc)$ for $i\geq0$.
The construction in \cite{BGMMZ} uses categories with covering families, whose covering families form a symmetric multicategory. One would therefore need to consider how to associate such a category with covering families to our non-symmetric multicategory $\text{H}_0\text{Cob}(\cc)$.

\clearpage

\phantomsection

\end{document}